\documentclass{cocv}
\usepackage{graphicx}
\usepackage{subfigure}
\begin{document}
%%-----------------------------
%%      the top matter
%%-----------------------------
\title{Social Optima in Leader-Follower Mean Field Linear Quadratic Control}\thanks{The work of Jianhui Huang was supported by by RGC Grants PolyU 153005/14P, 153275/16P; The work of Bing-Chang Wang was supported by the National Natural Science Foundation of China under Grant 61773241}%\thanks{...}% At most 5 thanks
\author{Jianhui Huang}\address{Department of Applied Mathematics, The Hong Kong Polytechnic University, Hong Kong, China;\\ \email{james.huang@polyu.edu.hk;\ tinghan.xie@connect.polyu.hk}}
\author{Bing-Chang Wang}\address{School of Control Science and Engineering, Shandong University, Jinan, China; \email{bcwang@sdu.edu.cn}}
\author{Tinghan Xie$\;^{\dagger,}$}\footnote[2]{\scriptsize Corresponding author} \sameaddress{1}
%
%\date{...}
%
\begin{abstract} This paper investigates a linear quadratic mean field leader-follower team problem, where the model involves one leader and a large number of weakly-coupled interactive followers. The leader and the followers cooperate to optimize the social cost. Specifically, for any strategy provided first by the leader, the followers would like to choose a strategy to minimize social cost functional.
Using variational analysis and person-by-person optimality, we construct two auxiliary control problems. By solving sequentially the auxiliary control problems with consistent mean field approximations, we can obtain a set of decentralized social optimality strategy with help of a class of forward-backward consistency systems. The relevant Stackelberg equilibrium is further proved under some proper conditions. \end{abstract}
%
%\begin{resume} ... \end{resume}
%
\maketitle
%%-----------------------------
%%      your text
%%-----------------------------
\section*{Introduction}
Mean field games have been studied by researchers from various aspects \cite{bfy2013,c2014,cdl2016,ll2007}. They involve a large number of population and the interaction between each individual is negligible. The mean field game approach has been applied in many fields such as finance \cite{fc2016}, economics \cite{wbvr2008}, information technology \cite{hcm2003a}, engineering \cite{cptd2012,km2014} and medicine \cite{cd2004}. %especially for specially for combining with linear quadratic (LQ) control problems.
 %mean field linear quadratic (LQ) control problem which has been introduced by many literatures. For example, $\varepsilon$-Nash equilibrium strategies in LQ mean field game are constructed in by using the Nash certainty equivalence method.
 The mean field linear quadratic (LQ) control problem is a special class of control problems, which can model many problems in applications and its solution exhibits elegant properties.
 For more work about the problem, readers can refer to \cite{bp2014,cd2013,jm2017,hl2018,hcm2012,lz2008,tzb2014,bj2012a,bj2017b,yz1999}.

For the model with one major player and $N$ minor players, the state of the major player has a significant influence on state equations and cost functionals of other minor individuals, which can be considered as a strong effect of the major player on minor ones. Mean field games (control) with major and minor players have been discussed in some literature, such as \cite{hn2016,nh2012} for LQ problems, \cite{bcy2016a,nc2013} for nonlinear problems, \cite{cz2014,hww2016} for probabilistic approaches and \cite{cw2016} for finite-state problems. In such types of problem, there is no hierarchical structure of decision making between the major player and the minor players.

In contrast to the model discussed above, leader-follower (Stackelberg) problems contain at least two hierarchies of players. One hierarchy of the players is defined as the leaders with a major position and another hierarchy of the players is defined as the followers with a minor position. The leader has priority to announce a strategy first and then the followers seek their strategies to minimize their cost functionals with response to the given strategy of leader. According to the followers' optimal points, the leader will choose his optimal strategy to minimize its cost functional. %The leader has priority to first give a strategy and then the follower chooses an optimal strategy by the given strategy of leader. Considering the follower will select an optimal strategy, the leader will take an optimal strategy to maximize its payoff.
Leader-follower problem also has been widely investigated. For example, the two-person leader-follower problem combines with stochastic LQ differential game had been studied by Yong in \cite{y02} and the problem of one leader and $N$ followers who play a noncooperative game under LQ stochastic differential game had been studied by Moon and Basar in \cite{mb2018}. For further literature related to Stackelberg games, readers can refer to \cite{bbs2010,bo1999,ljz2019,mzzgb2013,swx2016,sc1973,bj2017a} .
%Meanwhile Stackelberg game problem have been applied in an extensive range. For instance, in a communication network where the followers at lower level play a Nash game and the leader at upper level aims to choose an optimal price to make the revenue maximum [].

Different from noncooperative games, the social optimization (team optimization) problem is a joint decision problem which all the players have the same goal and work cooperatively. %The players can acquired information through observing the environment or communicating with other players.
The aim of each player is to select an optimal strategy and maximize the total payoff. %Team optimization models can be applied in many field, such as science, engineering, economic and finance, management science, and operations research. For example, traffic system of large cities or computer networks in geographical areas. The first person who began to consider team optimization problem was J. Marschak \cite{jm1955} and the method he used is based on game theory and decision theory. After that, the relevant work are limited to static team \cite{ra1962}.
Team optimization problem has been studied for many years. Marschak \cite{jm1955} first considered team optimization based on game theory. Ho and Chu studied team decision theory in optimal control problems \cite{hch1972}. Groves did the research of viewing the incentive problem as a team problem which the information for decisions is incomplete \cite{tg1973}. The team theory and person-by-person optimization with binary decision was investigated by Bauso and Pesenti \cite{dr2012} and the team problems under stochastic information structure with suboptimal solutions was studied in \cite{gmm2012}.

In this paper, we investigate social optimality of the leader-follower mean field LQ control problem. %Our work has the following features.
Our model contains one leader and $N$ followers. The leader's  state appears in both state equation and cost functional of each follower. It shows that the dynamics and cost functionals of the $N$ followers are directly influenced by the behavior of the leader. Unlike the model in \cite{hn2016,mb2018}, our model has a population state average term in all state equations and cost functionals. This implies that such state dynamics and cost functionals are highly interactive and coupled. %Because of existence of $x^{(N)}$, the dynamics and the cost functionals are highly coupled. %Especially, the cost functional of each follower depends on $u=\{u_1, u_2, \cdots, u_N\}$ which presents such weak-coupling.
%Thus, we can use mean field approximation to investigate the followers' problem. %Numerical computation of relevant dynamic optimization is difficult due to its dimensionality and complexity.
In reality, it is almost impossible for one player to obtain all the information of other players. Therefore, decentralized control which is based on the individual information set will be used instead of centralized control which is based on full information set and the information structure of each agent is different.

Compared with previous works, this paper mainly makes the following contributions: %of compared with previous work are summarized as follows
\begin{itemize}
%\item The second situation contains $N_1$ leaders and $N_2$ followers. The same as the first situation, both the leader and followers in the second situation have their own dynamics and cost functionals. And the form of each follower's dynamic contains a state term of leader in the drift part. However, the state average term in each dynamic of the leaders is no longer the state average of followers. It is changed to the state average of the leaders. The cost functional of both leaders and followers contain both the state average of the leaders and the state average of the followers which the cost functionals are highly coupled. And the second situation is much more complicated than the first situation.
\item A social optimum problem is studied for mean field models with hierarchical structure. Unlike the problem in \cite{mb2018} where the leader and followers play a noncooperative game and try to minimize their own individual cost functional, all individuals in our models aim to minimize the social cost functional which equals the summation of cost functionals of all players. The $N$ followers are coupled by the population state average term. Since the cost functional presents individual performance in the game problems, the order of magnitude of the perturbation is $\frac{1}{N}$ which can be ignored. The population state average term may be approximated by a stochastic process directly (see \cite{jm2017}). However, in the team problems, the order of magnitude of the perturbation cannot be ignored after summing up all the cost functionals, which makes the problem very complicated.
    To overcome such difficulties, we approximate some terms as $N$ goes to infinity and use a duality procedure combined with auxiliary equations to transform the variation of the social cost functional into a standard LQ control form. Then, we construct an auxiliary control problem and a forward-backward consistency system which contains four equations to help us obtain the decentralized form of the optimal controls for the $N$ followers.
\item The decentralized controls of the leader-follower problem are obtained and the solvability of a high-dimensional consistency condition system (CC system) is discussed. Since the leader's state equation and cost functional are fully coupled with the followers' state equations and cost functionals, it is more difficult to solve the leader's problem. Except constructing auxiliary problem by mean field approximation as in the former part, we need to construct six auxiliary equations and use duality relations to obtain the decentralized form of the optimal control for the leader. %To solve the all the indeterminate coefficients,
    Unlike the problem for $N$ followers, the final CC system of the leader's problem contains ten equations which becomes a high-dimensional problem. %This is more complex than the situations in \cite{hn2016} and \cite{mb2018} which only have six equations.
    To solve such equations directly is very difficult since they are fully coupled and have high-dimensional characteristics. We transform the high-dimensional CC system to a simple form of linear forward-backward stochastic differential equation (FBSDE; see \cite{my1999,y1999}) and discuss the solvability of the FBSDE through Ricatti equation method.
\item The decentralized strategies of leader-follower problem are proved to be Stackelberg equilibrium by perturbation analysis. Different from \cite{hcm2012,hn2016,mb2018,wh2017}, we discuss the Stackelberg equilibrium for the team optimization problem. %which compares the difference of the social cost functional between the centralized control and the decentralized control. %Before proving the asymptotic optimality, we first need to %prove the errors of the mean field approximations are equal to $O(\frac{1}{\sqrt{N}})$.for the mean field approximation.
    %Unlike \cite{hcm2012} and \cite{wh2017}, the proof of asymptotic social optimality may be not easy to tackle for leader-follower problem.
    First, we need to prove the decentralized strategies for the followers have asymptotic social optimality. Because of the Stackelberg problem contains two hierarchies, we consider two coupled cost functionals (the leader and the followers) when using the standard method (see \cite{hcm2012}). We prove the asymptotic optimality by decoupling them with two duality procedures and some arguments in error estimates. Second, we need to prove the decentralized strategies for the leader-follower problem is Stackelberg equilibrium. Also, some error estimates are very hard to be given directly since they are fully coupling. We decompose them by applying Ricatti equation method and then estimate them in proper order. %For the first part, it can be proved by applying the continuity dependence of social cost functional and some arguments in error estimates. %Combining these two parts, our result follows. %the difference of the social cost functional under between decentralized and centralized strategy is not great than $O(\frac{1}{\sqrt{N}})$.
\end{itemize}

In the real world, our model can be used to describe some examples. For an automatic machine, the major part first gives an information to the system and the minor parts will adjust their parameters automatically such that the whole system keeps in the best state. For the economic environment, the small companies may hesitate to make decisions and often follow a monopoly company when facing to the volatility of the market. A monopoly company announces a decision first. Once the small companies try to make decisions according to their own situations, the monopoly company adjusts its decision such that the sum of the social wealth can be maximized. Moreover, the relationship between the employer and the employees or the federal government and the government in each state also can be described by our model.

The paper is organized as follows. The problem is formulated in Section 2. %the social optima LQG mean field control involving a leader and $N$ followers.
In Section 3, we solve the optimal controls for followers based on person-by-person optimality  %and give the details of using mean field approximation method and social variational method to obtain the variations of personal state equations and social cost functional. We design an auxiliary optimal control problem with mean field approximation method
and obtain the CC system of the follower's problem.
In Section 4, we seek the social optimal solution of the leader's problem and give the CC system of the leader's problem. Then the CC system is transformed to a simple form of linear FBSDE in Section 5 and its wellposeness is discussed. %and the solvability of relevant Ricatti equations.
In Section 6, we give the details of proving the Stackelberg equilibrium. Also, some prior lemmas will be introduced and proved.
In Section 7, a numerical example is provided to simulate the efficiency of decentralized control. Section 8 is the conclusion of the paper.

%\subsection{Notation}
\emph{Notation}: Throughout this paper, $\mathbb{R}^{n\times m}$ and $\mathbb{S}^n$ denote the set of all $(n\times m)$ real matrices and the set of all $(n\times n)$ symmetric matrices, respectively. $\|\cdot\|$ is the standard Euclidean norm and $\langle\cdot,\cdot\rangle$ is the standard Euclidean inner product. For given symmetric matrix $S\geq 0$, the quadratic form $x^TSx$ may be defined as $\|x\|_S^2$, where $x^T$ is the transpose of $x$. $C^1([0,T];\mathbb{R}^{n\times m})$ be the space of all $\mathbb{R}^{n\times m}$-valued continuously differentiable functions on $[0, T]$. For notation $o(1)$, $\lim_{n\rightarrow \infty}o(1)=0$. %We denote $\|u\|_{L^2}=:\mathbb{E}[\int_{0}^{T}\|u\|^2dt]$, where $u$ is a vector-valued process.
By \cite{y1999}, for sake of notation simplicity, we will use $K$ to denote a generic constant in following discussion. The value of $K$ may be different at different places and it only depends on the coefficients and initial values.

\section{Problem formulation}

Let $(\Omega,\mathcal{F}, \mathbb{P})$ be a complete probability space which contains all $\mathbb{P}$-null sets in $\mathcal{F}$. $\xi_i\in \mathbb{R}^n$ are the values of the initial states and $W_i(\cdot)$ are $d$-dimensional standard Brownian motions, where $i=0,1,\cdots, N$. $\xi_i$ and $W_i(\cdot)$ are defined on $(\Omega,\mathcal{F}, \mathbb{P})$. Consider a large-population system which contains one leader and $N$ followers. The state processes of the leader and the $i$th follower, $i=1,2,\cdots,N,$ are modeled by the following linear stochastic differential equations (SDE) on a finite time horizon $[0,T]$:
\begin{equation}
\left\{
\begin{aligned}
&dx_0(t)=[A_0(t)x_0(t)+B_0(t)u_0(t)+C_0(t)x^{(N)}(t)]dt+D_0(t)dW_0(t),\quad x_0(0)=\xi_0,\\
&dx_i(t)=[A(t)x_i(t)+B(t)u_i(t)+C(t)x^{(N)}(t)+F(t)x_0(t)]dt+D(t)dW_i(t), \quad x_i(0)=\xi_{i},\\
%& \quad \quad
\end{aligned}
\right.
\end{equation}\\
where $x^{(N)}(t):=\frac{1}{N}\sum_{i=1}^{N}x_i(t)$ is the state average of the followers. Let $\sigma$-algebra $\mathcal{F}_t^i=\sigma(W_i(s), 0\leq s\leq t)$ and
$\mathcal{G}_t^i=\mathcal{F}_t^i\bigvee\sigma\{\xi_i, \xi_0, W_0(s), 0\leq s\leq t\}$, where $0\leq i\leq N$.
$\mathcal{F}_t=\sigma(W_i(s), 0\leq s\leq t, 0\leq i\leq N)$
%$\mathcal{G}_t^0=\mathcal{F}_t^0\bigvee\sigma\{\xi_0\}$,
%$\mathcal{H}_t^i=\mathcal{F}_t^i\bigvee\sigma\{x_i(s),0\leq s\leq t\}$ for $0\leq i\leq N$
and
$\mathcal{G}_t=\mathcal{F}_t\bigvee\sigma\{\xi_i,0\leq i\leq N\}$.
%$\mathcal{H}_t=\mathcal{F}_t\bigvee\sigma\{x_i(s),0\leq s\leq t,0\leq i\leq N\}$.
$\mathbb{F}^i=\{\mathcal{F}_t^i\}_{0\leq t\leq T}$ %$0\leq i\leq N$
is the natural filtration generated by $W_i(\cdot)$ and %$\mathbb{G}^0=\{\mathcal{G}_t^0\}_{0\leq t\leq T}$,
$\mathbb{G}^i=\{\mathcal{G}_t^i\}_{0\leq t\leq T}$, where $0\leq i\leq N$.
%$\mathbb{H}^i=\{\mathcal{H}_t^i\}_{0\leq t\leq T}$ for $0\leq i\leq N$.
Correspondingly, we denote
% and,
$\mathbb{F}=\{\mathcal{F}_t\}_{0\leq t\leq T}$, $\mathbb{G}=\{\mathcal{G}_t\}_{0\leq t\leq T}$.
%$\mathbb{H}=\{\mathcal{H}_t\}_{0\leq t\leq T}$.
Next we introduce the following spaces:
\begin{equation*}
\begin{aligned}
L^{\infty}(0,T;\mathbb{R}^{n\times m})=&\big\{\varphi : [0,T]\rightarrow \mathbb{R}^{n\times m}\big| \ \varphi(\cdot) \text{ is bounded and measurable}\big\},\\
L_{\mathbb{F}}^2(\Omega;\mathbb{R}^m)=&\big\{\xi : \Omega\rightarrow \mathbb{R}^m\big| \ \xi \text{ is }\mathbb{F}\text{-measurable, }\mathbb{E}\|\xi\|^2<\infty\big\},\\
%\end{aligned}
%\end{equation*}
%\begin{equation*}
%\begin{aligned}
L_{\mathbb{F}}^2(0,T;\mathbb{R}^m)=&\big\{x : [0,T]\times\Omega\rightarrow \mathbb{R}^m\big| \ x(\cdot) \text{ is }\mathbb{F}\text{-progressively}\\
&\hspace{0.5cm}\text{measurable, }\|x(t)\|_{L^2}^2:=\mathbb{E}\int_{0}^{T}\|x(t)\|^2dt<\infty\big\},\\
L_{\mathbb{F}}^2(\Omega;C([0,T];\mathbb{R}^m))=&\big\{x : [0,T]\times\Omega\rightarrow \mathbb{R}^m\big| \ x(\cdot) \text{ is }\mathbb{F}\text{-progressively} \\
&\hspace{0.5cm}\text{measurable, continuous, } \mathbb{E}\sup_{t\in[0,T]}\|x(t)\|^2<\infty\big\},\\
\mathcal{M}[0,T]:= L_{\mathbb{F}}^2(\Omega;&C([0,T];\mathbb{R}^{n}))\times L_{\mathbb{F}}^2(\Omega;C([0,T];\mathbb{R}^{m}))\times L_{\mathbb{F}}^2(0,T;\mathbb{R}^{m\times d}).
\end{aligned}
\end{equation*}
The set of admissible controls for the leader is defined as follows:
\begin{equation*}
\mathcal{U}_0=\left\{u_0|u_0(t)\in L_{\mathbb{G}^0}^2(0,T;\mathbb{R}^m)\right\},
\end{equation*}
and the set of admissible controls for the $i$th follower is defined as follows:
$$\mathcal{U}_i=\left\{u_i|u_i(t)\in L_{\mathbb{G}^i}^2(0,T;\mathbb{R}^m)\right\},\quad 1\leq i \leq N.$$
These are the decentralized control sets and we let
$\mathcal{U}=\mathcal{U}_1\times\mathcal{U}_2 \times\cdots\times\mathcal{U}_N$.
For comparison, the centralized control set is given by
\begin{equation*}
  \mathcal{U}_c = \Big\{(u_0,u_1,\cdots,u_N)|u_i(t)\in L_{\mathbb{G}}^2(0,T;\mathbb{R}^m), 0\leq i \leq N\Big\}.
\end{equation*}

Now we introduce the cost functionals of the leader and the $i$th follower, $1\leq i\leq N$. For the leader, the cost functional is defined as follows:
\begin{equation}
\begin{aligned}
\mathcal{J}_0(u_0(\cdot);u(\cdot))=&\mathbb{E}\bigg\{\int_{0}^{T}\big[\|x_0(t)- \Theta_0(t)x^{(N)}(t)-\eta_0(t)\|_{Q_0(t)}^{2}+ \|u_0(t)\|_{R_0(t)}^{2}\big]dt\\
&+ \|x_0(T)- \hat{\Theta}_0x^{(N)}(T)-\hat{\eta}_0\|_{G_0}^{2}\bigg\},
\end{aligned}
\end{equation}
where $u(\cdot)=(u_1(t),\cdots,u_N(t))\in\mathcal{U}_c$. $Q_0(\cdot)$, $R_0(\cdot)$ and $G_0(\cdot)$ are weight matrices. $Q_0(\cdot)$ and $\Theta_0(\cdot)$ represent the coupling between the leader and the population state average. This implies that the states of the followers can influence the cost functional of the leader. For the $i$th follower, the individual cost functional is defined as follows:
\begin{equation}
\begin{aligned}
\mathcal{J}_i(u_0(\cdot);u(\cdot))=&\mathbb{E}\bigg\{\int_{0}^{T}\big[\|x_i(t)- \Theta(t) x^{(N)}(t)-\Theta_1(t)x_0(t)-\eta(t)\|_{Q(t)}^{2} +\|u_i(t)\|_{R(t)}^{2}\big]dt\\
&+\|x_i(T)- \hat{\Theta} x^{(N)}(T)-\hat{\Theta}_1x_0(T)-\hat{\eta}\|_{G}^{2}\bigg\},
\end{aligned}
\end{equation}
where $Q(\cdot)$, $R(\cdot)$ and $G(\cdot)$ are weight matrices. $Q(\cdot)$, $\Theta(\cdot)$ and $\Theta_1(\cdot)$ represent the coupling between the $i$th follower, the population state average and the leader. This implies that the cost functional of the $i$th follower will be affected by the behavior of both the leader and the other followers. All the individuals in the system, including the leader and followers, aim to minimize the social cost functional, which is denoted by
\begin{equation*}
\mathcal{J}_{soc}^{(N)}(u_0(\cdot);u(\cdot))=\alpha N\mathcal{J}_0(u_0(\cdot);u(\cdot)) +\sum_{i=1}^{N}\mathcal{J}_i(u_0(\cdot);u(\cdot)),\quad \alpha>0.
\end{equation*}
%where $\alpha\in \mathbb{R}$ and $\alpha>1$.\\
Similar to \cite{hn2016} and \cite{nh2012}, we have a scaling factor $\alpha N$ before $\mathcal{J}_0(u_0(\cdot);u(\cdot))$ such that $\mathcal{J}_0(u_0(\cdot);u(\cdot))$ and $\mathcal{J}_i(u_0(\cdot);u(\cdot))$ have the same order of magnitude. Otherwise, if $\alpha N=1$, then the performance of the leader will be insensitive when $N$ becomes larger. %Thus, in our model, $\alpha$ should not be too small.
Now we introduce our assumptions.\\
%\begin{assumption}
\text{\rm(\textbf{A1})} The coefficients of \text{\rm(1)}, \text{\rm(2)} and \text{\rm(3)} satisfy
\begin{equation*}\left\{
\begin{aligned}
&A_0(\cdot), C_0(\cdot), A(\cdot), C(\cdot), F(\cdot)\in L^{\infty}(0,T;\mathbb{R}^{n\times n}),\\
&B_0(\cdot), B(\cdot)\in L^{\infty}(0,T;\mathbb{R}^{n\times m}),\quad D_0(\cdot), D(\cdot)\in L^{\infty}(0,T;\mathbb{R}^{n\times d}).\\
\end{aligned}\right.
\end{equation*}
\begin{equation*}\left\{
\begin{aligned}
&Q_0(\cdot), Q(\cdot)\in L^{\infty}(0,T;\mathbb{S}^n),\quad R_0(\cdot), R(\cdot)\in L^{\infty}(0,T;\mathbb{S}^m),\\
&\Theta_0(\cdot), \Theta_1(\cdot), \Theta(\cdot)\in L^{\infty}(0,T;\mathbb{R}^{n\times n}),\quad \eta_0(\cdot), \eta(\cdot)\in L^{2}(0,T;\mathbb{R}^{n}),\\
&\hat{\Theta}_0, \hat{\Theta}_1, \hat{\Theta}\in \mathbb{R}^{n\times n},\quad G_0, G\in\mathbb{S}^n,\quad \hat{\eta}_0, \hat{\eta}\in \mathbb{R}^{n}.
\end{aligned}\right.
\end{equation*}
%\end{assumption}
%\begin{assumption}
\text{\rm(\textbf{A2})} $x_0(0)$ and $W_0(\cdot)$ are mutually independent. $\{x_i(0),1\leq i\leq N\}$ and $\{W_i(\cdot),1\leq i\leq N\}$ are independent of each other. %$\mathbb{E}x_0(0)=\hat{\xi}_0$ and
$\mathbb{E}x_i(0)=\hat{\xi}$, $1\leq i\leq N$. For some constant $K$, which is independent of $N$, such that $\sup_{1\leq i \leq N}\mathbb{E}\|x_i(0)\|^2\leq K$. Furthermore, $x_0(0)$, $W_0(\cdot)$ and $\{x_i(0),1\leq i\leq N\}$, $\{W_i(t),1\leq i\leq N\}$ are independent of each other.\\
%\end{assumption}
%\begin{assumption}
\text{\rm(\textbf{A3})} $Q_0(\cdot)\geq 0$, $R_0(\cdot)>\delta I$, $G_0\geq 0$ and $Q(\cdot)\geq 0$, $R(\cdot)>\delta I$, $G\geq 0$, for some $\delta>0$.
%\end{assumption}

From now on, we may suppress the notation of time $t$ if necessary. We introduce our leader-follower problem:
\begin{prblm}
Under \text{\rm(\textbf{A1})-(\textbf{A3})}, for any $u_0\in \mathcal{U}_0$, to find a mapping $\mathcal{M}$: $\mathcal{U}_0\rightarrow \mathcal{U}$, and a control $\bar{u}_0\in \mathcal{U}_0$ such that
\begin{equation*}
\left\{
\begin{aligned}
&\mathcal{J}_{soc}^{(N)}(u_0;\mathcal{M}(u_0)) =\inf_{u\in\mathcal{U}_c}\mathcal{J}_{soc}^{(N)}(u_0;u),\\
&\mathcal{J}_{soc}^{(N)}(\bar{u}_0;\mathcal{M}(\bar{u}_0))= \inf_{u_0\in\mathcal{U}_0}\mathcal{J}_{soc}^{(N)}(u_0;\mathcal{M}(u_0)).
\end{aligned}
\right.
\end{equation*}
\end{prblm}
Note that the $\mathcal{M}$ here is a mapping, which is different from the notation $\mathcal{M}[0,T]$ we just introduced.

\section{The mean field LQ control problem for the $N$ followers}

\subsection{Person-by-person optimality}

Fix %$(u_0, x_0)$,
$u_0\in \mathcal{U}_0$. The leader firstly announces his own open-loop strategy. Let $\bar{u}=\{\bar{u}_1, \bar{u}_2,\cdots, \bar{u}_N\}$ be the centralized optimal control of the followers and $\bar{x}=\{\bar{x}_1, \bar{x}_2,\cdots, \bar{x}_N\}$ be the corresponding states. Now we perturb $\bar{u}_i$ and fix other $\bar{u}_j$, where $j\neq i$. Then we denote
$\delta u_i=u_i-\bar{u}_i$, $\delta x_i=x_i-\bar{x}_i$, where $u_i$ is the control after perturbing and $x_i$ is its corresponding state. The Fr\'{e}chet differential $\delta \mathcal{J}_0(\delta u_i)=\mathcal{J}_0(u_0;u)-\mathcal{J}_0(u_0;\bar{u})+o(\|\delta u_i\|)$ and $\delta \mathcal{J}_i(\delta u_i)=\mathcal{J}_i(u_0;u)-\mathcal{J}_i(u_0;\bar{u})+o(\|\delta u_i\|)$, where $i=1,\cdots, N$. Therefore, the variations of the state equations for the leader, the $i$th follower and the $j$th follower, where $j\neq i$, are
\begin{equation*}
\left\{
\begin{aligned}
&d\delta x_0=(A_0\delta x_0+C_0\delta x^{(N)})dt, \quad \delta x_0(0)=0,\\
&d\delta x_i=(A\delta x_i+B\delta u_i+C\delta x^{(N)}+F\delta x_0)dt, \quad \delta x_i(0)=0,\\
&d\delta x_j=(A\delta x_j+C\delta x^{(N)}+F\delta x_0)dt, \quad \delta x_j(0)=0, \quad j\neq i,\\
\end{aligned}
\right.
\end{equation*}
and the variations of their corresponding cost functionals are
\begin{equation*}%\resizebox{\linewidth}{!}{$
\begin{aligned}
\frac{1}{2}\delta \mathcal{J}_0(\delta u_i)=&\mathbb{E}\bigg\{\int_{0}^{T}\langle Q_0(\bar{x}_0- \Theta_0\bar{x}^{(N)}-\eta_0),\delta x_0- \Theta_0\delta x^{(N)}\rangle dt\\
&+ \langle G_0(\bar{x}_0(T)- \hat{\Theta}_0\bar{x}^{(N)}(T)-\hat{\eta}_0),\delta x_0(T)- \hat{\Theta}_0\delta x^{(N)}(T)\rangle \bigg\},\\
\end{aligned}
\end{equation*}
%and
\begin{equation*}
\begin{aligned}
\frac{1}{2}\delta \mathcal{J}_i(\delta u_i)=&\mathbb{E}\bigg\{\int_{0}^{T}\langle Q(\bar{x}_i- \Theta \bar{x}^{(N)}-\Theta_1\bar{x}_0-\eta),\delta x_i- \Theta \delta x^{(N)}-\Theta_1\delta x_0\rangle+\langle R\bar{u}_i,\delta u_i\rangle dt\\
&+\langle G(\bar{x}_i(T) -\hat{\Theta}\bar{x}^{(N)}(T) -\hat{\Theta}_1\bar{x}_0(T)-\hat{\eta}),\delta x_i(T)- \hat{\Theta} \delta x^{(N)}(T)-\hat{\Theta}_1\delta x_0(T)\rangle \bigg\},\\
\end{aligned}
\end{equation*}
%If $j\neq i$, then the cost functional of the other followers is:
\begin{equation*}
\begin{aligned}
\frac{1}{2}\delta \mathcal{J}_j(\delta u_i)=&\mathbb{E}\bigg\{\int_{0}^{T}\langle Q(\bar{x}_j- \Theta \bar{x}^{(N)}-\Theta_1\bar{x}_0-\eta),\delta x_j- \Theta \delta x^{(N)}-\Theta_1\delta x_0\rangle dt\\
&+\langle G(\bar{x}_j(T)-\hat{\Theta}\bar{x}^{(N)}(T) -\hat{\Theta}_1\bar{x}_0(T)-\hat{\eta}),\delta x_j(T)- \hat{\Theta} \delta x^{(N)}(T)-\hat{\Theta}_1\delta x_0(T)\rangle \bigg\},
\end{aligned}%$}
\end{equation*}
respectively. Consequently, we have the variation of the social cost functional as:
\begin{equation}\label{variation_soc}%\resizebox{\linewidth}{!}{$
\begin{aligned}
&\frac{1}{2}\delta \mathcal{J}_{soc}^{(N)}(\delta u_i)=\frac{1}{2}\bigg[\alpha N\delta \mathcal{J}_0(\delta u_i)+\sum_{j\neq i}\delta \mathcal{J}_j(\delta u_i)+\delta \mathcal{J}_i(\delta u_i)\bigg]\\
&=\mathbb{E}\bigg\{\int_{0}^{T}\alpha N\langle Q_0(\bar{x}_0- \Theta_0\bar{x}^{(N)}-\eta_0),\delta x_0\rangle-\alpha N\langle \Theta_0^{T}Q_0(\bar{x}_0- \Theta_0\bar{x}^{(N)}-\eta_0),\delta x^{(N)}\rangle+\langle Q(\bar{x}_i- \Theta \bar{x}^{(N)}\\
&-\Theta_1\bar{x}_0-\eta),\delta x_i\rangle-\langle \Theta^{T}Q(\bar{x}_i-\Theta \bar{x}^{(N)}-\Theta_1\bar{x}_0-\eta),\delta x^{(N)}\rangle-\langle\Theta_1^{T}Q(\bar{x}_i-\Theta\bar{x}^{(N)} -\Theta_1\bar{x}_0-\eta),\delta x_0\rangle\\
&+\langle R\bar{u}_i,\delta u_i\rangle+\sum_{j\neq i}\langle Q(\bar{x}_j- \Theta \bar{x}^{(N)}-\Theta_1\bar{x}_0-\eta),\delta x_j\rangle-\sum_{j\neq i}\langle \Theta^{T}Q(\bar{x}_j- \Theta \bar{x}^{(N)}-\Theta_1\bar{x}_0-\eta), \delta x^{(N)}\rangle\\
&-\sum_{j\neq i}\langle\Theta_1^{T}Q(\bar{x}_j- \Theta \bar{x}^{(N)}-\Theta_1\bar{x}_0-\eta),\delta x_0\rangle dt+ \alpha N\langle G_0(\bar{x}_0(T)- \hat{\Theta}_0\bar{x}^{(N)}(T)-\hat{\eta}_0),\delta x_0(T)\rangle\\
&- \alpha N\langle \hat{\Theta}_0^{T}G_0(\bar{x}_0(T)- \hat{\Theta}_0\bar{x}^{(N)}(T) -\hat{\eta}_0),\delta x^{(N)}(T)\rangle+\langle G(\bar{x}_i(T) -\hat{\Theta}\bar{x}^{(N)}(T) -\hat{\Theta}_1\bar{x}_0(T)-\hat{\eta}),\delta x_i(T)\rangle\\
&-\langle \hat{\Theta}^{T}G(\bar{x}_i(T) -\hat{\Theta}\bar{x}^{(N)}(T) -\hat{\Theta}_1\bar{x}_0(T)-\hat{\eta}),\delta x^{(N)}(T)\rangle-\langle \hat{\Theta}_1^{T}G(\bar{x}_i(T) -\hat{\Theta}\bar{x}^{(N)}(T) -\hat{\Theta}_1\bar{x}_0(T)-\hat{\eta}),\\
&\delta x_0(T)\rangle+\sum_{j\neq i}\langle G(\bar{x}_j(T)-\hat{\Theta}\bar{x}^{(N)}(T) -\hat{\Theta}_1\bar{x}_0(T)-\hat{\eta}),\delta x_j(T)\rangle-\sum_{j\neq i}\langle \hat{\Theta}^{T}G(\bar{x}_j(T)-\hat{\Theta}\bar{x}^{(N)}(T) -\hat{\Theta}_1\bar{x}_0(T) \\
&-\hat{\eta}),\delta x^{(N)}(T)\rangle-\sum_{j\neq i}\langle \hat{\Theta}_1^{T}G(\bar{x}_j(T)-\hat{\Theta}\bar{x}^{(N)}(T) -\hat{\Theta}_1\bar{x}_0(T)-\hat{\eta}),\delta x_0(T)\rangle \bigg\}.
\end{aligned}%$}
\end{equation}
When $N\rightarrow \infty$, it follows that
\begin{equation*}%\resizebox{\linewidth}{!}{$
\begin{aligned}
&\frac{1}{2}\delta \mathcal{J}_{soc}^{(N)}(\delta u_i)
=\mathbb{E}\bigg\{\int_{0}^{T}\alpha\langle Q_0(\bar{x}_0- \Theta_0\bar{x}^{(N)}-\eta_0),N\delta x_0\rangle-\alpha\langle \Theta_0^{T}Q_0(\bar{x}_0- \Theta_0\bar{x}^{(N)}-\eta_0),N\delta x^{(N)}\rangle\\
&+\langle Q(\bar{x}_i- \Theta \bar{x}^{(N)}-\Theta_1\bar{x}_0-\eta),\delta x_i\rangle+\langle R\bar{u}_i,\delta u_i\rangle
%&-\langle \Theta^{T}Q(\bar{x}_i-\Theta \bar{x}^{(N)}-\Theta_1\bar{x}_0-\eta),\delta x^{(N)}\rangle\\
%&-\langle\Theta_1^{T}Q(\bar{x}_i-\Theta\bar{x}^{(N)} -\Theta_1\bar{x}_0-\eta),\delta x_0\rangle\\
+\langle \frac{1}{N}\sum_{j\neq i}Q(\bar{x}_j- \Theta \bar{x}^{(N)}-\Theta_1\bar{x}_0-\eta),N\delta x_j\rangle\\
&-\langle \frac{1}{N}\sum_{j\neq i}\Theta^{T}Q(\bar{x}_j- \Theta \bar{x}^{(N)}-\Theta_1\bar{x}_0-\eta), N\delta x^{(N)}\rangle-\langle \frac{1}{N}\sum_{j\neq i}\Theta_1^{T}Q(\bar{x}_j- \Theta \bar{x}^{(N)}-\Theta_1\bar{x}_0-\eta),\\
&N\delta x_0\rangle dt+ \alpha\langle G_0(\bar{x}_0(T)- \hat{\Theta}_0\bar{x}^{(N)}(T)-\hat{\eta}_0),N\delta x_0(T)\rangle- \alpha\langle \hat{\Theta}_0^{T}G_0(\bar{x}_0(T)- \hat{\Theta}_0\bar{x}^{(N)}(T)-\hat{\eta}_0),\\
& N\delta x^{(N)}(T)\rangle+\langle G(\bar{x}_i(T) -\hat{\Theta}\bar{x}^{(N)}(T) -\hat{\Theta}_1\bar{x}_0(T)-\hat{\eta}),\delta x_i(T)\rangle+\langle \frac{1}{N}\sum_{j\neq i}G(\bar{x}_j(T)-\hat{\Theta}\bar{x}^{(N)}(T)\\
%&-\langle \hat{\Theta}^{T}G(\bar{x}_i(T) -\hat{\Theta}\bar{x}^{(N)}(T) -\hat{\Theta}_1\bar{x}_0(T)-\hat{\eta}),\delta x^{(N)}(T)\rangle\\
%&-\langle \hat{\Theta}_1^{T}G(\bar{x}_i(T) -\hat{\Theta}\bar{x}^{(N)}(T) -\hat{\Theta}_1\bar{x}_0(T)-\hat{\eta}),\delta x_0(T)\rangle\\
& -\hat{\Theta}_1\bar{x}_0(T)-\hat{\eta}),N\delta x_j(T)\rangle
-\langle\frac{1}{N}\sum_{j\neq i} \hat{\Theta}^{T}G(\bar{x}_j(T)-\hat{\Theta}\bar{x}^{(N)}(T) -\hat{\Theta}_1\bar{x}_0(T)-\hat{\eta}),N\delta x^{(N)}(T)\rangle\\
&-\langle \frac{1}{N}\sum_{j\neq i} \hat{\Theta}_1^{T}G(\bar{x}_j(T)-\hat{\Theta}\bar{x}^{(N)}(T) -\hat{\Theta}_1\bar{x}_0(T)-\hat{\eta}),N\delta x_0(T)\rangle \bigg\}+o(1).\\
\end{aligned}%$}
\end{equation*}
Note that $\mathbb{E}\sup_{0\leq t\leq T}\|\delta x_0\|^2 = O(\frac{1}{N^2})$, $\mathbb{E}\sup_{0\leq t\leq T}\|\delta x^{(N)}\|^2 = O(\frac{1}{N^2})$ and $\langle \Theta^{T}Q(\bar{x}_i-\Theta \bar{x}^{(N)}-\Theta_1\bar{x}_0-\eta),\delta x^{(N)}\rangle+ \langle\Theta_1^{T}Q(\bar{x}_i-\Theta\bar{x}^{(N)} -\Theta_1\bar{x}_0-\eta),\delta x_0\rangle+\langle \hat{\Theta}^{T}G(\bar{x}_i(T) -\hat{\Theta}\bar{x}^{(N)}(T) -\hat{\Theta}_1\bar{x}_0(T)-\hat{\eta}),\delta x^{(N)}(T)\rangle+\langle \hat{\Theta}_1^{T}G(\bar{x}_i(T) -\hat{\Theta}\bar{x}^{(N)}(T) -\hat{\Theta}_1\bar{x}_0(T)-\hat{\eta}),\delta x_0(T)\rangle=o(1)$ (the rigorous proof will be shown in Section 5).
Let
\begin{equation}\label{x_dagger}
\left\{
\begin{aligned}
&\delta x_0^{\dagger}=\lim_{N\rightarrow+\infty}(N\delta x_0),\\
&\delta x^{\dagger}=\lim_{N\rightarrow+\infty}(N\delta x_j)=\lim_{N\rightarrow+\infty}(\sum_{j\neq i}\delta x_j),\ j\neq i.
\end{aligned}
\right.
\end{equation}
%where the convergence are in $L^2$.
Here $N\delta x_0$ converges to $\delta x_0^{\dagger}$ such that $\mathbb{E}\int_{0}^{T}\|N\delta x_0-\delta x_0^{\dagger}\|^2=O(\frac{1}{N^2})$. Similarly, $\sum_{j\neq i}\delta x_j$ and $N\delta x_j$ converge to $\delta x^{\dagger}$ (see Section 5 for the detailed proof). Then one can obtain
\begin{equation}\label{dx_dagger}
\left\{
\begin{aligned}
&d\delta x_0^{\dagger}=(A_0\delta x_0^{\dagger}+C_0\delta x_i+C_0\delta x^{\dagger})dt, \quad \quad \quad \delta x_0^{\dagger}(0)=0,\\
&d\delta x^{\dagger}=(A\delta x^{\dagger}+C\delta x_i+C\delta x^{\dagger}+F\delta x_0^{\dagger})dt, \quad \delta x^{\dagger}(0)=0.\\
%&d\delta x_i=(A\delta x_i+B\delta u_i+C\delta x^{(N)}+F\delta x_0)dt, \quad \delta x_i(0)=0.\\
\end{aligned}
\right.
\end{equation}
When $N\rightarrow\infty$, by mean field approximation, we use $\hat{x}$ to approximate $\bar{x}^{(N)}$. Note that $\hat{x}$ will be affected by $u_0$ which is given by the leader. Moreover, the influence of individual follower on $\hat{x}$ may be negligible. Hence, by straightforward computation, we simplified the social cost functional as follows:
\begin{equation}\label{lf_ctr_cost}%\resizebox{\linewidth}{!}{$
\begin{aligned}
\frac{1}{2}\delta \mathcal{J}_{soc}^{(N)}(\delta u_i)
&=\mathbb{E}\bigg\{\int_{0}^{T}\langle \alpha Q_0\Psi_1-\Theta_1^{T}Q\Psi_3,\delta x_0^{\dagger}\rangle+\langle Q\Psi_2^i-\Theta^{T}Q\Psi_3-\alpha\Theta_0^{T}Q_0\Psi_1,\delta x_i\rangle\\
&+\langle R\bar{u}_i,\delta u_i\rangle+\langle Q\Psi_3-\Theta^{T}Q\Psi_3-\alpha\Theta_0^{T}Q_0\Psi_1,\delta x^{\dagger}\rangle dt+\langle \alpha G_0\Psi_4(T)\\
&-\hat{\Theta}_1^{T}G\Psi_6(T),\delta x_0^{\dagger}(T)\rangle
+\langle G\Psi_6(T)-\hat{\Theta}^{T}G\Psi_6(T)-\alpha\hat{\Theta}_0^{T}G_0\Psi_4(T),\\
&\delta x^{\dagger}(T)\rangle+\langle G\Psi_5^i(T)-\hat{\Theta}^{T}G\Psi_6(T)-\alpha\hat{\Theta}_0^{T}G_0\Psi_4(T),\delta x_i(T)\rangle \bigg\},
\end{aligned}%$}
\end{equation}
where
\begin{equation*}
\left\{
\begin{aligned}
&\Psi_1(\cdot):=\bar{x}_0- \Theta_0\hat{x}-\eta_0,\quad \Psi_2^i(\cdot):=\bar{x}_i- \Theta \hat{x}-\Theta_1\bar{x}_0-\eta,\\
&\Psi_3(\cdot):=(I-\Theta)\hat{x}-\Theta_1\bar{x}_0-\eta,\\
\end{aligned}
\right.
\end{equation*}
are related to time $t$, and
\begin{equation*}
\left\{
\begin{aligned}
&\Psi_4(T):=\bar{x}_0(T)- \hat{\Theta}_0\hat{x}(T)-\hat{\eta}_0,\quad \Psi_5^i(T):=\bar{x}_i(T) -\hat{\Theta}\hat{x}(T) -\hat{\Theta}_1\bar{x}_0(T)-\hat{\eta},\\
&\Psi_6(T):=(I-\hat{\Theta})\hat{x}(T) -\hat{\Theta}_1\bar{x}_0(T)-\hat{\eta},\\
\end{aligned}
\right.
\end{equation*}
are related to time $T$ which are terminal terms.

It is very important to formulate an auxiliary control problem to obtain the decentralized optimal control for analyzing the problem of social optimality. Usually, an auxiliary control problem is a standard LQ control problem (see \cite{hcm2012,wh2017}). However, \eqref{lf_ctr_cost} contains $\delta x_0^{\dagger}$ and $\delta x^{\dagger}$, which are the terms we do not want them appear in the social cost functional. Therefore, we need to use a duality procedure (see \cite[Chapter 3]{yz1999}) to get off the dependence of $\delta \mathcal{J}_{soc}^{(N)}(\delta u_i)$ on $\delta x_0^{\dagger}$ and $\delta x^{\dagger}$.
To this end, we introduce two auxiliary equations %$k_1$ and $k_2$.
\begin{equation}\label{two_aux_origin}
\left\{
\begin{aligned}
&dk_1=\alpha_1dt+\beta_1dW_0,\ k_1(T)=\alpha G_0\Psi_4(T)-\hat{\Theta}_1^{T}G\Psi_6(T),\\
%\end{aligned}
%\right.
%\end{equation}
%\begin{equation}
%\left\{
%\begin{aligned}
&dk_2=\alpha_2dt+\beta_2dW_0,\ k_2(T)=(I-\hat{\Theta}^{T})G\Psi_6(T)-\alpha\hat{\Theta}_0^{T}G_0\Psi_4(T).
\end{aligned}
\right.
\end{equation}
Using It\^{o} formula, we have the following duality relations
\begin{equation}\label{ito_two_aux}
\begin{aligned}
&\mathbb{E}\langle \alpha G_0\Psi_4(T)-\hat{\Theta}_1^{T}G\Psi_6(T), \delta x_0^{\dagger}(T)\rangle\\
&\hspace{0cm}=\mathbb{E}\langle k_1(0), \delta x_0^{\dagger}(0)\rangle+\mathbb{E}\int_{0}^{T}\langle k_1,A_0\delta x_0^{\dagger}+C_0\delta x_i+C_0\delta x^{\dagger}\rangle+\langle \alpha_1, \delta x_0^{\dagger}\rangle dt,\\
%\end{aligned}
%\end{equation}
%and
%\begin{equation}
%\begin{aligned}
&\mathbb{E}\langle (I-\hat{\Theta}^{T})G\Psi_6(T)-\alpha\hat{\Theta}_0^{T}G_0\Psi_4(T), \delta x^{\dagger}(T)\rangle\\
&\hspace{0cm}=\mathbb{E}\langle k_2(0), \delta x^{\dagger}(0)\rangle+\mathbb{E}\int_{0}^{T}\langle k_2,A\delta x^{\dagger}+C\delta x_i+C\delta x^{\dagger}+F\delta x_0^{\dagger}\rangle+\langle \alpha_2, \delta x^{\dagger}\rangle dt.
\end{aligned}
\end{equation}
Putting \eqref{lf_ctr_cost} and \eqref{ito_two_aux} together, we obtain
\begin{equation*}%\resizebox{\linewidth}{!}{$
\begin{aligned}
\frac{1}{2}\delta \mathcal{J}_{soc}^{(N)}(\delta u_i)
%=&\mathbb{E}\bigg\{\int_{0}^{T}\langle \alpha Q_0\Psi_1-\Theta_1^{T}Q\Psi_3,\delta x_0^{\dagger}\rangle+\langle Q\Psi_2^i-\Theta^{T}Q\Psi_3-\alpha\Theta_0^{T}Q_0\Psi_1,\delta x_i\rangle\\
%&+\langle Q\Psi_3-\Theta^{T}Q\Psi_3-\alpha\Theta_0^{T}Q_0\Psi_1,\delta x^{\dagger}\rangle+\langle R\bar{u}_i,\delta u_i\rangle +\langle k_1,A_0\delta x_0^{\dagger}\\
%&+C_0\delta x_i+C_0\delta x^{\dagger}\rangle+\langle \alpha_1, \delta x_0^{\dagger}\rangle +\langle k_2,A\delta x^{\dagger}+C\delta x_i+C\delta x^{\dagger}+F\delta x_0^{\dagger}\rangle\\
%&+\langle \alpha_2, \delta x^{\dagger}\rangle dt
%+\langle G\Psi_5^i(T)-\hat{\Theta}^{T}G\Psi_6(T)-\alpha\hat{\Theta}_0^{T}G_0\Psi_4(T),\delta x_i(T)\rangle \bigg\}\\
=&\mathbb{E}\bigg\{\int_{0}^{T}\langle R\bar{u}_i,\delta u_i\rangle+\langle \alpha Q_0\Psi_1-\Theta_1^{T}Q\Psi_3+\alpha_1+F^Tk_2+A_0^Tk_1,\delta x_0^{\dagger}\rangle\\
&+\langle Q\Psi_3-\Theta^{T}Q\Psi_3-\alpha\Theta_0^{T}Q_0\Psi_1+C_0^Tk_1 +C^Tk_2+\alpha_2+A^Tk_2,\delta x^{\dagger}\rangle \\
&+\langle Q\Psi_2^i-\Theta^{T}Q\Psi_3-\alpha\Theta_0^{T}Q_0\Psi_1+C_0^Tk_1 +C^Tk_2,\delta x_i\rangle dt\\
&+\langle G\Psi_5^i(T)-\hat{\Theta}^{T}G\Psi_6(T)-\alpha\hat{\Theta}_0^{T}G_0\Psi_4(T),\delta x_i(T)\rangle \bigg\}.
\end{aligned}%$}
\end{equation*}
Comparing the coefficients, it follows that
\begin{equation*}
\left\{
\begin{aligned}
&\alpha_1=-(\alpha Q_0\Psi_1-\Theta_1^{T}Q\Psi_3+F^Tk_2+A_0^Tk_1),\\
&\alpha_2=-( Q\Psi_3-\Theta^{T}Q\Psi_3-\alpha\Theta_0^{T}Q_0\Psi_1+C_0^Tk_1 +C^Tk_2+A^Tk_2).
\end{aligned}
\right.
\end{equation*}
Then, according to above discussion, the two auxiliary equations can be rewritten as:
\begin{equation*}
\left\{
\begin{aligned}
&dk_1=-(\alpha Q_0\Psi_1-\Theta_1^{T}Q\Psi_3+F^Tk_2+A_0^Tk_1)dt+\beta_1dW_0,\\
&dk_2=-( Q\Psi_3-\Theta^{T}Q\Psi_3-\alpha\Theta_0^{T}Q_0\Psi_1+C_0^Tk_1 +C^Tk_2+A^Tk_2)dt+\beta_2dW_0,\\
&k_1(T)=\alpha G_0\Psi_4(T)-\hat{\Theta}_1^{T}G\Psi_6(T), \ k_2(T)=(I-\hat{\Theta}^{T})G\Psi_6(T)-\alpha\hat{\Theta}_0^{T}G_0\Psi_4(T),
\end{aligned}
\right.
\end{equation*}
and the variation of social cost functional is equivalent to
\begin{equation}\label{final_soc}
\begin{aligned}
\frac{1}{2}\delta &\mathcal{J}_{soc}^{(N)}(\delta u_i)
%=&\mathbb{E}\bigg\{\int_{0}^{T}\!\!\!\langle R\bar{u}_i,\delta u_i\rangle \!+\! \langle Q\Psi_2^i \!-\! \Theta^{T}Q\Psi_3 \!-\! \alpha\Theta_0^{T}Q_0\Psi_1 \!+\! C_0^Tk_1 \\ &+C^Tk_2,\delta x_i\rangle dt+\langle G\Psi_5^i(T)-\hat{\Theta}^{T}G\Psi_6(T)-\alpha\hat{\Theta}_0^{T}G_0\Psi_4(T),\delta x_i(T)\rangle \bigg\}.
%\end{aligned}
%\end{equation}
%Hence, we can construct a decentralized auxiliary cost functional %$\hat{\mathcal{J}}_i$
%$\hat{\mathcal{J}}_i((u_0,\hat{x});u_i)$ such that
%\begin{equation*}
%\begin{aligned}
%\hat{\mathcal{J}}_i((u_0,\hat{x});u_i)
=\mathbb{E}\bigg\{\int_{0}^{T}\langle Q\bar{x}_i,\delta x_i\rangle+\langle R\bar{u}_i,\delta u_i\rangle+\langle -Q(\Theta \hat{x}+\Theta_1\bar{x}_0+\eta)\\
&-\Theta^{T}Q\Psi_3-\alpha\Theta_0^{T}Q_0\Psi_1 +C_0^Tk_1 +C^Tk_2,\delta x_i\rangle dt
+\langle G\bar{x}_i(T),\delta x_i(T)\rangle\\
&+\langle -G(\hat{\Theta}\hat{x}(T)+\hat{\Theta}_1\bar{x}_0(T) +\hat{\eta})-\hat{\Theta}^{T}G\Psi_6(T)-\alpha\hat{\Theta}_0^{T}G_0\Psi_4(T),\delta x_i(T)\rangle \bigg\}.
\end{aligned}
\end{equation}

\subsection{Decentralized strategy design for followers}

As discussed in previous subsection, when $N$ is sufficiently large, a stochastic process $\hat{x}$ can
be used to approximate $x^{(N)}$. Now, we can introduce the following auxiliary control problem for the $i$th follower.
\begin{prblm}
\text{\rm(\textbf{P2})} Minimize $\hat{\mathcal{J}_i}((u_0,\hat{x});u_i)$ over $u_i\in\mathcal{U}_i$, where
\begin{equation}\label{auxi_ctrl_problem_state}
%\left\{
%\begin{aligned}
dx_i=[Ax_i+Bu_i+C\hat{x}+F\bar{x}_0(u_0)]dt+DdW_i, \  x_i(0)=\xi_{i}, \ i=1,2,\cdots,N,\\
\end{equation}
\begin{equation}\label{auxi_ctrl_problem_cost}
\hat{\mathcal{J}_i}((u_0,\hat{x});u_i)=\mathbb{E}\bigg\{\int_{0}^{T}\|x_i\|_Q^2 +\|u_i\|_R^2+2\langle \chi_1,x_i\rangle dt+\|x_i(T)\|_G^2+2\langle \chi_2,x_i(T)\rangle \bigg\},
%\end{aligned}
%\right.
\end{equation}
with
\begin{equation*}
\left\{
\begin{aligned}
&\chi_1=-Q(\Theta \hat{x}+\Theta_1\bar{x}_0(u_0)+\eta)-\Theta^{T}Q\Psi_3-\alpha\Theta_0^{T}Q_0\Psi_1 +C_0^Tk_1 +C^Tk_2,\\
&\chi_2=-G(\hat{\Theta}\hat{x}(T)+\hat{\Theta}_1\bar{x}_0(u_0)(T) +\hat{\eta})-\hat{\Theta}^{T}G\Psi_6(T)-\alpha\hat{\Theta}_0^{T}G_0\Psi_4(T).
\end{aligned}
\right.
\end{equation*}
Here, $\bar{x}_0(u_0)$ means $\bar{x}_0$ is related to $u_0$. $\bar{x}_0$, $\hat{x}$, $k_1$ and $k_2$ are determined by
\begin{equation}\label{pre_ccsys}%\resizebox{\linewidth}{!}{$
\left\{
\begin{aligned}
&d\bar{x}_0=[A_0\bar{x}_0+B_0u_0+C_0\hat{x}]dt+D_0dW_0, \quad \bar{x}_0(0)=\xi_0,\\
&d\hat{x}=[A\hat{x}+B\hat{u}+C\hat{x}+F\bar{x}_0(u_0)]dt,\quad  \hat{x}(0)=\hat{\xi},\\
&dk_1=-(\alpha Q_0\Psi_1-\Theta_1^{T}Q\Psi_3+F^Tk_2+A_0^Tk_1)dt+\beta_1dW_0,\\
&dk_2=-( Q\Psi_3-\Theta^{T}Q\Psi_3-\alpha\Theta_0^{T}Q_0\Psi_1+C_0^Tk_1 +C^Tk_2+A^Tk_2)dt+\beta_2dW_0,\\
&k_1(T)=\alpha G_0\Psi_4(T)-\hat{\Theta}_1^{T}G\Psi_6(T), \ k_2(T)=(I-\hat{\Theta}^{T})G\Psi_6(T)-\alpha\hat{\Theta}_0^{T}G_0\Psi_4(T),
\end{aligned}
\right.%$}
\end{equation}
where $\hat{x}$ and $\hat{u}$ are the approximations of $x^{(N)}$ and $\frac{1}{N}\sum_{i=1}^{N}u_i$, respectively.
\end{prblm}

In what follows, we let $\bar{u}=\mathcal{M}(u_0)=\{\bar{u}_1, \bar{u}_2,\cdots, \bar{u}_N\}\in \mathcal{U}$. %Thus, when $u_0$ is changed, $\bar{u}_i$ will be changed correspondingly.
Note that $\bar{u}$ here represents the decentraliezd optimal control, which is different from the same notation in the beginning of Section 2.
\begin{prpstn}
Assume that \text{\rm(\textbf{A1})-(\textbf{A3})} hold. For given $u_0\in \mathcal{U}_0$, \text{\rm(\textbf{P2})} has a unique optimal control
\begin{equation*}
\begin{aligned}
\bar{u}_i=-R^{-1}B^Tp_i,
\end{aligned}
\end{equation*}
where $p_i$ is an adaptive solution to the following backward stochastic differential equation \text{\rm(BSDE)}
\begin{equation*}
%\left\{
\begin{aligned}
dp_i=-(A^Tp_i+Q\bar{x}_i+\chi_1)dt+\zeta_0dW_0+\zeta_idW_i,\
%&p_i(T)=Gx_i(T)-G(\hat{\Theta}\hat{x}(T)+\hat{\Theta}_1\bar{x}_0(T) +\hat{\eta})\\
%&\hspace{1.2cm}-\hat{\Theta}^{T}G[(I-\hat{\Theta})\hat{x}(T) -\hat{\Theta}_1\bar{x}_0(T)-\hat{\eta}]-\alpha\hat{\Theta}_0^{T}G_0[\bar{x}_0(T)- \hat{\Theta}_0\hat{x}(T)-\hat{\eta}_0].
p_i(T)=G\bar{x_i}(T)+\chi_2.
\end{aligned}
%\right.
\end{equation*}
\end{prpstn}
\begin{proof}
%{\it Proof }
The variation of the state equation in \eqref{auxi_ctrl_problem_state} is
\begin{equation*}
%\left\{
\begin{aligned}
d\delta x_i=(A\delta x_i+B\delta u_i)dt,\
\delta x_i(0)=0, \quad i=1,2,\cdots,N,
\end{aligned}
%\right.
\end{equation*}
and the variation of the corresponding cost functional is
\begin{equation}\label{var_auxi_ctrl_problem_cost}
%\left\{
\begin{aligned}
&\frac{1}{2}\delta\hat{\mathcal{J}_i}((u_0,\hat{x});u_i)
%\frac{\tilde{J}((u_0,\hat{x});u_i)}{\delta u_i}\\
=\mathbb{E}\bigg\{\int_{0}^{T}\langle Q\bar{x}_i,\delta x_i\rangle +\langle R\bar{u}_i,\delta u_i\rangle+\langle \chi_1,\delta x_i\rangle dt+\langle G\bar{x}_i(T),\delta x_i(T)\rangle+\langle \chi_2,\delta x_i(T)\rangle \bigg\}.
\end{aligned}
%\right.
\end{equation}
Using a similar argument from \eqref{two_aux_origin} to \eqref{final_soc}, we construct the following auxiliary equation
%and use $p_i$ to substitute $\delta x_i$ in (24)
\begin{equation}\label{aux_plb_aux_eq}
%\left\{
\begin{aligned}
dp_i=\theta_idt+\zeta_0dW_0+\zeta_idW_i,\
p_i(T)=G\bar{x}_i(T)+\chi_2,\\
%&\hspace{0.9cm}=G\bar{x}_i(T)-G(\hat{\Theta}\hat{x}(T) +\hat{\Theta}_1\bar{x}_0(u_0)(T) +\hat{\eta})-\hat{\Theta}^{T}G\Psi_6-\alpha\hat{\Theta}_0^{T}G_0\Psi_4,
\end{aligned}
%\right.
\end{equation}
and have the following duality relation between $p_i$ and $\delta x_i$ by using It\^{o} formula
\begin{equation*}
\begin{aligned}
\mathbb{E}\langle p_i(T),\delta x_i(T)\rangle=&\mathbb{E}\langle p_i(0),\delta x_i(0)\rangle+\mathbb{E}\int_{0}^{T}\langle A^Tp_i+\theta_i,\delta x_i\rangle+\langle B^Tp_i,\delta u_i\rangle dt.
\end{aligned}
\end{equation*}
For given $u_0\in \mathcal{U}_0$, since $Q\geq 0$ and $R>\delta I$, for some $\delta>0$, \eqref{auxi_ctrl_problem_cost} is uniformly convex and it has a unique optimal control. Let \eqref{var_auxi_ctrl_problem_cost} equal zero, we have
\begin{equation}\label{bar_u_i}
%\left\{
\begin{aligned}
\theta_i=-(A^Tp_i+Q\bar{x}_i+\chi_1),\quad
%\end{aligned}
%\end{equation}
%and
%\begin{equation*}
%\begin{aligned}
\bar{u}_i=-R^{-1}B^Tp_i.
\end{aligned}
%\right.
\end{equation}
Then, the proposition follows.
%\qed
\end{proof}

Substituting \eqref{bar_u_i} into \eqref{auxi_ctrl_problem_state} and \eqref{aux_plb_aux_eq}, we have the following FBSDE
\begin{equation}\label{x_p_fbsde_1}%\resizebox{\linewidth}{!}{$
\left\{
\begin{aligned}
&d\bar{x}_i=[A\bar{x}_i-BR^{-1}B^Tp_i+C\hat{x}+F\bar{x}_0]dt+DdW_i,\ x_i(0)=\xi_{i}, \ i=1,2,\cdots,N,\\
&dp_i=-(A^Tp_i+Q\bar{x}_i+\chi_1)dt+\zeta_0dW_0+\zeta_idW_i, \ p_i(T)=Gx_i(T)+\chi_2.
\end{aligned}
\right.%$}
\end{equation}
By taking limits, the above FBSDE can be rewritten as:
\begin{equation}\label{lim_x_p_fbsde_1}
\left\{
\begin{aligned}
&d\hat{x}=[(A+C)\hat{x}+F\bar{x}_0-BR^{-1}B^T\hat{p}]dt,\ \hat{x}(0)=\hat{\xi}, \\
&d\hat{p}=-(A^T\hat{p}+Q\hat{x}+\chi_1)dt+\zeta_0dW_0, \ \hat{p}(T)=G\hat{x}(T)+\chi_2.
\end{aligned}
\right.
\end{equation}

\subsection{The consistency condition of the follower problem}

Let %By using the following notations to simplify our later presentations
\begin{equation*}%\resizebox{\linewidth}{!}{$
\left\{
 \begin{aligned}
   &\Xi_1:=(I-\Theta^T)Q(I-\Theta)+\alpha\Theta_0^TQ_0\Theta_0, \ \Xi_1^G:=(I-\hat{\Theta}^T)G(I-\hat{\Theta}) +\alpha\hat{\Theta}_0^TG_0\hat{\Theta}_0, \\
   &\Xi_2:=(I-\Theta^T)Q\Theta_1+\alpha\Theta_0^TQ_0, \ \Xi_2^G:=(I-\hat{\Theta}^T)G\hat{\Theta}_1+\alpha\hat{\Theta}_0^TG_0,\\
   &\Xi_3:=(I-\Theta^T)Q\eta-\alpha\Theta_0^TQ_0\eta_0, \ \Xi_3^G:=(I-\hat{\Theta}^T)G\hat{\eta}-\alpha\hat{\Theta}_0^TG_0\hat{\eta}_0,\\
   &\Xi_4:=\Theta_1^TQ\Theta_1+\alpha Q_0, \ \Xi_4^G:=\hat{\Theta}_1^TG\hat{\Theta}_1+\alpha G_0,\\
   &\Xi_5:=\Theta_1^TQ\eta-\alpha Q_0\eta_0, \ \Xi_5^G:=\hat{\Theta}_1^TG\hat{\eta}-\alpha G_0\hat{\eta}_0.
 \end{aligned}
\right.%$}
\end{equation*}
Combining \eqref{pre_ccsys} and \eqref{lim_x_p_fbsde_1}, we can obtain the CC system
\begin{equation}\label{ccsys1_3}%\resizebox{\linewidth}{!}{$
\left\{
\begin{aligned}
&d\hat{x}=[(A+C)\hat{x}+F\bar{x}_0-BR^{-1}B^Tk_2]dt,\ \hat{x}(0)=\hat{\xi}, \\
&d\bar{x}_0=[A_0\bar{x}_0+B_0u_0+C_0\hat{x}]dt+D_0dW_0,\ \bar{x}_0(0)=\xi_0,\\
%&d\hat{p}=-[A^T\hat{p}+\Xi_1\hat{x} -\Xi_2\bar{x}_0+C_0^Tk_1 +C^Tk_2-\Xi_3]dt+\zeta_0dW_t^0,\\
%&\hspace{0.8cm}\\
&dk_1=-[\Xi_4\bar{x}_0-\Xi_2^T\hat{x}+A_0^Tk_1 +F^Tk_2+\Xi_5]dt+\beta_1dW_0,\\
%&\hspace{1cm}\\
&dk_2=-[\Xi_1\hat{x} -\Xi_2\bar{x}_0+C_0^Tk_1+(A+C)^Tk_2-\Xi_3]dt+\beta_2dW_0,\\
%&\hspace{1cm}\\
%&\hat{p}(T)=\Xi_1^G\hat{x}(T) -\Xi_2^G\bar{x}_0(T)-\Xi_3^G,\\
%&\hspace{1.2cm}\\
&k_1(T)=\Xi_4^G\bar{x}_0(T)- (\Xi_2^G)^T\hat{x}(T)+\Xi_5^G, \ k_2(T)=\Xi_1^G\hat{x}(T) -\Xi_2^G\bar{x}_0(T)-\Xi_3^G,
%&\hspace{1.2cm}
\end{aligned}
\right.%$}
\end{equation}
where $\hat{p}=k_2$ can be easily verified.

\section{The optimal control problem for the leader}

Now, let (\textbf{P2}) have a unique solution. Then, for $u_0\in\mathcal{U}_0$ given by leader,
%In the auxiliary control problem, for the fix $u_0$ given by leader,
the followers choose their optimal control $\bar{u}=\mathcal{M}(u_0)=\{\bar{u}_1, \bar{u}_2,\cdots, \bar{u}_N\}\in \mathcal{U}$, where $\bar{u}_i$ is shown in \eqref{bar_u_i}. Now we consider the optimal control of the leader to further minimize the social cost functional. In the infinite population system, $x^{(N)}$ may be approximated by $\hat{x}$. Hence, we can construct the following auxiliary optimal control problem for the leader.
\begin{prblm}
\text{\rm(\textbf{P3})} Minimize $\hat{\mathcal{J}}_{soc}^{(N)}(u_0; \bar{u})$ over $u_0\in\mathcal{U}_0$, where
\begin{equation}\label{leader_prob_state}
\begin{aligned}
&dx_0=[A_0x_0+B_0u_0+C_0\hat{x}]dt+D_0dW_0, \quad x_0(0)=\xi_0,\\
&\hat{\mathcal{J}}_{soc}^{(N)}(u_0; \bar{u})=\alpha N\hat{\mathcal{J}}_0(u_0; \bar{u})+\sum_{i=1}^{N}\hat{\mathcal{J}}_i(u_0; \bar{u}).
\end{aligned}
\end{equation}
\text{\rm(\textbf{P3})} is based on \text{\rm(\textbf{P2})}. Therefore, combining \eqref{bar_u_i}, the equations below \eqref{lf_ctr_cost} and the equations (2), (3) with mean field approximations, the cost functionals of the leader and the $i$th follower are
\begin{equation*}
\begin{aligned}
\hat{\mathcal{J}}_0(u_0; \bar{u})=&\mathbb{E}\bigg\{\int_{0}^{T}\langle Q_0\Psi_1,\Psi_1\rangle+\langle R_0u_0,u_0\rangle dt+\langle G_0\Psi_4,\Psi_4\rangle\bigg\},\\
%\end{aligned}
%\end{equation}
%and
%\begin{equation}
%\begin{aligned}
%\sum_{i=1}^{N}
\hat{\mathcal{J}}_i(u_0; \bar{u})=&\mathbb{E}\bigg\{\int_{0}^{T}\langle Q\Psi_2^i,\Psi_2^i\rangle+\langle B^Tp_i,R^{-1}B^Tp_i\rangle dt+\langle G\Psi_5^i,\Psi_5^i\rangle\bigg\},
\end{aligned}
\end{equation*}
where $\hat{x}$, $\bar{x}_0$, $k_1$, $k_2$, $\bar{x}_i$, $p_i$ are determined by \eqref{ccsys1_3} and \eqref{x_p_fbsde_1}.
%with
%\begin{equation*}
%\left\{
%\begin{aligned}
%&d\hat{x}=[(A+C)\hat{x}+F\bar{x}_0-BR^{-1}B^Tk_2]dt,\ \hat{x}(0)=\hat{\xi}, \\
%&d\bar{x}_0=[A_0\bar{x}_0+B_0u_0+C_0\hat{x}]dt+D_0dW_0,\ \bar{x}_0(0)=\xi_0,\\
%%&d\hat{p}=-[A^T\hat{p}+\Xi_1\hat{x} -\Xi_2\bar{x}_0+C_0^Tk_1 +C^Tk_2-\Xi_3]dt+\zeta_0dW_t^0,\\
%%&\hspace{0.8cm}\\
%&dk_1=-[\Xi_4\bar{x}_0-\Xi_2^T\hat{x}+A_0^Tk_1 +F^Tk_2+\Xi_5]dt+\beta_1dW_0,\\
%%&\hspace{1cm}\\
%&dk_2=-[\Xi_1\hat{x} -\Xi_2\bar{x}_0+C_0^Tk_1+(A+C)^Tk_2-\Xi_3]dt+\beta_2dW_0,\\
%%&\hspace{1cm}\\
%%&\hat{p}(T)=\Xi_1^G\hat{x}(T) -\Xi_2^G\bar{x}_0(T)-\Xi_3^G,\\
%%&\hspace{1.2cm}\\
%&k_1(T)=\Xi_4^G\bar{x}_0(T)- (\Xi_2^G)^T\hat{x}(T)+\Xi_5^G, \ k_2(T)=\Xi_1^G\hat{x}(T) -\Xi_2^G\bar{x}_0(T)-\Xi_3^G,
%%&\hspace{1.2cm}
%\end{aligned}
%\right.
%\end{equation*}
%and
%\begin{equation*}
%\left\{
%\begin{aligned}
%&d\bar{x}_i=[A\bar{x}_i-BR^{-1}B^Tp_i+C\hat{x}+F\bar{x}_0]dt+DdW_i,\ x_i(0)=\xi_{i}, \ i=1,2,\cdots,N,\\
%&dp_i=-(A^Tp_i+Q\bar{x}_i+\chi_1)dt+\zeta_0dW_0+\zeta_idW_i, \ p_i(T)=Gx_i(T)+\chi_2.
%\end{aligned}
%\right.
%\end{equation*}
\end{prblm}
Using a similar argument in Section 2, we let $\bar{u}_0$ be the optimal strategy of the leader and perturb $u_0$ in \eqref{leader_prob_state}, where $\delta u_0=u_0-\bar{u}_0$. Since $\bar{x}_0$, $\hat{x}$, $\bar{x}_i$ and $p_i$ are determined by $u_0$, we denote their corresponding perturbations as: $\delta \bar{x}_0=\bar{x}_0(u_0)-\bar{x}_0(\bar{u}_0)$, $\delta \hat{x}=\hat{x}(u_0)-\hat{x}(\bar{u}_0)$, $\delta \bar{x}_i=\bar{x}_i(u_0)-\bar{x}_i(\bar{u}_0)$ and $\delta p_i=p_i(u_0)-p_i(\bar{u}_0)$. For sake of notation simplicity, we drop $(\bar{u}_0)$ in the following $\bar{x}_0(\bar{u}_0)$, $\hat{x}(\bar{u}_0)$, $\bar{x}_i(\bar{u}_0)$ and $p_i(\bar{u}_0)$, etc. Then, one can obtain %the following variation of $\bar{x}_0$ of the leader in auxiliary control problem
\begin{equation*}
%\left\{
\begin{aligned}
d\delta\bar{x}_0=[A_0\delta\bar{x}_0+B_0\delta u_0+ C_0\delta\hat{x}]dt,\ \delta\bar{x}_0(0)=0,\\
\end{aligned}
%\right.
\end{equation*}
and the variations of corresponding cost functionals
\begin{equation*}%\resizebox{\linewidth}{!}{$
\begin{aligned}
\frac{1}{2}\delta\hat{\mathcal{J}}_0(\delta u_0)=&\mathbb{E}\bigg\{\int_{0}^{T}\langle Q_0\Psi_1,\delta \bar{x}_0- \Theta_0\delta\hat{x}\rangle+\langle R_0\bar{u}_0,\delta u_0\rangle dt+\langle G_0\Psi_4,\delta \bar{x}_0(T)- \hat{\Theta}_0\delta\hat{x}(T)\rangle\bigg\},\\
%\end{aligned}
%\end{equation}
%and
%\begin{equation}
%\begin{aligned}
\frac{1}{2}\sum_{i=1}^{N}\delta\hat{\mathcal{J}}_i(\delta &u_0)=\sum_{i=1}^{N}\mathbb{E}\bigg\{\int_{0}^{T}\langle Q\Psi_2^i,\delta \bar{x}_i- \Theta\delta\hat{x}-\Theta_1\delta \bar{x}_0\rangle+\langle R^{-1}B^Tp_i,B^T\delta p_i\rangle dt\\
&\hspace{1cm}+\langle G\Psi_5^i,\delta \bar{x}_i(T)- \hat{\Theta}\delta\hat{x}(T)-\hat{\Theta}_1\delta \bar{x}_0(T)\rangle\bigg\}.
\end{aligned}%$}
\end{equation*}
Here $\Psi_1$, $\Psi_2^i$, $\Psi_4(T)$, $\Psi_5^i(T)$, are related to $\bar{u}_0$. In what follows, $\Psi_1$, $\Psi_2^i$, $\Psi_3$, $\Psi_4(T)$, $\Psi_5^i(T)$, $\Psi_6(T)$ will be related to $\bar{u}_0$. Therefore, the variation of the social cost functional is
\begin{equation}\label{leader_soc}%\resizebox{\linewidth}{!}{$
\begin{aligned}
&\frac{1}{2}\delta\hat{\mathcal{J}}_{soc}^{(N)}(\delta u_0)\\
=&\alpha N\mathbb{E}\int_{0}^{T}\langle Q_0\Psi_1,\delta \bar{x}_0\rangle-\langle \Theta_0^TQ_0\Psi_1, \delta\hat{x}\rangle+\langle R_0\bar{u}_0,\delta u_0\rangle dt\\
&+\sum_{i=1}^{N}\mathbb{E}\int_{0}^{T}\langle Q\Psi_2^i,\delta \bar{x}_i\rangle-\langle \Theta^TQ\Psi_2^i, \delta\hat{x}\rangle-\langle \Theta_1^TQ\Psi_2^i,\delta \bar{x}_0\rangle+\langle BR^{-1}B^Tp_i,\delta p_i\rangle dt\\
&+\alpha N\langle G_0\Psi_4(T),\delta \bar{x}_0(T)\rangle-\alpha N\langle \hat{\Theta}_0^TG_0\Psi_4(T),\delta\hat{x}(T)\rangle-\sum_{i=1}^{N}\langle \hat{\Theta}_1^TG\Psi_5^i(T),\delta \bar{x}_0(T)\rangle \\
& -\sum_{i=1}^{N}\langle \hat{\Theta}^TG\Psi_5^i(T),\delta\hat{x}(T)\rangle+\sum_{i=1}^{N}\langle G\Psi_5^i(T),\delta \bar{x}_i(T)\rangle.
\end{aligned}%$}
\end{equation}
Similarly, the variations of those equations in \eqref{x_p_fbsde_1} and \eqref{ccsys1_3} are given by
\begin{equation*}%\resizebox{\linewidth}{!}{$
\left\{
\begin{aligned}
&d\delta \bar{x}_i=[A\delta \bar{x}_i-BR^{-1}B^T\delta p_i+C\delta\hat{x}+F\delta \bar{x}_0]dt,\quad \delta \bar{x}_i(0)=0, \quad i=1,2,\cdots,N,\\
&d\delta p_i=-(A^T\delta p_i+Q\delta x_i+[\Xi_1-Q]\delta\hat{x}-\Xi_2\delta \bar{x}_0+C_0^T\delta k_1 +C^T\delta k_2)dt+\delta\zeta_0dW_0+\delta\zeta_idW_i,\\
&\delta p_i(T)=G\delta x_i(T)+[\Xi_1^G-G]\delta\hat{x}(T)-\Xi_2^G\delta \bar{x}_0(T),
\end{aligned}
\right.%$}
\end{equation*}
and
\begin{equation*}
\left\{
\begin{aligned}
&d\delta\hat{x}=[(A+C)\delta\hat{x}+F\delta \bar{x}_0-BR^{-1}B^T\delta k_2]dt,\quad
\delta\hat{x}(0)=0,\\
%&d\delta \bar{x}_0=[A_0\delta \bar{x}_0+C_0\delta\hat{x}+B_0\delta u_0]dt, \quad
%\delta \bar{x}_0(0)=0,\\
&d\delta k_1=-[\Xi_4\delta \bar{x}_0-\Xi_2^T\delta\hat{x}+A_0^T\delta k_1 +F^T\delta k_2]dt+\delta\beta_1dW_0,\\
&d\delta k_2=-[\Xi_1\delta\hat{x} -\Xi_2\delta \bar{x}_0+C_0^T\delta k_1+(A+C)^T\delta k_2]dt+\delta\beta_2dW_0,\\
&\delta k_1(T)=\Xi_4^G\delta \bar{x}_0(T)- (\Xi_2^G)^T\delta\hat{x}(T),\quad \delta k_2(T)=\Xi_1^G\delta\hat{x}(T) -\Xi_2^G\delta \bar{x}_0(T).\\
\end{aligned}
\right.
\end{equation*}
Since \eqref{leader_soc} contains many terms that we do not want them appear in the social cost functional, we will use a similar argument in Section 2 to get off the dependence of $\delta\hat{\mathcal{J}}_{soc}^{(N)}(\delta u_0)$ on those terms. Therefore, we need to construct six auxiliary equations to help us obtain the optimal control of the leader. We introduce the first three auxiliary equations:
\begin{equation*}
\left\{
\begin{aligned}
&dq_i=m_idt+n_i^0dW_0+n_idW_i,\quad q_i(0)=0, \quad i=1,2,\cdots,N,\\
%\end{aligned}
%\right.
%\end{equation}
%and
%\begin{equation}
%\left\{
%\begin{aligned}
%&d\hat{y}=\hat{\alpha}dt+0dW^0,\quad \hat{y}(T)= ,\\
%&dy_0=\hat{\alpha}_0dt+\hat{\beta}_0dW^0,\quad y_0(T)=
%\end{aligned}
%\right.
%\end{equation}
%\begin{equation}
%\left\{
%\begin{aligned}
&dl_1=s_1dt+r_1dW_0,\quad l_1(0)=0,\\
&dl_2=s_2dt+r_2dW_0,\quad l_2(0)=0.
\end{aligned}
\right.
\end{equation*}
Here $q_i$, $l_1$ and $l_2$ are used to free $\delta\hat{\mathcal{J}}_{soc}^{(N)}(\delta u_0)$ from the dependence on $p_i$, $k_1$ and $k_2$, respectively. By a similar argument from \eqref{two_aux_origin} to \eqref{final_soc}, we can rewrite the variation of the social cost functional as follows:
\begin{equation*}%\resizebox{\linewidth}{!}{$
\begin{aligned}
&\frac{1}{2}\delta\hat{\mathcal{J}}_{soc}^{(N)}(\delta u_0)\\
&=\mathbb{E}\int_{0}^{T}\sum_{i=1}^{N}\langle \alpha Q_0\Psi_1- \Theta_1^TQ\Psi_2^i,\delta \bar{x}_0\rangle+\sum_{i=1}^{N}\langle -\alpha \Theta_0^TQ_0\Psi_1-\Theta^TQ\Psi_2^i, \delta\hat{x}\rangle+\sum_{i=1}^{N}\langle Q\Psi_2^i,\delta \bar{x}_i\rangle\\
&+N\langle \alpha R_0\bar{u}_0,\delta u_0\rangle +\sum_{i=1}^{N}\langle BR^{-1}B^Tp_i,\delta p_i\rangle+ \sum_{i=1}^{N}\langle l_1, -(A_0^T\delta k_1+F^T\delta k_2+\Xi_4\delta \bar{x}_0-\Xi_2^T\delta \hat{x})\rangle\\
&+\sum_{i=1}^{N}\langle s_1,\delta k_1\rangle+\sum_{i=1}^{N}\langle s_2,\delta k_2\rangle+\sum_{i=1}^{N}\langle l_2, -((A+C)^T\delta k_2+C_0^T\delta k_1+\Xi_1\delta\hat{x} -\Xi_2\delta \bar{x}_0)\rangle \\
&+\sum_{i=1}^{N}\langle -(A^T\delta p_i+Q\delta \bar{x}_i+[\Xi_1-Q]\delta\hat{x}-\Xi_2\delta \bar{x}_0+C_0^T\delta k_1 +C^T\delta k_2), q_i\rangle +\sum_{i=1}^{N}\langle \delta p_i, m_i\rangle\\
&+\sum_{i=1}^{N}[\langle \delta \beta_1, r_1\rangle+\langle \delta \beta_2, r_2\rangle]+\sum_{i=1}^{N}[\langle \delta\zeta_0, n_i^0\rangle+\langle\delta\zeta_i, n_i\rangle] dt+\sum_{i=1}^{N}\langle G\Psi_5^i(T)-Gq_i(T),\delta \bar{x}_i(T)\rangle\\
&+\sum_{i=1}^{N}\langle \alpha  G_0\Psi_4(T)-\hat{\Theta}_1^TG\Psi_5^i(T) -(\Xi_4^G)^Tl_1(T) +(\Xi_2^G)^Tl_2(T)+(\Xi_2^G)^Tq_i(T) ,\delta \bar{x}_0(T)\rangle\\
&-\sum_{i=1}^{N}\langle\alpha \hat{\Theta}_0^TG_0\Psi_4(T)+\hat{\Theta}^TG\Psi_5^i(T) -(\Xi_2^G)^Tl_1(T)+(\Xi_1^G)^Tl_2(T) +(\Xi_1^G-G)^Tq_i(T),\delta\hat{x}(T)\rangle.\\
\end{aligned}%$}
\end{equation*}
Next, we introduce another three auxiliary equations:
\begin{equation*}%\resizebox{\linewidth}{!}{$
\left\{
\begin{aligned}
&d\hat{y}^i=\hat{\alpha}dt+\hat{\beta}dW_0+\sum_{i=1}^{N}\hat{\beta}_idW_i,\\
&\hat{y}^i(T)=\alpha \hat{\Theta}_0^TG_0\Psi_4(T)+\hat{\Theta}^TG\Psi_5^i(T) -(\Xi_2^G)^Tl_1(T)+(\Xi_1^G)^Tl_2(T)+(\Xi_1^G-G)^Tq_i(T),\\
&dy_0^i=\hat{\alpha}_0dt+\hat{\beta}_0dW_0+\sum_{i=1}^{N}\hat{\beta}_i^0dW_i,\\ &y_0^i(T)=\alpha G_0\Psi_4(T)-\hat{\Theta}_1^TG\Psi_5^i(T)-(\Xi_4^G)^Tl_1(T) +(\Xi_2^G)^Tl_2(T)+(\Xi_2^G)^Tq_i(T),\\
%\end{aligned}
%\right.
%\end{equation}
%and
%\begin{equation}
%\left\{
%\begin{aligned}
&dy_i=\alpha_idt+\beta_0dW_0+\beta_idW_i,\quad y_i(T)=G\Psi_5^i(T)-Gq_i(T),  \quad i=1,2,\cdots,N.
\end{aligned}
\right.%$}
\end{equation*}
Here $\hat{y}^i$, $y_0^i$ and $y_i$ are used to free $\delta\hat{\mathcal{J}}_{soc}^{(N)}(\delta u_0)$ from the dependence on $\delta \hat{x}$, $\delta \bar{x}_0$ and $\delta \bar{x}_i$, respectively. Similarly, by It\^{o} formula and the duality relations,
the variation of the social cost functional can be further rewritten as follows:
\begin{equation*}%\resizebox{\linewidth}{!}{$
\begin{aligned}
&\frac{1}{2}\delta\hat{\mathcal{J}}_{soc}^{(N)}(\delta u_0)\\
=&\mathbb{E}\int_{0}^{T}\sum_{i=1}^{N}\langle \alpha Q_0\Psi_1- \Theta_1^TQ\Psi_2^i+F^Ty_i-F^T\hat{y}^i+A_0^Ty_0^i +\hat{\alpha}_0-\Xi_4^Tl_1+\Xi_2^Tl_2 \\
&+\Xi_2^Tq_i,\delta \bar{x}_0\rangle+\sum_{i=1}^{N}\langle -\alpha \Theta_0^TQ_0\Psi_1-\Theta^TQ\Psi_2^i +C^Ty_i-(A+C)^T\hat{y}^i+C_0^Ty_0^i\\
%& \\
&-\hat{\alpha}+\Xi_2l_1-\Xi_1^Tl_2-(\Xi_1-Q)^Tq_i, \delta\hat{x}\rangle+\sum_{i=1}^{N}\langle Q\Psi_2^i+A^Ty_i+\alpha_i-Qq_i,\delta \bar{x}_i\rangle\\
%& \\
&+\sum_{i=1}^{N}\langle BR^{-1}B^Tp_i-BR^{-1}B^Ty_i+m_i-Aq_i,\delta p_i\rangle +\sum_{i=1}^{N}[\langle n_i, \delta \zeta_i\rangle+\langle n_i^0, \delta \zeta_0\rangle]\\
&+\sum_{i=1}^{N}\langle s_1-C_0l_2-A_0l_1+C_0q_i,\delta k_1\rangle+\sum_{i=1}^{N}\langle s_2-(A+C)l_2+BR^{-1}B^T\hat{y}^i-Fl_1\\
&+Cq_i,\delta k_2\rangle+\sum_{i=1}^{N}[\langle r_1,\delta\beta_1\rangle+\langle r_2,\delta\beta_2\rangle]+\langle \alpha NR_0\bar{u}_0+\sum_{i=1}^{N}B_0^Ty_0^i,\delta u_0\rangle dt.\\
%&
\end{aligned}%$}
\end{equation*}
Using a similar discussion in Proposition 2.1 and comparing the coeffcients, we have
\begin{equation*}%\resizebox{\linewidth}{!}{$
\left\{
\begin{aligned}
&\hat{\alpha}_0=-(\alpha Q_0\Psi_1- \Theta_1^TQ\Psi_2^i+F^Ty_i-F^T\hat{y}^i+A_0^Ty_0^i -\Xi_4^Tl_1+\Xi_2^Tl_2+\Xi_2^Tq_i),\\
&\hat{\alpha}=-\alpha \Theta_0^TQ_0\Psi_1-\Theta^TQ\Psi_2^i +C^Ty_i-(A+C)^T\hat{y}^i+C_0^Ty_0^i +\Xi_2l_1-\Xi_1^Tl_2-(\Xi_1-Q)^Tq_i,\\
&\alpha_i=-(Q\Psi_2^i+A^Ty_i-Q^Tq_i),\quad m_i=-(BR^{-1}B^Tp_i-BR^{-1}B^Ty_i-Aq_i),\\
&s_1=C_0l_2+A_0l_1-C_0q_i,\quad s_2=(A+C)l_2-BR^{-1}B^T\hat{y}^i+Fl_1-Cq_i,\\
&n_i=0,\quad n_i^0=0,\quad r_1=0, \quad r_2=0.\\
\end{aligned}
\right.%$}
\end{equation*}
%Hence,we denote $u_0^{(N)}=\bar{u}_0$.by (53) and (54),
%Note that the $\bar{u}_0$ is
Then, we have following FBSDE
\begin{equation}\label{y_l_fbsde_1}%\resizebox{\linewidth}{!}{$
\left\{
\begin{aligned}
&dy_i=-(A^Ty_i-Q^Tq_i+Q\Psi_2^i)dt+\beta_0dW_0+\beta_idW_i,\quad y_i(T)=G\Psi_5^i(T)-Gq_i(T), \\
%\end{aligned}
%\right.
%\end{equation}
%\begin{equation}
%\left\{
%\begin{aligned}
&dq_i=(BR^{-1}B^Ty_i+Aq_i-BR^{-1}B^Tp_i)dt,\quad q_i(0)=0, \quad i=1,2,\cdots,N,\\
%\end{aligned}
%\right.
%\end{equation}
%and
%\begin{equation}
%\left\{
%\begin{aligned}
&d\hat{y}^i=(-\alpha \Theta_0^TQ_0\Psi_1-\Theta^TQ\Psi_2^i +C^Ty_i-(A+C)^T\hat{y}^i+C_0^Ty_0^i +\Xi_2l_1-\Xi_1^Tl_2\\ &\hspace{0.8cm}-(\Xi_1-Q)^Tq_i)dt+\hat{\beta}dW_0+\sum_{i=1}^{N}\hat{\beta}_idW_i,\\
&\hat{y}^i(T)=\alpha \hat{\Theta}_0^TG_0\Psi_4(T)+\hat{\Theta}^TG\Psi_5^i(T) -(\Xi_2^G)^Tl_1(T)+(\Xi_1^G)^Tl_2(T)+(\Xi_1^G-G)^Tq_i(T) ,\\
&dy_0^i=-(\alpha Q_0\Psi_1- \Theta_1^TQ\Psi_2^i+F^Ty_i-F^T\hat{y}^i+A_0^Ty_0^i -\Xi_4^Tl_1+\Xi_2^Tl_2+\Xi_2^Tq_i)dt\\
&\hspace{1cm}+\hat{\beta}_0dW_0+\sum_{i=1}^{N}\hat{\beta}_i^0dW_i,\\ &y_0^i(T)=\alpha G_0\Psi_4(T)-\hat{\Theta}_1^TG\Psi_5^i(T)-(\Xi_4^G)^Tl_1(T) +(\Xi_2^G)^Tl_2(T)+(\Xi_2^G)^Tq_i(T),\\
%\end{aligned}
%\right.
%\end{equation}
%and
%\begin{equation}
%\left\{
%\begin{aligned}
&dl_1=(A_0l_1+C_0l_2-C_0q_i)dt,\quad l_1(0)=0,\\
&dl_2=[Fl_1+(A+C)l_2-BR^{-1}B^T\hat{y}^i-Cq_i]dt,\quad l_2(0)=0,
\end{aligned}
\right.%$}
\end{equation}
with the centralized form of the optimal control for the leader
%\begin{equation*}
%\begin{aligned}
%\alpha NR_0\bar{u}_0+\sum_{i=1}^{N}B_0^Ty_0^i=0.
%\end{aligned}
%\end{equation*}
\begin{equation}\label{u0*_origin}
\begin{aligned}
%u_0^{(N)}=\bar{u}_0=-\frac{\alpha^{-1}}{N}R_0^{-1}B_0^Ty_0.
\bar{u}_0=-\frac{\alpha^{-1}}{N}R_0^{-1}B_0^T\sum_{i=1}^{N}y_0^i := u_0^{(N)}.
\end{aligned}
\end{equation}
Note that $\bar{u}_0$ relies on $N$.
Denote
\begin{equation*}
\left\{
\begin{aligned}
&y^*=\lim_{N\rightarrow+\infty}\frac{1}{N}\sum_{i=1}^{N}y_i, \quad \hat{y}^*=\lim_{N\rightarrow+\infty}\frac{1}{N}\sum_{i=1}^{N}\hat{y}^i,\quad y_0^*=\lim_{N\rightarrow+\infty}\frac{1}{N}\sum_{i=1}^{N}y_0^i,\\ &q^*=\lim_{N\rightarrow+\infty}\frac{1}{N}\sum_{i=1}^{N}q_i,\quad
l_1^*=\lim_{N\rightarrow+\infty}\frac{1}{N}\sum_{i=1}^{N}l_1,\quad
l_2^*=\lim_{N\rightarrow+\infty}\frac{1}{N}\sum_{i=1}^{N}l_2.
\end{aligned}
\right.
\end{equation*}
%where the convergence are in $L^2$.
Here, using a similar argument of \eqref{x_dagger}, we can easily prove that $\frac{1}{N}\sum_{i=1}^{N}y_i$, $\frac{1}{N}\sum_{i=1}^{N}\hat{y}^i$, $\frac{1}{N}\sum_{i=1}^{N}y_0^i$, $\frac{1}{N}\sum_{i=1}^{N}q_i$, $\frac{1}{N}\sum_{i=1}^{N}l_1$ and $\frac{1}{N}\sum_{i=1}^{N}l_2$ converge to $y^*$, $\hat{y}^*$, $y_0^*$, $q^*$, $l_1^*$ and $l_2^*$, respectively. Thus, combining \eqref{ccsys1_3} and \eqref{y_l_fbsde_1}, when $N\rightarrow \infty$,
we can obtain the CC system for the leader-follower problem
\begin{equation}\label{ccsys2_1}%\resizebox{\linewidth}{!}{$
\left\{
\begin{aligned}
%&dx_i=[Ax_i-BR^{-1}B^Tp_i+C\hat{x}+Fx_0]dt+DdW^i,\quad x_i(0)=\xi_{i}, \quad i=1,2,\cdots,N,\\
%&dp_i=-(A^Tp_i+Qx_i+[\Xi_1-Q]\hat{x}-\Xi_2x_0-\Xi_3+C_0^Tk_1 +C^Tk_2)dt+\zeta_0dW^0+\zeta_idW^i,\\
%&p_i(T)=Gx_i(T)+[\Xi_1^G-G]\hat{x}(T)-\Xi_2^Gx_0(T)-\Xi_3^G,\\
%\end{aligned}
%\right.
%\end{equation}
%where
%\begin{equation}
%\left\{
%\begin{aligned}
&d\hat{x}=[(A+C)\hat{x}+F\bar{x}_0-BR^{-1}B^Tk_2]dt,\quad \hat{x}(0)=\hat{\xi}, \\
&d\bar{x}_0=[A_0\bar{x}_0+C_0\hat{x}-B_0(\alpha R_0)^{-1}B_0^Ty_0^*]dt+D_0dW_0,\quad \bar{x}_0(0)=\xi_0,\\
&dk_1=-[\Xi_4\bar{x}_0-\Xi_2^T\hat{x}+A_0^Tk_1 +F^Tk_2+\Xi_5]dt+\beta_1dW_0,\\
&dk_2=-[\Xi_1\hat{x} -\Xi_2\bar{x}_0+C_0^Tk_1+(A+C)^Tk_2-\Xi_3]dt+\beta_2dW_0,\\
&k_1(T)=\Xi_4^G\bar{x}_0(T)- (\Xi_2^G)^T\hat{x}(T)+\Xi_5^G, \quad k_2(T)=\Xi_1^G\hat{x}(T) -\Xi_2^G\bar{x}_0(T)-\Xi_3^G,\\
%\end{aligned}
%\right.
%\end{equation}
%\begin{equation}
%\left\{
%\begin{aligned}
%&dy^*=-(A^Ty^*-Q^Tq^*+Q\Psi_3)dt+\beta^*dW_0,\\ &y^*(T)=G\Psi_6(T)-Gq^*(T), \\
%&dq^*=(BR^{-1}B^Ty^*+Aq^*-BR^{-1}B^Tk_2)dt,\quad q^*(0)=0,\\
%&d\hat{y}^*=[-\alpha \Theta_0^TQ_0(\bar{x}_0- \Theta_0\hat{x}-\eta_0)-\Theta^TQ((I-\Theta)\hat{x}-\Theta_1\bar{x}_0-\eta) +C^Ty^* \\ &\hspace{0.8cm}-(A+C)^T\hat{y}^*+C_0^Ty_0^*+\Xi_2^Tl_1^*-\Xi_1^Tl_2^*-(\Xi_1-Q)^Tq^*]dt+\hat{\beta}^*dW_0,\\
%&\hat{y}^*(T)=\alpha \hat{\Theta}_0^TG_0(\bar{x}_0(T)- \hat{\Theta}_0\hat{x}(T)-\hat{\eta}_0)+\hat{\Theta}^TG((I-\hat{\Theta})\hat{x}(T) -\hat{\Theta}_1\bar{x}_0(T)-\hat{\eta})\\
%&\hspace{1.4cm}-(\Xi_2^G)^Tl_1^*(T)+(\Xi_1^G)^Tl_2^*(T)+(\Xi_1^G-G)^Tq^*(T) ,\\
%&dy_0^*=-[\alpha Q_0(\bar{x}_0- \Theta_0\hat{x}-\eta_0)- \Theta_1^TQ((I-\Theta)\hat{x}-\Theta_1\bar{x}_0-\eta)+F^Ty^*-F^T\hat{y}^*\\
%&\hspace{1cm}+A_0^Ty_0^* -\Xi_4^Tl_1^*+\Xi_2^Tl_2^*+\Xi_2^Tq^*]dt+\hat{\beta}_0^*dW_0,\\
%&y_0^*(T)=\alpha  G_0(\bar{x}_0(T)- \hat{\Theta}_0\hat{x}(T)-\hat{\eta}_0)-\hat{\Theta}_1^TG((I-\hat{\Theta})\hat{x}(T) -\hat{\Theta}_1\bar{x}_0(T)-\hat{\eta})\\
%&\hspace{1.4cm}-(\Xi_4^G)^Tl_1^*(T) +(\Xi_2^G)^Tl_2^*(T)+(\Xi_2^G)^Tq^*(T),\\
%&dl_1^*=(A_0l_1^*+C_0l_2^*-C_0q^*)dt,\quad l_1^*(0)=0,\\
%&dl_2^*=[Fl_1^*+(A+C)^Tl_2^*-BR^{-1}B^T\hat{y}^*-Cq^*]dt,\quad l_2^*(0)=0.
&dy^*=-(A^Ty^*-Q^Tq^*+Q\Psi_3)dt+\beta^*dW_0,\quad y^*(T)=G\Psi_6(T)-Gq^*(T), \\
&dq^*=(BR^{-1}B^Ty^*+Aq^*-BR^{-1}B^Tk_2)dt,\quad q^*(0)=0,\\
&d\hat{y}^*=[-\alpha \Theta_0^TQ_0\Psi_1-\Theta^TQ\Psi_3 +C^Ty^*-(A+C)^T\hat{y}^*+C_0^Ty_0^* +\Xi_2l_1^*-\Xi_1^Tl_2^*\\ &\hspace{0.8cm}-(\Xi_1-Q)^Tq^*]dt+\hat{\beta}^*dW_0,\\
&\hat{y}^*(T)=\alpha \hat{\Theta}_0^TG_0\Psi_4(T)+\hat{\Theta}^TG\Psi_6(T) -(\Xi_2^G)^Tl_1^*(T)+(\Xi_1^G)^Tl_2^*(T)+(\Xi_1^G-G)^Tq^*(T) ,\\
&dy_0^*=-(\alpha Q_0\Psi_1- \Theta_1^TQ\Psi_3+F^Ty^*-F^T\hat{y}^*+A_0^Ty_0^* -\Xi_4^Tl_1^*+\Xi_2^Tl_2^*+\Xi_2^Tq^*)dt\\
&\hspace{1cm}+\hat{\beta}_0^*dW_0,\\
&y_0^*(T)=\alpha  G_0\Psi_4(T)-\hat{\Theta}_1^TG\Psi_6(T)-(\Xi_4^G)^Tl_1^*(T) +(\Xi_2^G)^Tl_2^*(T)+(\Xi_2^G)^Tq^*(T),\\
&dl_1^*=(A_0l_1^*+C_0l_2^*-C_0q^*)dt,\quad l_1^*(0)=0,\\
&dl_2^*=[Fl_1^*+(A+C)l_2^*-BR^{-1}B^T\hat{y}^*-Cq^*]dt,\quad l_2^*(0)=0.
\end{aligned}
\right.%$}
\end{equation}
and the decentralized optimal control for the leader
\begin{equation}\label{u0*}
\begin{aligned}
u_0^*=-(\alpha R_0)^{-1}B_0^Ty_0^*.
\end{aligned}
\end{equation}
The final CC system is highly coupled with five forward equations and five backward equations. The existence and uniqueness of \eqref{ccsys2_1} is very important for obtaining the optimal control, however it is very difficult to solve such high-dimensional system. We need to simplify the CC system to a FBSDE using block matrices and these will be discussed in next section.

\section{Well-posedness of the CC system}

Note that in \eqref{ccsys2_1}, the equations of $(\hat{x},\bar{x}_0,k_1,k_2)$ form a coupled FBSDE and $(y^*,q^*,\hat{y}^*,y_0^*,l_1^*,l_2^*)$ form another coupled FBSDE. The two FBSDEs are also fully coupled with each other. Therefore, we try to look at the above FBSDEs in a different way. To this end, we set
%\begin{equation*}%\resizebox{\linewidth}{!}{$
%\mathbb{X}=\left(
%  \begin{array}{c}
%    \hat{x}\\
%    \bar{x}_0 \\
%    q^* \\
%    l_1^* \\
%    l_2^* \\
%  \end{array}
%\right),\
%\mathbb{Y}=\left(
%  \begin{array}{c}
%    y^* \\
%    \hat{y}^* \\
%    y_0^* \\
%    k_1 \\
%    k_2 \\
%  \end{array}
%\right),\
%\mathbb{X}(0)=\left(
%  \begin{array}{c}
%    \hat{\xi} \\
%    \xi_0 \\
%    0 \\
%    0 \\
%    0 \\
%  \end{array}
%\right),
%\end{equation*}
%\begin{equation*}
%\mathbb{Y}(T)=\left(
%  \begin{array}{c}
%    G\Psi_6-Gq^*(T) \\
%    \alpha \hat{\Theta}_0^TG_0\Psi_4-\hat{\Theta}^TG\Psi_6 -(\Xi_2^G)^Tl_1^*(T)+(\Xi_1^G)^Tl_2^*(T)+(\Xi_1^G-G)^Tq^*(T)  \\
%    \alpha  G_0\Psi_4-\hat{\Theta}_1^TG\Psi_6-(\Xi_4^G)^Tl_1^*(T) +(\Xi_2^G)^Tl_2^*(T)+(\Xi_2^G)^Tq^*(T)\\
%    \Xi_4^G\bar{x}_0(T)- (\Xi_2^G)^T\hat{x}(T)+\Xi_5^G \\
%    \Xi_1^G\hat{x}(T) -\Xi_2^G\bar{x}_0(T)-\Xi_3^G \\
%  \end{array}
%\right).%$}
%\end{equation*}
\begin{equation*}\resizebox{\linewidth}{!}{$
\mathbb{X}=\left(
  \begin{array}{c}
    \hat{x}\\
    \bar{x}_0 \\
    q^* \\
    l_1^* \\
    l_2^* \\
  \end{array}
\right),\
\mathbb{Y}=\left(
  \begin{array}{c}
    y^* \\
    \hat{y}^* \\
    y_0^* \\
    k_1 \\
    k_2 \\
  \end{array}
\right),\
\mathbb{X}(0)=\left(
  \begin{array}{c}
    \hat{\xi} \\
    \xi_0 \\
    0 \\
    0 \\
    0 \\
  \end{array}
\right),
%\end{equation*}
%\begin{equation*}
\mathbb{Y}(T)=\left(
  \begin{array}{c}
    G\Psi_6-Gq^*(T) \\
    \alpha \hat{\Theta}_0^TG_0\Psi_4-\hat{\Theta}^TG\Psi_6 -(\Xi_2^G)^Tl_1^*(T)+(\Xi_1^G)^Tl_2^*(T)+(\Xi_1^G-G)^Tq^*(T)  \\
    \alpha  G_0\Psi_4-\hat{\Theta}_1^TG\Psi_6-(\Xi_4^G)^Tl_1^*(T) +(\Xi_2^G)^Tl_2^*(T)+(\Xi_2^G)^Tq^*(T)\\
    \Xi_4^G\bar{x}_0(T)- (\Xi_2^G)^T\hat{x}(T)+\Xi_5^G \\
    \Xi_1^G\hat{x}(T) -\Xi_2^G\bar{x}_0(T)-\Xi_3^G \\
  \end{array}
\right).$}
\end{equation*}
Then \eqref{ccsys2_1} is equivalent to
\begin{equation}\label{fbsde_1}
\left\{
\begin{aligned}
&d\mathbb{X}=[\mathbb{A}\mathbb{X}+\mathbb{B}\mathbb{Y}+b]dt +\mathbb{D}dW_0,\quad \mathbb{X}(0)=(\hat{\xi}^T\quad \xi_0^T\quad 0\quad 0\quad 0)^T,\\
&d\mathbb{Y}=[\hat{\mathbb{A}}\mathbb{X}+ \hat{\mathbb{B}}\mathbb{Y}+\hat{b}]dt+\hat{\mathbb{D}}dW_0,\quad \mathbb{Y}(T)=\mathbb{G}\mathbb{X}(T)+g,
\end{aligned}
\right.
\end{equation}
with
\begin{equation*}\resizebox{\linewidth}{!}{$
\mathbb{A}=\left(
  \begin{array}{ccccc}
    A+C & F & 0 & 0 & 0\\
    C_0 & A_0 & 0 & 0 & 0\\
    0 & 0 & A & 0 & 0\\
    0 & 0 & -C_0 & A_0 & C_0\\
    0 & 0 & -C & F & (A+C)\\
  \end{array}
\right),\
%\end{equation*}
%\begin{equation*}
\mathbb{B}=\left(
  \begin{array}{ccccc}
    0 & 0 & 0 & 0 & -BR^{-1}B^T\\
    0 & 0 & -B_0(\alpha R_0)^{-1}B_0^T & 0 & 0\\
    BR^{-1}B^T & 0 & 0 & 0 & -BR^{-1}B^T\\
    0 & 0 & 0 & 0 & 0\\
    0 & -BR^{-1}B^T & 0 & 0 & 0\\
  \end{array}
\right),$}
\end{equation*}
\begin{equation*}%\resizebox{\linewidth}{!}{$
b=\left(
  \begin{array}{c}
    0 \\
    0 \\
    0 \\
    0 \\
    0 \\
  \end{array}
\right),\quad
\mathbb{D}=\left(
  \begin{array}{c}
    0 \\
    D_0 \\
    0 \\
    0 \\
    0 \\
  \end{array}
\right), \quad
\hat{\mathbb{A}}=\left(
  \begin{array}{ccccc}
    -Q(I-\Theta) & Q\Theta_1 & Q^T & 0 & 0\\
    \Xi_1-Q(I-\Theta) & -\Xi_2+Q\Theta_1 & -(\Xi_1-Q)^T & \Xi_2 & -\Xi_1^T\\
    \Xi_2^T & -\Xi_4 & -\Xi_2^T & \Xi_4^T & -\Xi_2^T\\
    \Xi_2^T & -\Xi_4 & 0 & 0 & 0\\
    -\Xi_1 & \Xi_2 & 0 & 0 & 0\\
  \end{array}
\right),%$}
\end{equation*}
\begin{equation*}%\resizebox{\linewidth}{!}{$
\hat{b}=\left(
  \begin{array}{c}
    Q\eta \\
    -\Xi_3+Q\eta \\
    -\Xi_5 \\
    -\Xi_5 \\
    \Xi_3 \\
  \end{array}
\right), \quad
\hat{\mathbb{D}}=\left(
  \begin{array}{c}
    \beta^* \\
    \hat{\beta}^* \\
    \hat{\beta}_0^* \\
    \beta_1 \\
    \beta_2 \\
  \end{array}
\right), \quad
%\end{equation*}
%\begin{equation*}
\hat{\mathbb{B}}=\left(
  \begin{array}{ccccc}
    -A^T & 0 & 0 & 0 & 0\\
    C^T & -(A+C)^T & C_0^T & 0 & 0\\
    -F^T & F^T & -A_0^T & 0 & 0\\
    0 & 0 & 0 & -A_0^T & -F^T\\
    0 & 0 & 0 & -C_0^T & -(A+C)^T\\
  \end{array}
\right),%$}
\end{equation*}
\begin{equation*}%\resizebox{\linewidth}{!}{$
g=\left(
  \begin{array}{c}
    -G\hat{\eta} \\
    \Xi_3^G-G\hat{\eta} \\
    \Xi_5^G \\
    \Xi_5^G \\
    -\Xi_3^G \\
  \end{array}
\right), \quad
\mathbb{G}=\left(
  \begin{array}{ccccc}
    G(I-\hat{\Theta}) & -G\hat{\Theta}_1 & -G & 0 & 0\\
    -\Xi_1^G+G(I-\hat{\Theta}) & \Xi_2^G-G\hat{\Theta}_1 & (\Xi_1^G-G)^T & -(\Xi_2^G)^T & (\Xi_1^G)^T \\
    -(\Xi_2^G)^T & \Xi_4^G & (\Xi_2^G)^T & -(\Xi_4^G)^T & (\Xi_2^G)^T \\
    -(\Xi_2^G)^T& \Xi_4^G & 0 & 0 & 0\\
    \Xi_1^G & -\Xi_2^G & 0 & 0 & 0\\
  \end{array}
\right).%$}
\end{equation*}
Denote
\begin{equation*}
\begin{aligned}
\bar{\mathbb{A}}=\left(
  \begin{array}{cc}
    \mathbb{A}+\mathbb{B}\mathbb{G} & \mathbb{B} \\
    \hat{\mathbb{A}}-\mathbb{G}\mathbb{A}+\hat{\mathbb{B}}\mathbb{G} -\mathbb{G}\hat{\mathbb{B}}\mathbb{G} & \hat{\mathbb{B}}-\mathbb{G}\mathbb{B} \\
  \end{array}
\right),\quad
\bar{b}=\left(
  \begin{array}{c}
    b  \\
    \hat{b}-\mathbb{G}b  \\
  \end{array}
\right),\quad
\bar{\mathbb{D}}=\left(
  \begin{array}{c}
    \mathbb{D}  \\
    \hat{\mathbb{D}}  \\
  \end{array}
\right),\quad  \bar{\mathbb{Y}}=\mathbb{Y}-\mathbb{G}\mathbb{X}.
\end{aligned}
\end{equation*}
Then \eqref{fbsde_1} can be rewritten as:
\begin{equation}\label{fbsde_2}
\left\{
\begin{aligned}
&d\left(
  \begin{array}{c}
    \mathbb{X}  \\
    \bar{\mathbb{Y}}  \\
  \end{array}
\right)=\left\{\bar{\mathbb{A}}\left(
  \begin{array}{c}
    \mathbb{X}  \\
    \bar{\mathbb{Y}}  \\
  \end{array}
\right)+\bar{b}\right\}dt+\bar{\mathbb{D}}dW_0,\\
&\mathbb{X}(0)=(\hat{\xi}^T\quad \xi_0^T\quad 0\quad 0\quad 0)^T,\quad \bar{\mathbb{Y}}(T)=g.
\end{aligned}
\right.
\end{equation}
This is a fully coupled FBSDE. From the Theorem 3.7 of Chapter 2 in \cite{my1999}, the FBSDE \eqref{fbsde_2} is solvable for all $g\in L_{\mathbb{F}}^2(\Omega ;\mathbb{R}^{5n})$ if and only if the following condition holds:
\begin{equation}\label{fbsde_ricatti_solve_1}
\begin{aligned}
\det\left\{(0,I)e^{\bar{\mathbb{A}}t}\left(\begin{array}{c}0\\I\end{array} \right)\right\}> 0, \quad \forall \ t\in [0,T].
\end{aligned}
\end{equation}
In the case, \eqref{fbsde_1} admits an unique solution for any given $g\in L_{\mathbb{F}}^2(\Omega ;\mathbb{R}^{5n})$.

Under the condition \eqref{fbsde_ricatti_solve_1}, we may decouple the FBSDE \eqref{fbsde_2} by  %It can be verified that $\mathcal{X}$ and $\mathcal{Y}$ satisfy the following equation
$$\bar{\mathbb{Y}}=\mathbb{K}\mathbb{X}+\kappa, \quad \ t\in [0,T],$$
where $\mathbb{K}\in C^1([0,T];\mathbb{S}^{5n})$ is a solution of the following Ricatti equation
\begin{equation*}
%\left\{
\begin{aligned}
&\dot{\mathbb{K}}+\mathbb{K}(\mathbb{A}+\mathbb{B}\mathbb{G}) +\mathbb{K}\mathbb{B}\mathbb{K} -(\hat{\mathbb{B}}-\mathbb{G}\mathbb{B})\mathbb{K} -(\hat{\mathbb{A}}-\mathbb{G}\mathbb{A}+\hat{\mathbb{B}}\mathbb{G} -\mathbb{G}\hat{\mathbb{B}}\mathbb{G})=0,  \ t\in [0,T],\quad
\mathbb{K}(T)=0,
\end{aligned}
%\right.
\end{equation*}
and $\kappa\in C^1([0,T];\mathbb{R}^{5n})$ satisfies
\begin{equation}\label{ricattikappa}
%\left\{
\begin{aligned}
&\dot{\kappa}+(\mathbb{K}\mathbb{B}-(\hat{\mathbb{B}}-\mathbb{G}\mathbb{B}))\kappa+ \mathbb{K}b-(\hat{b}-\mathbb{G}b)=0, \ t\in [0,T],\quad
\kappa(T)=g.
\end{aligned}
%\right.
\end{equation}
By the Theorem 3.7 and Theorem 4.3 of Chapter 2 in \cite{my1999}, if \eqref{fbsde_ricatti_solve_1} hold, then the Ricatti equation admits a unique solution $\mathbb{K}(\cdot)$ which has the following representation:
\begin{equation}\label{ricattik}
\begin{aligned}
\mathbb{K}=-\bigg[(0,I)e^{\bar{\mathbb{A}}(T-t)}\left(\begin{array}{c}0\\I\end{array} \right)\bigg]^{-1}\bigg[(0,I)e^{\bar{\mathbb{A}}(T-t)}\left(\begin{array}{c}I\\0\end{array} \right)\bigg], \quad  t\in [0,T].
\end{aligned}
\end{equation}
\begin{xmpl}\label{example_5.1}
Consider the system \eqref{fbsde_2} with parameters $A_0=0.1$, $B_0=1$, $C_0=0.01$, $D_0=1$, $A=0.05$, $B=1$, $C=0.05$, $D=1$, $F=0.3$, $\Theta_0=1$, $Q_0=1$, $R_0=10$, $G_0=0$, $\Theta=0.1$, $\Theta_1=1$, $Q=0.9$, $R=15$, $G=0$, $\alpha=1.02$, $T=12$, $\eta_0=\eta=0$. Then, we have
\begin{equation*}
\begin{aligned}
&\mathbb{A}=\small{\left(
  \begin{array}{ccccc}
    0.10 & 0.30 & 0 & 0 & 0\\
    0.01 & 0.10 & 0 & 0 & 0\\
    0 & 0 & 0.05 & 0 & 0\\
    0 & 0 & -0.01 & 0.10 & 0.01\\
    0 & 0 & -0.05 & 0.30 & 0.10\\
  \end{array}
\right)},
\mathbb{B}=\small{\left(
  \begin{array}{ccccc}
    0 & 0 & 0 & 0 & -0.0667\\
    0 & 0 & -0.0980 & 0 & 0\\
    0.0667 & 0 & 0 & 0 & -0.0667\\
    0 & 0 & 0 & 0 & 0\\
    0 & -0.0667 & 0 & 0 & 0\\
  \end{array}
\right)}, \\
&\hat{\mathbb{A}}=\small{\left(
  \begin{array}{ccccc}
    -0.81 & 0.90 & 0.90 & 0 & 0\\
    0.939 & -0.93 & -2.649 & 1.83 & -1.749\\
    1.83 & -1.92 & -1.83 & 1.92 & -1.83\\
    1.83 & -1.92 & 0 & 0 & 0\\
    -1.749 & 1.83 & 0 & 0 & 0\\
  \end{array}
\right)},
\hat{\mathbb{B}}=\small{\left(
  \begin{array}{ccccc}
    -0.05 & 0 & 0 & 0 & 0\\
    0.05 & -0.10 & 0.01 & 0 & 0\\
    -0.30 & 0.30 & -0.10 & 0 & 0\\
    0 & 0 & 0 & -0.10 & -0.30\\
    0 & 0 & 0 & -0.01 & -0.10\\
  \end{array}
\right)}.
\end{aligned}
\end{equation*}
Hence, according to the simulation through Matlab software, for any $t\in[0,T]$, we obtain
\begin{equation*}
\begin{aligned}
\bar{\mathbb{A}}=\left(
  \begin{array}{cc}
    \mathbb{A}+\mathbb{B}\mathbb{G} & \mathbb{B} \\
    \hat{\mathbb{A}}-\mathbb{G}\mathbb{A}+\hat{\mathbb{B}}\mathbb{G} -\mathbb{G}\hat{\mathbb{B}}\mathbb{G} & \hat{\mathbb{B}}-\mathbb{G}\mathbb{B} \\
  \end{array}
\right),\quad \quad
\det\left\{(0,I)e^{\bar{\mathbb{A}}t}\left(\begin{array}{c}0\\I\end{array} \right)\right\}> 0,
\end{aligned}
\end{equation*}
(e.g. for $t=6$,  $\det\left\{(0,I)e^{\bar{\mathbb{A}}t}\left(\begin{array}{c}0\\I\end{array} \right)\right\}=12.7053> 0$). By the argument above, \text{\rm FBSDE} \eqref{fbsde_1} is solvable.
\end{xmpl}

For further analysis, we make the following assumption:\\
%\begin{assumption}
\text{\rm(\textbf{A4})} The equation \eqref{fbsde_2} has a unique solution and the solution $(\mathbb{X},\bar{\mathbb{Y}},\bar{\mathbb{D}})$ belongs to $\mathcal{M}[0,T]$.
%\end{assumption}

For the following equation
\begin{equation}\label{barxi_barpi}%\resizebox{\linewidth}{!}{$
\left\{
\begin{aligned}
&d\bar{x}_i=[A\bar{x}_i-BR^{-1}B^Tp_i+C\hat{x}+F\bar{x}_0]dt+DdW_i, \ x_i(0)=\xi_{i},\ i=1,2,\cdots,N,\\
&dp_i=-[A^Tp_i+Q\bar{x}_i+\chi_1]dt+\zeta_0dW_0+\zeta_idW_i,
\ p_i(T)=Gx_i(T)+\chi_2,\\
\end{aligned}
\right.%$}
\end{equation}
where $\chi_1$ and $\chi_2$ are related to $\bar{u}_0$. We let $p_i=\bar{P}\bar{x}_i+\bar{\varphi},\ t\in [0,T]$, where $\bar{P}\in C^1([0,T];\mathbb{S}^{n})$ is a solution of the following Ricatti equation and $\bar{\varphi}\in C^1([0,T];\mathbb{R}^{n})$ satisfies
%\begin{equation*}
%%\left\{
%\begin{aligned}
%&\dot{\bar{P}}+\bar{P}A-\bar{P}BR^{-1}B^T\bar{P}+A^T\bar{P}+Q=0, \ t\in [0,T],\ \bar{P}(T)=G,
%\end{aligned}
%%\right.
%\end{equation*}
%
%\begin{equation*}
%%\left\{
%\begin{aligned}
%&\dot{\bar{\varphi}}+(A^T-\bar{P}BR^{-1}B^T)\bar{\varphi} +\chi_1+\bar{P}C\hat{x}+\bar{P}F\bar{x}_0=0,\ t\in [0,T],\
%\bar{\varphi}(T)=\chi_2.\\
%\end{aligned}
%%\right.
%\end{equation*}
\begin{equation*}
\left\{
\begin{aligned}
&\dot{\bar{P}}+\bar{P}A-\bar{P}BR^{-1}B^T\bar{P}+A^T\bar{P}+Q=0, \ t\in [0,T],\ \bar{P}(T)=G,\\
%\end{aligned}
%%\right.
%\end{equation*}
%and $\bar{\varphi}\in C^1([0,T];\mathbb{R}^{n})$ satisfies
%\begin{equation*}
%%\left\{
%\begin{aligned}
&\dot{\bar{\varphi}}+(A^T-\bar{P}BR^{-1}B^T)\bar{\varphi} +\chi_1+\bar{P}C\hat{x}+\bar{P}F\bar{x}_0=0,\ t\in [0,T],\
\bar{\varphi}(T)=\chi_2.\\
\end{aligned}
\right.
\end{equation*}
Since the Ricatti equation is standard, it has a unique solution. Hence, the FBSDE \eqref{barxi_barpi} is uniquely solvable and the solution belongs to $\mathcal{M}[0,T]$.

%The optimal control of the leader and the followers are %$\bar{u}_i=-R^{-1}B^Tp_i$ and $u_0^*=-(\alpha R_0)^{-1}B_0^Ty_0^*$
%\begin{equation}\label{u0*}
%%\left\{
%\begin{aligned}
%&u_0^*=-(\alpha R_0)^{-1}B_0^Ty_0^*,\\
%\end{aligned}
%%\right.
%\end{equation}
%where $y_0^*$ %, $\hat{x}$, $x_0$, $y^*$, $\hat{y}^*$, $l_1^*$ and $l_2^*$ are
%is determined by the CC system \eqref{ccsys2_1}, and
%\begin{equation}\label{ui*}
%\left\{
%\begin{aligned}
%&u^*_i=-R^{-1}B^Tp_i,\quad i=1,2,\cdots,N,\\
%&d\bar{x}_0=[A_0\bar{x}_0+B_0u_0^*+C_0\hat{x}]dt+D_0dW_0,\quad \bar{x}_0(0)=\xi_0,\\
%&d\bar{x}_i=[A\bar{x}_i-BR^{-1}B^Tp_i+C\hat{x}+F\bar{x}_0]dt+DdW_i,\quad x_i(0)=\xi_{i},\\
%&dp_i=-[A^Tp_i+Q\bar{x}_i+\chi_1]dt+\zeta_0dW_0+\zeta_idW_i,
%\quad p_i(T)=Gx_i(T)+\chi_2,\\
%\end{aligned}
%\right.
%\end{equation}
%where $\hat{x}$, %$\bar{x}_0$,
%$k_1$ and $k_2$ are determined by the CC system \eqref{ccsys2_1}.

\section{Asymptotically social optimality}

In this section, we discuss that if the leader announces $u^*_0$ obtained in \eqref{u0*} to the $N$ followers, then the set of the optimal decentralized controls for the leader and the followers will constitutes an approximated Stackelberg equilibrium. First,
for the open-loop decentralized strategy $(u^*_0,u^*)$ in \eqref{u0*} and \eqref{bar_u_i}, %given as $\tilde{u}_i  =  \Theta_1\tilde{x}_i + \Theta_2$, where
we have the realized decentralized state $x_0^*$ and $x_i^*$, satisfies
\begin{equation}\label{realized state}
\left\{
\begin{aligned}
&dx_0^*(t)=[A_0x_0^*(t)-B_0(\alpha R_0)^{-1}B_0^Ty_0^*(t)+C_0(x^*)^{(N)}(t)]dt +D_0dW_0(t),\\
&dx^*_i(t)=[Ax^*_i(t)-BR^{-1}B^Tp_i(t)+C(x^*)^{(N)}(t) +Fx_0^*(t)]dt+DdW_i(t),\\
&x_0^*(0)=\xi_0, \quad \quad \quad x^*_i(0)=\xi_{i}, \quad i=1,2,\cdots,N,\\
\end{aligned}
\right.
\end{equation}
where $y_0^*, p_i$ satisfy \eqref{ccsys2_1} and \eqref{barxi_barpi}, respectively. Then, by \cite{bo1999} and \cite{mb2018}, we give the definition of the Stackelberg equilibrium.
\begin{dfntn}
  A set of control laws $\mathcal{M}(\check{u}_0)\in \mathcal{U}$ has asymptotic social optimality if \\
  \begin{equation*}
  %\bigg|\bigg(\frac{1}{N}\bigg)J_{soc}^{(N)}(\bar{u}_0,\bar{u}(\bar{u}_0))-\inf_{u_0\in\mathcal{U}_0,u\in \mathcal{U}}\bigg(\frac{1}{N}\bigg)J_{soc}^{(N)}(u_0,u) \bigg|< \varepsilon
  \bigg|\frac{1}{N}\mathcal{J}_{soc}^{(N)}(\check{u}_0; \mathcal{M}(\check{u}_0))-\frac{1}{N}\inf_{(\check{u}_0,\check{u})\in \mathcal{U}_c}\mathcal{J}_{soc}^{(N)}(\check{u}_0; \check{u})\bigg|=O(\frac{1}{\sqrt{N}}),\quad \quad %\lim\limits_{N\rightarrow+\infty}o(1)=0,
\end{equation*}
where $\mathcal{M}$ is a mapping and $\mathcal{M}:\mathcal{U}_0\rightarrow \mathcal{U}$.  $\mathcal{U}_c$ is defined in Section 2 as a set of centralized information-based control.
\end{dfntn}
\begin{dfntn}
A set of control laws $(u_0^*,u^*)\in \mathcal{U}_0\times \mathcal{U}$, where $u^*=\mathcal{M}(u_0^*)$, is an Stackelberg equilibrium with respect to $\mathcal{J}_{soc}^{(N)}(u_0,u)$ if the following two properties hold:
\begin{enumerate}
\item $\mathcal{M}(\check{u}_0)$ has asymptotic social optimality under $\check{u}_0$.
\item The following equation is satisfied
\begin{equation*}
  \bigg|\frac{1}{N}\mathcal{J}_{soc}^{(N)}(u_0^*; \mathcal{M}(u_0^*))-\frac{1}{N}\inf_{\check{u}_0\in \mathcal{U}_c}\mathcal{J}_{soc}^{(N)}(\check{u}_0; \mathcal{M}(\check{u}_0))\bigg|=O(\frac{1}{\sqrt{N}}).
\end{equation*}
\end{enumerate}
\end{dfntn}

We first need to introduce some lemmas before proving the Stackelberg equilibrium.
\begin{lmm}\label{lma_6.1}
Assume that \text{\rm(\textbf{A1})-(\textbf{A4})} hold. %and the CC systems \eqref{ccsys1_1} and \eqref{ccsys2_1} admit solutions.
Then
   \begin{equation*}
  \mathbb{E}\int_{0}^{T}\|(x^*)^{(N)}-\hat{x}\|^2dt+ \mathbb{E}\int_{0}^{T}\|p^{(N)}-\hat{p}\|^2dt+ \mathbb{E}\int_{0}^{T}\|x^*_0-\bar{x}_0\|^2dt=O(\frac{1}{N}).
  \end{equation*}
\end{lmm}
\begin{proof}
%{\it Proof }
See Appendix A.
%\qed
\end{proof}
\begin{lmm}\label{lma_6.2}
Assume that \text{\rm(\textbf{A1})-(\textbf{A4})} hold. There exists a constant $K$, which is independent of $N$, such that
$$\mathcal{J}_{soc}^{(N)}(u_0^*; u^*)\leq N K.$$
\end{lmm}
\begin{proof}
%{\it Proof }
See Appendix B.
%\qed
\end{proof}
\begin{prpstn}\label{prop_6.1}
Assume that \text{\rm(\textbf{A1})-(\textbf{A4})} hold. For all $(\check{u}_0;\check{u})\in \mathcal{U}_c$, there exists a constant $K$, which is independent of $N$, such that
\begin{equation*}
  \alpha N\|\check{u}_0\|_{L^2}^2+\|\check{u}\|_{L^2}^2\leq NK.
\end{equation*}
\end{prpstn}
\begin{proof}
%{\it Proof }
By Lemma \ref{lma_6.2}, we have
\begin{equation*}
\mathbb{E}\int_{0}^{T}\alpha N\|\check{u}_0\|^2+\|\check{u}\|^2dt\leq \inf_{(\check{u}_0;\check{u})}\mathcal{J}_{soc}^{(N)}(\check{u}_0; \check{u})\leq \mathcal{J}_{soc}^{(N)}(u_0^*; u^*)\leq NK,\quad (\alpha>0).
\end{equation*}
Therefore, $\alpha N\|\check{u}_0\|_{L^2}^2+\|\check{u}\|_{L^2}^2\leq NK$, where $K$ is independent of $N$. The proposition follows.
%\qed
\end{proof}
The following two propositions will give the rigorous proofs for the approximations in Section 3.
\begin{prpstn}\label{prop_6.2}
Assume that \text{\rm(\textbf{A1})-(\textbf{A4})} hold. Then, for \eqref{variation_soc}, $\mathbb{E}\sup_{0\leq t\leq T}\|\delta x_0\|^2 = O(\frac{1}{N^2})$, $\mathbb{E}\sup_{0\leq t\leq T}\|\delta x^{(N)}\|^2 = O(\frac{1}{N^2})$ and
$\langle \Theta^{T}Q(\bar{x}_i-\Theta \bar{x}^{(N)}-\Theta_1\bar{x}_0-\eta),\delta x^{(N)}\rangle+ \langle\Theta_1^{T}Q(\bar{x}_i-\Theta\bar{x}^{(N)} -\Theta_1\bar{x}_0-\eta),\delta x_0\rangle+\langle \hat{\Theta}^{T}G(\bar{x}_i(T) -\hat{\Theta}\bar{x}^{(N)}(T) -\hat{\Theta}_1\bar{x}_0(T)-\hat{\eta}),\delta x^{(N)}(T)\rangle+\langle \hat{\Theta}_1^{T}G(\bar{x}_i(T) -\hat{\Theta}\bar{x}^{(N)}(T) -\hat{\Theta}_1\bar{x}_0(T)-\hat{\eta}),\delta x_0(T)\rangle=o(1)$.
\end{prpstn}
\begin{proof}
%{\it Proof }
 See Appendix C.
 %\qed
\end{proof}
\begin{prpstn}\label{prop_6.3}
Assume that \text{\rm(\textbf{A1})-(\textbf{A4})} hold. Then, $N\delta x_j$, $N\delta x_0$, $N\delta x_j$ converge to $\sum_{j\neq i}\delta x_j$, $\delta x_0^{\dagger}$, $\delta x^{\dagger}$ such that %$N\delta x_0$, $\sum_{j\neq i}\delta x_j$ and $N\delta x_j$ converge to $\delta x_0^{\dagger}$, $N\delta x_j$ and $\delta x^{\dagger}$, respectively such that
\begin{equation*}
%\left\{
  \begin{aligned}
  &\mathbb{E}\int_{0}^{T}\|N\delta x_j-\sum_{j\neq i}\delta x_j\|^2=O(\frac{1}{N^2}),\
  \mathbb{E}\int_{0}^{T}\|N\delta x_0-\delta x_0^{\dagger}\|^2=O(\frac{1}{N^2}),\ \mathbb{E}\int_{0}^{T}\|N\delta x_j-\delta x^{\dagger}\|^2=O(\frac{1}{N^2}).
  \end{aligned}
 % \right.
\end{equation*}
\end{prpstn}
\begin{proof}
%{\it Proof }
See Appendix C.
%\qed
\end{proof}

By the lemmas and propositions we discussed above, we give the main result.
\begin{thrm}
Assume that \text{\rm(\textbf{A1})-(\textbf{A4})} hold. %and the CC systems \eqref{ccsys1_1} and \eqref{ccsys2_1} admit solutions.
Then $(u_0^*, u^*)$ given in \eqref{u0*} and \eqref{bar_u_i} is an Stackelberg equilibrium with respect to the social cost functional.
\end{thrm}
\begin{proof}
%{\it Proof }
For $(\check{u}_0;\check{u})\in \mathcal{U}_c$, let
\begin{equation*}%\resizebox{\linewidth}{!}{$
\begin{aligned}
&\frac{1}{N}\mathcal{J}_{soc}^{(N)}(u_0^*; u^*)- \frac{1}{N}\mathcal{J}_{soc}^{(N)}(\check{u}_0; \check{u})
%=&\frac{1}{N}\mathcal{J}_{soc}^{(N)}(u_0^*; u^*)- \frac{1}{N}\mathcal{J}_{soc}^{(N)}(\check{u}_0; \mathcal{M}(\check{u}_0))+\frac{1}{N}\mathcal{J}_{soc}^{(N)}(\check{u}_0; \mathcal{M}(\check{u}_0))- \frac{1}{N}\mathcal{J}_{soc}^{(N)}(\check{u}_0; \check{u})\\
=\frac{1}{N}\mathcal{J}_{soc}^{(N)}(u_0^*; \mathcal{M}(u_0^*))- \frac{1}{N}\mathcal{J}_{soc}^{(N)}(\check{u}_0; \mathcal{M}(\check{u}_0))\\
&+\frac{1}{N}\mathcal{J}_{soc}^{(N)}(\check{u}_0; \mathcal{M}(\check{u}_0))- \frac{1}{N}\mathcal{J}_{soc}^{(N)}(\check{u}_0; \check{u})
:= \Delta_1+\Delta_2,
\end{aligned}%$}
\end{equation*}
where $\Delta_1=\frac{1}{N}\mathcal{J}_{soc}^{(N)}(u_0^*; \mathcal{M}(u_0^*))- \frac{1}{N}\mathcal{J}_{soc}^{(N)}(\check{u}_0; \mathcal{M}(\check{u}_0))$, $\Delta_2=\frac{1}{N}\mathcal{J}_{soc}^{(N)}(\check{u}_0; \mathcal{M}(\check{u}_0))- \frac{1}{N}\mathcal{J}_{soc}^{(N)}(\check{u}_0; \check{u})$.
%with
%\begin{equation*}
%\left\{
%\begin{aligned}
%&\Delta_1=\frac{1}{N}\mathcal{J}_{soc}^{(N)}(u_0^*; \mathcal{M}(u_0^*))- \frac{1}{N}\mathcal{J}_{soc}^{(N)}(\check{u}_0; \mathcal{M}(\check{u}_0)),\\
%&\Delta_2=\frac{1}{N}\mathcal{J}_{soc}^{(N)}(\check{u}_0; \mathcal{M}(\check{u}_0))- \frac{1}{N}\mathcal{J}_{soc}^{(N)}(\check{u}_0; \check{u}).
%\end{aligned}
%\right.
%\end{equation*}
Since $\check{u}_0$ is fixed, by following the standard method in Huang \emph{et al} \cite{hcm2012}, we obtain $\|\Delta_2\|^2\leq c(\|\check{u}_0\|_{L^2}^2)\frac{1}{N}$.   %is $\varepsilon$-social optimality for $\check{u}_0$ parameterized process.
Specifically, we denote $\grave{x}_i %=\mathcal{M}(\check{x}_i)
$ as the state of the $i$th follower when its control is $\mathcal{M}_i(\check{u}_0)$, thus $\grave{x}_i$ is equivalent to $\bar{x}_i$ in Section 3.
Let
\begin{equation*}
\left\{
\begin{aligned}
&\tilde{u}_0=\check{u}_0-\check{u}_0=0, \quad \tilde{u}=\check{u}-\mathcal{M}(\check{u}_0), \quad\tilde{u}_i=\check{u}_i-\mathcal{M}_i(\check{u}_0),\\
&\tilde{x}_0=\check{x}_0-\grave{x}_0, \quad \tilde{x}_i=\check{x}_i-\grave{x}_i.\\
\end{aligned}
\right.
\end{equation*}
Then we have
\begin{equation*}
\begin{aligned}
&\mathcal{J}_{soc}^{(N)}(\check{u}_0; \check{u})=\alpha N\mathcal{J}_0(\check{u}_0;\check{u})+\sum_{i=1}^{N}\mathcal{J}_i(\check{u}_0;\check{u})\\
=&\alpha N\mathcal{J}_0(\check{u}_0;\mathcal{M}(\check{u}_0))+\alpha NH_0+\alpha NI_0+\sum_{i=1}^{N}\mathcal{J}_i(\check{u}_0;\mathcal{M}(\check{u}_0))
+\sum_{i=1}^{N}H_i+\sum_{i=1}^{N}I_i,
\end{aligned}
\end{equation*}
where
\begin{equation*}%\resizebox{\linewidth}{!}{$
\begin{aligned}
&\mathcal{J}_0(\check{u}_0;\mathcal{M}(\check{u}_0))=\mathbb{E}\bigg\{\int_{0}^{T}\|\grave{x}_0- \Theta_0\grave{x}^{(N)}-\eta_0\|_{Q_0}^{2}+\|\check{u}_0\|_{R_0}^{2}dt+ \|\grave{x}_0(T)- \hat{\Theta}_0\grave{x}^{(N)}(T)-\hat{\eta}_0\|_{G_0}^{2}\bigg\},\\
%\end{aligned}$}
%\end{equation*}
%\begin{equation*}\resizebox{\linewidth}{!}{$
%\begin{aligned}
&H_0=\mathbb{E}\bigg\{\int_{0}^{T}\|\tilde{x}_0- \Theta_0\tilde{x}^{(N)}\|_{Q_0}^{2}dt+\|\tilde{x}_0(T)- \hat{\Theta}_0\tilde{x}^{(N)}(T)\|_{G_0}^{2}\bigg\},\\
&\mathcal{J}_i(\check{u}_0;\mathcal{M}(\check{u}_0))=\mathbb{E}\bigg\{\int_{0}^{T}\|\grave{x}_i- \Theta \grave{x}^{(N)}-\Theta_1\grave{x}_0-\eta\|_{Q}^{2} +\|\mathcal{M}_i(\check{u}_0)\|_{R}^{2}dt\\
&\hspace{2.4cm}+\|\grave{x}_i(T)- \hat{\Theta} \grave{x}^{(N)}(T)-\hat{\Theta}_1\grave{x}_0(T)-\hat{\eta}\|_{G}^{2}\bigg\},\\
%\end{aligned}$}
%\end{equation*}
%\begin{equation*}\resizebox{\linewidth}{!}{$
%\begin{aligned}
&H_i=\mathbb{E}\bigg\{\int_{0}^{T}\|\tilde{x}_i- \Theta \tilde{x}^{(N)}-\Theta_1\tilde{x}_0\|_{Q}^{2}+\|\tilde{u}_i\|_{R}^2dt+\|\tilde{x}_i(T)- \hat{\Theta} \tilde{x}^{(N)}(T)-\hat{\Theta}_1\tilde{x}_0(T)\|_{G}^{2}\bigg\},\\
%\end{aligned}
%\end{equation*}
%\begin{equation*}
%\begin{aligned}
&I_0=\mathbb{E}\bigg\{\int_{0}^{T}(\grave{x}_0- \Theta_0\grave{x}^{(N)}-\eta_0)^TQ_0(\tilde{x}_0- \Theta_0\tilde{x}^{(N)})dt\\
&\hspace{1cm}+(\grave{x}_0(T)- \hat{\Theta}_0\grave{x}^{(N)}(T)-\hat{\eta}_0)^TG_0(\tilde{x}_0(T)- \hat{\Theta}_0\tilde{x}^{(N)}(T))\bigg\},\\
%\end{aligned}
%\end{equation*}
%\begin{equation*}
%\begin{aligned}
&I_i=\mathbb{E}\bigg\{\int_{0}^{T}(\grave{x}_i- \Theta \grave{x}^{(N)}-\Theta_1\grave{x}_0-\eta)^TQ(\tilde{x}_i- \Theta \tilde{x}^{(N)}-\Theta_1\tilde{x}_0)+\mathcal{M}_i^T(\check{u})R\tilde{u}_idt\\
&\hspace{1cm}+(\grave{x}_i(T)- \hat{\Theta} \grave{x}^{(N)}(T)-\hat{\Theta}_1\grave{x}_0(T)-\hat{\eta})^TG(\tilde{x}_i(T)- \hat{\Theta} \tilde{x}^{(N)}(T)-\hat{\Theta}_1\tilde{x}_0(T))\bigg\}.
\end{aligned}%$}
\end{equation*}
%By straightforward computation
%\begin{equation}\label{I_0}%\resizebox{\linewidth}{!}{$
%\begin{aligned}
%&\alpha NI_0
%%=\alpha N\mathbb{E}\bigg\{\int_{0}^{T}(\grave{x}_0- \Theta_0\grave{x}^{(N)}-\eta_0)^TQ_0(\tilde{x}_0- \Theta_0\tilde{x}^{(N)})dt\\
%%&\hspace{1cm}+(\grave{x}_0(T)- \hat{\Theta}_0\grave{x}^{(N)}(T)-\hat{\eta}_0)^TG_0(\tilde{x}_0(T)- \hat{\Theta}_0\tilde{x}^{(N)}(T))\bigg\}\\
%%&\hspace{1cm}=\mathbb{E}\bigg\{\int_{0}^{T}\alpha N(\grave{x}_0- \Theta_0\grave{x}^{(N)}-\eta_0)^TQ_0\tilde{x}_0-\alpha(\grave{x}_0- \Theta_0\grave{x}^{(N)}-\eta_0)^TQ_0\hat{\Theta}_0\sum_{i=1}^{N}\tilde{x}_idt\\ &\hspace{1cm}+\alpha N(\grave{x}_0(T)- \hat{\Theta}_0\grave{x}^{(N)}(T)-\hat{\eta}_0)^TG_0\tilde{x}_0(T)-\alpha(\grave{x}_0(T)- \hat{\Theta}_0\grave{x}^{(N)}(T)-\hat{\eta}_0)^TG_0\hat{\Theta}_0\sum_{i=1}^{N}\tilde{x}_i(T)\bigg\}\\
%%&\hspace{0.8cm}
%=\mathbb{E}\bigg\{\int_{0}^{T}\alpha N\Psi_1^TQ_0\tilde{x}_0-\alpha N(\Theta_0\upsilon_1)^TQ_0\tilde{x}_0-\alpha \Psi_1^TQ_0\Theta_0\sum_{i=1}^{N}\tilde{x}_i\\
%%\end{aligned}$}
%%\end{equation}
%%\begin{equation}\label{I_i}%\resizebox{\linewidth}{!}{$
%%\begin{aligned}
%&\hspace{1cm}+\alpha (\Theta_0\upsilon_1)^TQ_0\Theta_0\sum_{i=1}^{N}\tilde{x}_idt +\alpha N\Psi_4(T)^TG_0\tilde{x}_0(T)-\alpha N(\hat{\Theta}_0\upsilon_1(T))^TG_0\tilde{x}_0(T)\\
%&\hspace{1cm}-\alpha \Psi_4(T)^TG_0\hat{\Theta}_0\sum_{i=1}^{N}\tilde{x}_i(T)+\alpha (\hat{\Theta}_0\upsilon_1(T))^TG_0\hat{\Theta}_0\sum_{i=1}^{N}\tilde{x}_i(T)\bigg\},\\
%\end{aligned}%$}
%\end{equation}
By straightforward computation
\begin{equation}\label{I_0}%\resizebox{\linewidth}{!}{$
\begin{aligned}
&\alpha NI_0
%=\alpha N\mathbb{E}\bigg\{\int_{0}^{T}(\grave{x}_0- \Theta_0\grave{x}^{(N)}-\eta_0)^TQ_0(\tilde{x}_0- \Theta_0\tilde{x}^{(N)})dt\\
%&\hspace{1cm}+(\grave{x}_0(T)- \hat{\Theta}_0\grave{x}^{(N)}(T)-\hat{\eta}_0)^TG_0(\tilde{x}_0(T)- \hat{\Theta}_0\tilde{x}^{(N)}(T))\bigg\}\\
%&\hspace{1cm}=\mathbb{E}\bigg\{\int_{0}^{T}\alpha N(\grave{x}_0- \Theta_0\grave{x}^{(N)}-\eta_0)^TQ_0\tilde{x}_0-\alpha(\grave{x}_0- \Theta_0\grave{x}^{(N)}-\eta_0)^TQ_0\hat{\Theta}_0\sum_{i=1}^{N}\tilde{x}_idt\\ &\hspace{1cm}+\alpha N(\grave{x}_0(T)- \hat{\Theta}_0\grave{x}^{(N)}(T)-\hat{\eta}_0)^TG_0\tilde{x}_0(T)-\alpha(\grave{x}_0(T)- \hat{\Theta}_0\grave{x}^{(N)}(T)-\hat{\eta}_0)^TG_0\hat{\Theta}_0\sum_{i=1}^{N}\tilde{x}_i(T)\bigg\}\\
%&\hspace{0.8cm}
=\mathbb{E}\bigg\{\int_{0}^{T}\alpha N[\Psi_1^TQ_0-(\Theta_0\upsilon_1)^TQ_0]\tilde{x}_0-\alpha [\Psi_1^TQ_0\Theta_0- (\Theta_0\upsilon_1)^TQ_0\Theta_0]\sum_{i=1}^{N}\tilde{x}_idt \\
%\end{aligned}$}
%\end{equation}
%\begin{equation}\label{I_i}\resizebox{\linewidth}{!}{$
%\begin{aligned}
&+\alpha N[\Psi_4(T)^TG_0-(\hat{\Theta}_0\upsilon_1(T))^TG_0]\tilde{x}_0(T)-\alpha [\Psi_4(T)^TG_0\hat{\Theta}_0- (\hat{\Theta}_0\upsilon_1(T))^TG_0\hat{\Theta}_0]\sum_{i=1}^{N}\tilde{x}_i(T)\bigg\},\\
\end{aligned}%$}
\end{equation}
%\begin{equation}\label{I_i}%\resizebox{\linewidth}{!}{$
%\begin{aligned}
%&\sum_{i=1}^{N} I_i
%%=\sum_{i=1}^{N} \mathbb{E}\bigg\{\int_{0}^{T}(\grave{x}_i- \Theta \grave{x}^{(N)}-\Theta_1\grave{x}_0-\eta)^TQ(\tilde{x}_i- \Theta \tilde{x}^{(N)}-\Theta_1\tilde{x}_0)+\mathcal{M}_i^T(\check{u})R\tilde{u}_idt\\
%%&\hspace{1cm}+(\grave{x}_i(T)- \hat{\Theta} \grave{x}^{(N)}(T)-\hat{\Theta}_1\grave{x}_0(T)-\hat{\eta})^TG(\tilde{x}_i(T)- \hat{\Theta} \tilde{x}^{(N)}(T)-\hat{\Theta}_1\tilde{x}_0(T))\bigg\}\\
%%&\hspace{0.8cm}
%=\mathbb{E}\bigg\{\int_{0}^{T}\sum_{i=1}^{N}(\Psi_2^i)^TQ\tilde{x}_i -(\Theta\upsilon_1)^TQ\sum_{i=1}^{N}\tilde{x}_i -\Psi_3^TQ\Theta\sum_{i=1}^{N}\tilde{x}_i +[(I-\Theta)\upsilon_1]^TQ\Theta\sum_{i=1}^{N}\tilde{x}_i\\
%&\hspace{1cm} -N\Psi_3^TQ\Theta_1\tilde{x}_0+N[(I-\Theta)\upsilon_1]^TQ\Theta_1\tilde{x}_0 +\mathcal{M}_i^T(\check{u})R\tilde{u}_idt +\sum_{i=1}^{N}(\Psi_5^i(T))^TG\tilde{x}_i(T)\\
%&\hspace{1cm} -(\hat{\Theta}\upsilon_1(T))^TG\sum_{i=1}^{N}\tilde{x}_i(T)
%-\Psi_6(T)^TG\hat{\Theta}\sum_{i=1}^{N}\tilde{x}_i(T) +[(I-\hat{\Theta})\upsilon_1(T)]^TG\hat{\Theta}\sum_{i=1}^{N}\tilde{x}_i(T)\\
%%\end{aligned}
%%\end{equation}
%%\begin{equation*}
%%\begin{aligned}
%&\hspace{1cm}
%-N\Psi_6(T)^TG\hat{\Theta}_1\tilde{x}_0(T) +N[(I-\hat{\Theta})\upsilon_1(T)]^TG\hat{\Theta}_1\tilde{x}_0(T)\bigg\}.\\
%\end{aligned}%$}
%\end{equation}
\begin{equation}\label{I_i}%\resizebox{\linewidth}{!}{$
\begin{aligned}
&\sum_{i=1}^{N} I_i
%=\sum_{i=1}^{N} \mathbb{E}\bigg\{\int_{0}^{T}(\grave{x}_i- \Theta \grave{x}^{(N)}-\Theta_1\grave{x}_0-\eta)^TQ(\tilde{x}_i- \Theta \tilde{x}^{(N)}-\Theta_1\tilde{x}_0)+\mathcal{M}_i^T(\check{u})R\tilde{u}_idt\\
%&\hspace{1cm}+(\grave{x}_i(T)- \hat{\Theta} \grave{x}^{(N)}(T)-\hat{\Theta}_1\grave{x}_0(T)-\hat{\eta})^TG(\tilde{x}_i(T)- \hat{\Theta} \tilde{x}^{(N)}(T)-\hat{\Theta}_1\tilde{x}_0(T))\bigg\}\\
%&\hspace{0.8cm}
=\mathbb{E}\bigg\{\int_{0}^{T}\sum_{i=1}^{N}(\Psi_2^i)^TQ\tilde{x}_i -[(\Theta\upsilon_1)^TQ +\Psi_3^TQ\Theta -((I-\Theta)\upsilon_1)^TQ\Theta]\sum_{i=1}^{N}\tilde{x}_i-N[\Psi_3^TQ\Theta_1\\
&-[(I-\Theta)\upsilon_1]^TQ\Theta_1]\tilde{x}_0 +\mathcal{M}_i^T(\check{u})R\tilde{u}_idt +\sum_{i=1}^{N}(\Psi_5^i(T))^TG\tilde{x}_i(T) -[(\hat{\Theta}\upsilon_1(T))^TG+\Psi_6(T)^TG\hat{\Theta}\\
& -((I-\hat{\Theta})\upsilon_1(T))^TG\hat{\Theta}]\sum_{i=1}^{N}\tilde{x}_i(T)
-N[\Psi_6(T)^TG\hat{\Theta}_1 -[(I-\hat{\Theta})\upsilon_1(T)]^TG\hat{\Theta}_1]\tilde{x}_0(T)\bigg\}.\\
\end{aligned}%$}
\end{equation}
where $\upsilon_1=\grave{x}^{(N)}-\hat{x}$. By \eqref{ccsys2_1}, \eqref{barxi_barpi} and
%have following BSDEs
%\begin{equation}\label{grave_p_k}
%\left\{
%\begin{aligned}
%&dk_1=-[\alpha Q_0\Psi_1-\Theta_1^{T}Q\Psi_3+F^Tk_2+A_0^Tk_1]dt+\beta_1dW_0,\quad k_1(T)=\alpha G_0\Psi_4-\hat{\Theta}_1^{T}G\Psi_6,\\
%&dp_i=-[A^Tp_i+Q\grave{x}_i+\chi_1]dt+\zeta_0dW_0+\zeta_idW_i,\quad p_i(T)=G\grave{x}_i(T)+\chi_2.\\
%\end{aligned}
%\right.
%\end{equation}
%Then, using It\^{o} formula, one can obtain following relations:
It\^{o} formula, we obtain following relations:
\begin{equation}
\begin{aligned}
&N\langle k_1(T),\tilde{x}_0(T)\rangle=\langle\alpha N G_0\Psi_4,\tilde{x}_0(T)\rangle-\langle N\hat{\Theta}_1^{T}G\Psi_6,\tilde{x}_0(T)\rangle\\
%&=N\langle dk_1,\tilde{x}_0\rangle+N\langle k_1,d\tilde{x}_0\rangle\\
%&=\mathbb{E}\int_{0}^{T}\langle-\alpha Q_0\Psi_1+\Theta_1^{T}Q\Psi_3-F^Tk_2-A_0^Tk_1,N\tilde{x}_0\rangle+N\langle k_1,A_0\tilde{x}_0+C_0\tilde{x}^{(N)}\rangle dt\\
&=\mathbb{E}\int_{0}^{T}-\langle\alpha NQ_0\Psi_1,\tilde{x}_0\rangle+\langle N\Theta_1^{T}Q\Psi_3,\tilde{x}_0\rangle-\langle k_2,NF\tilde{x}_0\rangle+\langle C_0^Tk_1,\sum_{i=1}^{N}\tilde{x}_i\rangle dt,\\
\end{aligned}
\end{equation}
and
\begin{equation}
\begin{aligned}
&\sum_{i=1}^{N}\langle p_i(T),\tilde{x}_i(T)\rangle
%=\langle G\grave{x}_i(T)+\chi_2,\sum_{i=1}^{N}\tilde{x}_i(T)\rangle\\
%&=\langle G\Psi_5^i(T),\sum_{i=1}^{N}\tilde{x}_i(T)\rangle-\langle\hat{\Theta}^{T}G\Psi_6,\sum_{i=1}^{N}\tilde{x}_i(T)\rangle -\langle\alpha\hat{\Theta}_0^{T}G_0\Psi_4,\sum_{i=1}^{N}\tilde{x}_i(T)\rangle\\
%&=\sum_{i=1}^{N}\langle dp_i,\tilde{x}_i\rangle+\sum_{i=1}^{N}\langle p_i,d\tilde{x}_i\rangle\\
%&=\mathbb{E}\int_{0}^{T}\langle-A^Tp_i-Q\grave{x}_i-\chi_1, \sum_{i=1}^{N}\tilde{x}_i\rangle+\sum_{i=1}^{N}\langle p_i,A\tilde{x}_i+B\tilde{u}_i+C\tilde{x}^{(N)}+F\tilde{x}_0\rangle dt\\
=\mathbb{E}\int_{0}^{T}\langle\Theta^{T}Q\Psi_3,\sum_{i=1}^{N}\tilde{x}_i\rangle -\langle Q\Psi_2^i,\sum_{i=1}^{N}\tilde{x}_i\rangle +\langle\alpha\Theta_0^{T}Q_0\Psi_1,\sum_{i=1}^{N}\tilde{x}_i\rangle\\ &\hspace{0.5cm}-\langle C_0^Tk_1,\sum_{i=1}^{N}\tilde{x}_i\rangle -\langle p^{(N)}-k_2,C\sum_{i=1}^{N}\tilde{x}_i\rangle
+\sum_{i=1}^{N}\langle p_i,B\tilde{u}_i\rangle+\langle p^{(N)},N F\tilde{x}_0\rangle dt.\\
\end{aligned}
\end{equation}
Meanwhile, by \eqref{bar_u_i}, we have
\begin{equation}\label{I_u}
\begin{aligned}
&\sum_{i=1}^{N}\langle\mathcal{M}_i(\check{u}),R\tilde{u}_i\rangle +\sum_{i=1}^{N}\langle p_i,B\tilde{u}_i\rangle
=\sum_{i=1}^{N}\langle R\mathcal{M}_i(\check{u})+B^Tp_i,\tilde{u}_i\rangle=\sum_{i=1}^{N}\langle R(-R^{-1}B^Tp_i)+B^Tp_i,\tilde{u}_i\rangle
%=\sum_{i=1}^{N}\langle -B^Tp_i+B^Tp_i,\tilde{u}_i\rangle
=0.
\end{aligned}
\end{equation}
Combining \eqref{I_0}-\eqref{I_u}, Lemma \ref{lma_6.1} and Lemma \ref{lma_6.2}, it follows that
\begin{equation*}
\begin{aligned}
&\frac{1}{N}(\alpha NI_0+\sum_{i=1}^{N}I_i)
%=\mathbb{E}\bigg\{\int_{0}^{T}\alpha\langle Q_0 \Theta_0\upsilon_1,\tilde{x}_0+\frac{\Theta_0}{N}\sum_{i=1}^{N}\tilde{x}_i\rangle -\langle Q\Theta\upsilon_1,\sum_{i=1}^{N}\tilde{x}_i\rangle\\ &\hspace{3cm} +\langle Q(I-\Theta)\upsilon_1,\Theta\sum_{i=1}^{N}\tilde{x}_i +\Theta_1\tilde{x}_0\rangle+\langle p^{(N)}-k_2,NF\tilde{x}_0-C\sum_{i=1}^{N}\tilde{x}_i\rangle dt\\ &\hspace{3cm}-\alpha\langle G_0\hat{\Theta}_0\upsilon_1(T),\tilde{x}_0(T)+\frac{\hat{\Theta}_0}{N}\sum_{i=1}^{N}\tilde{x}_i(T)\rangle -\langle G\hat{\Theta}\upsilon_1(T),\sum_{i=1}^{N}\tilde{x}_i(T)\rangle\\
%&\hspace{3cm}+\langle G(I-\hat{\Theta})\upsilon_1(T),\hat{\Theta}\sum_{i=1}^{N}\tilde{x}_i(T) +\hat{\Theta}_1\tilde{x}_0(T)\rangle\bigg\}
%=\sqrt{O(\varepsilon_1^2+ \varepsilon_2^2+\varepsilon_3^2)\|\check{u}_0\|_{L^2}}
=O(\frac{1}{\sqrt{N}}).
\end{aligned}
\end{equation*}
Moreover, $\frac{1}{N}(\alpha NH_0+\sum_{i=1}^{N}H_i)\geq 0$. Thus, we have
\begin{equation}\label{delta_2}
%\left\{
\begin{aligned}
\Delta_2=\frac{1}{N}\mathcal{J}_{soc}^{(N)}(\check{u}_0; \mathcal{M}(\check{u}_0))- \frac{1}{N}\mathcal{J}_{soc}^{(N)}(\check{u}_0; \check{u})\leq c(\|\check{u}_0\|_{L^2})\frac{1}{\sqrt{N}}.
\end{aligned}
%\right.
\end{equation}
%We let $u_0^{(N)}$ be the centralized optimal control, $\lim\limits_{N\rightarrow+\infty}u_0^{(N)}=u_0^*$ be the decentralized social optimality.
For $\Delta_1$, we decompose it as follows:
\begin{equation*}
%\left\{
\begin{aligned}
&\Delta_1=\frac{1}{N}\mathcal{J}_{soc}^{(N)}(u_0^*; \mathcal{M}(u_0^*))- \frac{1}{N}\mathcal{J}_{soc}^{(N)}(\check{u}_0; \mathcal{M}(\check{u}_0))
=\frac{1}{N}\mathcal{J}_{soc}^{(N)}(u_0^*; \mathcal{M}(u_0^*))\\ &-\frac{1}{N}\mathcal{J}_{soc}^{(N)}(u_0^{(N)}; \mathcal{M}(u_0^{(N)}))%\quad (\leq O(\frac{1}{\sqrt{N}}))\\
+\frac{1}{N}\mathcal{J}_{soc}^{(N)}(u_0^{(N)}; \mathcal{M}(u_0^{(N)}))-\frac{1}{N}\mathcal{J}_{soc}^{(N)}(\check{u}_0; \mathcal{M}(\check{u}_0)).%\quad (\leq 0)\\
%\leq& O(\frac{1}{\sqrt{N}})+(\leq 0)\leq O(\frac{1}{\sqrt{N}})
\end{aligned}
%\right.
\end{equation*}
Note that $u_0^{(N)}$ is the centralized social optimal control in \eqref{u0*_origin}, thus one can easily obtain that
\begin{equation}\label{delta_1_2}
%\left\{
\begin{aligned}
\frac{1}{N}\mathcal{J}_{soc}^{(N)}(u_0^{(N)}; \mathcal{M}(u_0^{(N)}))\leq \frac{1}{N}\mathcal{J}_{soc}^{(N)}(\check{u}_0; \mathcal{M}(\check{u}_0)).
\end{aligned}
%\right.
\end{equation}
We know that $\mathcal{J}_{soc}^{(N)}(u_0; \mathcal{M}(u_0))$ continuously depends on $u_0$. %Thus we have:
Since $\mathcal{M}(u_0)$ is the solution of FBSDE \eqref{barxi_barpi} which continuously depends on parameters, we have $\mathcal{M}(u_0)$ is continuous in $u_0$. Note that $\mathcal{J}_{soc}^{(N)}(u_0; \mathcal{M}(u_0))$ is a quadratic functional and $u^*_0$ is fixed. %by following standard method in Huang \emph{et al} (2012), we obtain $|J_{soc}^{(N)}(u_0^{(N)}; \mathcal{M}(u_0^{(N)}))-J_{soc}^{(N)}(u_0^*; \mathcal{M}(u_0^*))|\leq N\cdot O(\frac{1}{\sqrt{N}})$.   %is $\varepsilon$-social optimality for $\check{u}_0$ parameterized process.
Let $\check{x}_0^{(N)}$ and $\check{x}_i^{(N)}$ %and $x_{(N)}^{(N)}$
be the state of the leader and the $i$th follower %and the population state of average
when the control of the leader is $u_0^{(N)}$. Denote
\begin{equation*}
\left\{
\begin{aligned}
&\acute{u}_0=u_0^{(N)}-u_0^*, \quad \delta\mathcal{M}(u_0)=\mathcal{M}(u_0^{(N)})-\mathcal{M}(u_0^*), \\
&\delta\mathcal{M}_i(u_0)=\mathcal{M}_i(u_0^{(N)})-\mathcal{M}_i(u_0^*),\quad \acute{x}_0=\check{x}_0^{(N)}-x_0^*, \quad \acute{x}_i=\check{x}_i^{(N)}-x_i^*.
\end{aligned}
\right.
\end{equation*}
Then we have
\begin{equation*}
\begin{aligned}
&\bigg|\mathcal{J}_{soc}^{(N)}(u_0^{(N)}; \mathcal{M}(u_0^{(N)}))-\mathcal{J}_{soc}^{(N)}(u_0^*; \mathcal{M}(u_0^*))\bigg|\\=&\bigg|\mathcal{J}_{soc}^{(N)}(u_0^{(N)}-u_0^*+u_0^*; \mathcal{M}(u_0^{(N)})-\mathcal{M}(u_0^*) +\mathcal{M}(u_0^*))-\mathcal{J}_{soc}^{(N)}(u_0^*; \mathcal{M}(u_0^*))\bigg|,
\end{aligned}
\end{equation*}
and
%\begin{equation*}
%\begin{aligned}
%&\mathcal{J}_{soc}^{(N)}(u_0^{(N)}; \mathcal{M}(u_0^{(N)}))=\alpha N\mathcal{J}_0(u_0^{(N)}; \mathcal{M}(u_0^{(N)}))+\sum_{i=1}^{N}\mathcal{J}_i(u_0^{(N)}; \mathcal{M}(u_0^{(N)}))\\
%=&\alpha N\mathcal{J}_0(u_0^*; \mathcal{M}(u_0^*))+\alpha NH'_0+\alpha NI'_0+\sum_{i=1}^{N}\mathcal{J}_i(u_0^*; \mathcal{M}(u_0^*))
%+\sum_{i=1}^{N}H'_i+\sum_{i=1}^{N}I'_i,
%\end{aligned}
%\end{equation*}
\begin{equation*}%\resizebox{\linewidth}{!}{$
\begin{aligned}
&\mathcal{J}_{soc}^{(N)}(u_0^{(N)}; \mathcal{M}(u_0^{(N)}))
%=\alpha N\mathcal{J}_0(u_0^{(N)}; \mathcal{M}(u_0^{(N)}))+\sum_{i=1}^{N}\mathcal{J}_i(u_0^{(N)}; \mathcal{M}(u_0^{(N)}))\\
=\alpha N[\mathcal{J}_0(u_0^*; \mathcal{M}(u_0^*))+H'_0+I'_0]+\sum_{i=1}^{N}[\mathcal{J}_i(u_0^*; \mathcal{M}(u_0^*))
+H'_i+I'_i],
\end{aligned}%$}
\end{equation*}
where
\begin{equation*}%\resizebox{\linewidth}{!}{$
\begin{aligned}
&\mathcal{J}_0(u_0^*; \mathcal{M}(u_0^*))=\mathbb{E}\bigg\{\int_{0}^{T}\|x_0^*- \Theta_0(x^*)^{(N)}-\eta_0\|_{Q_0}^{2}+\|u_0^*\|_{R_0}^{2}dt+ \|x_0^*(T)- \hat{\Theta}_0(x^*)^{(N)}(T)-\hat{\eta}_0\|_{G_0}^{2}\bigg\},\\
&H'_0=\mathbb{E}\bigg\{\int_{0}^{T}\|\acute{x}_0- \Theta_0\acute{x}^{(N)}\|_{Q_0}^{2}+\|\acute{u}_0\|_{R_0}^2dt+\|\acute{x}_0(T)- \hat{\Theta}_0\acute{x}^{(N)}(T)\|_{G_0}^{2}\bigg\},\\
&\mathcal{J}_i(u_0^*; \mathcal{M}(u_0^*))=\mathbb{E}\bigg\{\int_{0}^{T}\|x_i^*- \Theta (x^*)^{(N)}-\Theta_1x_0^*-\eta\|_{Q}^{2} +\|\mathcal{M}_i(u_0^*)\|_{R}^{2}dt\\
&\hspace{2.4cm}+\|x_i^*(T)- \hat{\Theta} (x^*)^{(N)}(T)-\hat{\Theta}_1x_0^*(T)-\hat{\eta}\|_{G}^{2}\bigg\},\\
&H'_i=\mathbb{E}\bigg\{\int_{0}^{T}\|\acute{x}_i- \Theta \acute{x}^{(N)}-\Theta_1\acute{x}_0\|_{Q}^{2} +\|\delta\mathcal{M}_i(u_0)\|_{R}^2dt+\|\acute{x}_i(T)- \hat{\Theta} \acute{x}^{(N)}(T)-\hat{\Theta}_1\acute{x}_0(T)\|_{G}^{2}\bigg\},\\
&I'_0=\mathbb{E}\bigg\{\int_{0}^{T}(x_0^*- \Theta_0(x^*)^{(N)}-\eta_0)^TQ_0(\acute{x}_0- \Theta_0\acute{x}^{(N)})dt\\
&\hspace{1cm}+(x_0^*(T)- \hat{\Theta}_0(x^*)^{(N)}(T)-\hat{\eta}_0)^TG_0(\acute{x}_0(T)- \hat{\Theta}_0\acute{x}^{(N)}(T))\bigg\},\\
%\end{aligned}
%\end{equation*}
%\begin{equation*}
%\begin{aligned}
&I'_i=\mathbb{E}\bigg\{\int_{0}^{T}(x^*_i- \Theta (x^*)^{(N)}-\Theta_1x^*_0-\eta)^TQ(\acute{x}_i- \Theta \acute{x}^{(N)}-\Theta_1\acute{x}_0)+\mathcal{M}_i^T(u_0^*)R\delta\mathcal{M}_i(u_0)dt\\
&\hspace{1cm}+(x^*_i(T)- \hat{\Theta} (x^*)^{(N)}(T)-\hat{\Theta}_1x^*_0(T)-\hat{\eta})^TG(\acute{x}_i(T)- \hat{\Theta} \acute{x}^{(N)}(T)-\hat{\Theta}_1\acute{x}_0(T))\bigg\}.
\end{aligned}%$}
\end{equation*}
By using similar arguments in Lemma \ref{lma_a.1} to Lemma \ref{lma_a.2} and $\|\Delta_2\|^2\leq c(\|\check{u}_0\|_{L^2}^2)\frac{1}{N}$, we obtain
\begin{equation*}
\begin{aligned}
&\frac{1}{N}H'_0+\frac{1}{N}H'_i+\alpha I'_0+\frac{1}{N}\sum_{i=1}^{N}I'_i%=\sqrt{O(\varepsilon_4^2+ \varepsilon_5^2+\varepsilon_6^2)\|u^*_0\|_{L^2}}
=O(\frac{1}{\sqrt{N}}).
\end{aligned}
\end{equation*}
%where
%\begin{equation*}
%\left\{
%\begin{aligned}
%&\varepsilon_4^2=\|\frac{1}{N}\sum_{i=1}^{N}x_{i(N)}-(x^*)^{(N)}\|_{L^2}^2\\
%&\varepsilon_5^2=\|\mathcal{M}_i(u_0^{(N)})-\mathcal{M}_i(u_0^*)\|_{L^2}^2\\
%&\varepsilon_6^2=\|x_{0(N)}-x_0^*\|_{L^2}^2.
%\end{aligned}
%\right.
%\end{equation*}
Hence, we have
\begin{equation}\label{delta_1_1}
%\left\{
\begin{aligned}
-\frac{1}{N}\mathcal{J}_{soc}^{(N)}(u_0^{(N)}; \mathcal{M}(u_0^{(N)}))+ \frac{1}{N}\mathcal{J}_{soc}^{(N)}(u_0^*; \mathcal{M}(u_0^*))\leq K(\frac{1}{\sqrt{N}})= O(\frac{1}{\sqrt{N}}),
\end{aligned}
%\right.
\end{equation}
where $K$ is independent of $N$.
By \eqref{delta_1_1} and \eqref{delta_1_2}, it follows that
\begin{equation*}
%\left\{
\begin{aligned}
\frac{1}{N}\mathcal{J}_{soc}^{(N)}(u_0^*; \mathcal{M}(u_0^*)) -\frac{1}{N}\mathcal{J}_{soc}^{(N)}(u_0^{(N)}; \mathcal{M}(u_0^{(N)}))= O(\frac{1}{\sqrt{N}}),
\end{aligned}
%\right.
\end{equation*}
and
\begin{equation*}
%\left\{
\begin{aligned}
\frac{1}{N}\mathcal{J}_{soc}^{(N)}(u_0^{(N)}; \mathcal{M}(u_0^{(N)}))-\frac{1}{N}\mathcal{J}_{soc}^{(N)}(\check{u}_0; \mathcal{M}(\check{u}_0))\leq 0,
\end{aligned}
%\right.
\end{equation*}
respectively. Thus, we have
\begin{equation}\label{delta_1}
%\left\{
\begin{aligned}
\Delta_1=&\frac{1}{N}\mathcal{J}_{soc}^{(N)}(u_0^*; \mathcal{M}(u_0^*))- \frac{1}{N}\mathcal{J}_{soc}^{(N)}(\check{u}_0; \mathcal{M}(\check{u}_0))\leq O(\frac{1}{\sqrt{N}}).
\end{aligned}
%\right.
\end{equation}
%Then, combining \eqref{delta_2}, \eqref{delta_1} and , we can obtain:
%\begin{equation*}
%%\left\{
%\begin{aligned}
%\Delta_1+\Delta_2\leq O(\frac{1}{\sqrt{N}})+c(\|\check{u}_0\|_{L^2})\frac{1}{\sqrt{N}}.
%\end{aligned}
%%\right.
%\end{equation*}
By Proposition \ref{prop_6.1}, there exists $K$ independent of $N$ such that $\|\check{u}_0\|_{L^2}\leq K$. Then, combining \eqref{delta_2}, \eqref{delta_1}, we can obtain:
%This implies
\begin{equation*}
%\left\{
\begin{aligned}
\Delta_1+\Delta_2\leq O(\frac{1}{\sqrt{N}})+c(\|\check{u}_0\|_{L^2})\frac{1}{\sqrt{N}})\leq K\cdot O(\frac{1}{\sqrt{N}})= O(\frac{1}{\sqrt{N}}),
\end{aligned}
%\right.
\end{equation*}
where $K$ is independent of $N$. The theorem follows. %Then the control $(u_0^*, u^*)$ is asymptotic social optimality.
%\qed
\end{proof}

\section{Numerical examples}

We now give a numerical example for Lemma \ref{lma_6.1}. By \eqref{ricattik} and \eqref{ricattikappa}, $\mathbb{K}$ and $\kappa$ can be easily computed. Consider $\mathbb{Y}=\mathbb{K}\mathbb{X}+\kappa$, we can obtain that
\begin{equation*}
%\left\{
\begin{aligned}
&d\mathbb{X}=[(\mathbb{A}+\mathbb{B}\mathbb{K})\mathbb{X} +\mathbb{B}\kappa+b]dt+\mathbb{D}dW_0, \quad
\mathbb{Y}=\mathbb{K}\mathbb{X}+\kappa,
\end{aligned}
%\right.
\end{equation*}
where $\mathbb{X}=((\hat{x})^T\ (\bar{x}_0)^T\ (q^*)^T\ (l_1^*)^T\ (l_2^*)^T)^T$, $\mathbb{Y}=((y^*)^T\ (\hat{y}^*)^T\ (y_0^*)^T\ (k_1)^T\ (k_2)^T)^T$.
%Similarly, by the following equations below \eqref{barxi_barpi}:
Since $p_i=\bar{P}\bar{x}_i+\bar{\varphi}$, by the following equations below \eqref{barxi_barpi}, we have
%\begin{equation*}
%\left\{
%\begin{aligned}
%&\dot{\bar{P}}+\bar{P}A-\bar{P}BR^{-1}B^T\bar{P}+A^T\bar{P}+Q=0, \ t\in [0,T],\
%\bar{P}(T)=G,\\
%%\end{aligned}
%%\right.
%%\end{equation*}
%%and
%%\begin{equation*}
%%\left\{
%%\begin{aligned}
%&\dot{\bar{\varphi}}+(A^T-\bar{P}BR^{-1}B^T)\bar{\varphi} +\chi_1+\bar{P}C\hat{x}+\bar{P}F\bar{x}_0=0, \ t\in [0,T],\
%%]dt-\zeta_0dW_0+(\bar{P}D-\zeta_i)dW_i=0,\\
%\bar{\varphi}(T)=\chi_2,
%\end{aligned}
%\right.
%\end{equation*}
%$\bar{P}$ and $\bar{\varphi}$ can be easily calculated. Since $p_i=\bar{P}\bar{x}_i+\bar{\varphi}$, we have
\begin{equation*}
%\left\{
\begin{aligned}
&d\bar{x}_i=[(A-BR^{-1}B^T\bar{P})\bar{x}_i-BR^{-1}B^T\bar{\varphi}+C\hat{x}+F\bar{x}_0]dt+DdW_i.
%p_i=\bar{P}\bar{x}_i+\bar{\varphi}.
\end{aligned}
%\right.
\end{equation*}
The realized decentralized state $x_0^*$ and $(x^*)^{(N)}$, can be derived by \eqref{realized state}.
Combining them with \eqref{ccsys2_1}, one can obtain
\begin{equation*}\resizebox{\linewidth}{!}{$
\left\{
\begin{aligned}
&d\left(
  \begin{array}{c}
    x_0^*-\bar{x}_0 \\
    (x^*)^{(N)}-\hat{x}  \\
  \end{array}
\right)=\bigg[\left(
  \begin{array}{cc}
    A_0 & C_0 \\
    F & A+C \\
  \end{array}
\right)\left(
  \begin{array}{c}
    x_0^*-\bar{x}_0  \\
    (x^*)^{(N)}-\hat{x}  \\
  \end{array}
\right)-\left(
  \begin{array}{c}
    0  \\
    BR^{-1}B^T  \\
  \end{array}
\right)(p^{(N)}-\hat{p})\bigg]dt
%&\hspace{3.5cm}
+\frac{1}{N}\left(
  \begin{array}{c}
    0 \\
    \sum_{1}^{N}D \\
  \end{array}\right)dW_i,\\
&\left(
  \begin{array}{c}
    x_0^*-\bar{x}_0  \\
    (x^*)^{(N)}-\hat{x}  \\
  \end{array}
\right)(0)=\left(
  \begin{array}{c}
    0  \\
    \frac{1}{N}\sum_{1}^{N}\xi_i-\hat{\xi}  \\
  \end{array}
\right),
\end{aligned}
\right.$}
\end{equation*}
where $\hat{p}=k_2$.

We continuously use the parameters in Example \ref{example_5.1}. The population $N=100$ and the time interval is $[0,12]$. By Matlab computation, the trajectories of the realized state $x^*_i$ is shown in Figure 1(a).
%\begin{figure}
%\centering
%%\begin{minipage}[c]{0.5\textwidth}
%\centering
%\includegraphics[height=4cm,width=6.8cm]{xstar1.jpg}
%%\end{minipage}%
%%\begin{minipage}[c]{0.5\textwidth}
%%\centering
%%\includegraphics[height=5cm,width=7.5cm]{pi1.jpg}
%%\end{minipage}
%\caption{The trajectories of $x^*_i$, $i=1,\cdots,100$.}
%\end{figure}

We defined
%$\varepsilon_1^2=\|(x^*)^{(N)}-\hat{x}\|_{L^2}^2$,
%$\varepsilon_2^2=\|p^{(N)}-\hat{p}\|_{L^2}^2$,
%$\varepsilon_3^2=\|x_0^*-\bar{x}_0\|_{L^2}^2$. Here, the integral is over the interval $[0,12]$ such that
$\varepsilon_1^2=\mathbb{E}\int_{0}^{12}\|(x*)^{(N)}-\hat{x}\|^2dt$, $\varepsilon_2^2=\mathbb{E}\int_{0}^{12}\|x^*_0-\bar{x}_0\|^2dt$, $\varepsilon_3^2=\mathbb{E}\int_{0}^{12}\|p^{(N)}-\hat{p}\|^2dt$. When $N$ increase from 1 to 100, the curves of $\varepsilon_1^2$, $\varepsilon_2^2$ and $\varepsilon_3^2$ are shown in Figure 1(b). The $X$ axis indicates $N$ and the $Y$ axis indicates $\varepsilon_i^2, i=1,2,3$.
%More specifically, $(a)$ is for $\varepsilon_1^2$ along with $N$, (b) is for $\varepsilon_2^2$ along with $N$ and $(c)$ is for $\varepsilon_3^2$ along with $N$.
It can be seen that they are approaching to zero when $N$ is growing larger and larger.
%\begin{figure}
%\centering
%\subfigure[]{
%\label{figa} %% label for first subfigure
%\includegraphics[width=1.4in]{XstarN-xhat2.jpg}}
%%\hspace{1in}
%\subfigure[]{
%\label{fig:subfig:b} %% label for second subfigure
%\includegraphics[width=1.4in]{Xstar0-barx02.jpg}}
%%\hspace{1in}
%\subfigure[]{
%\label{fig:subfig:c} %% label for second subfigure
%\includegraphics[width=1.4in]{pN-phat2.jpg}}
%\caption{The curves of $\varepsilon_1^2$, $\varepsilon_2^2$ and $\varepsilon_3^2$ when time interval is $[0,12]$.}
%%\label{figb} %% label for entire figure
%\end{figure}
\begin{figure}
\centering
\subfigure[]{
\label{figa} %% label for first subfigure
\includegraphics[height=4.6cm,width=7.3cm]{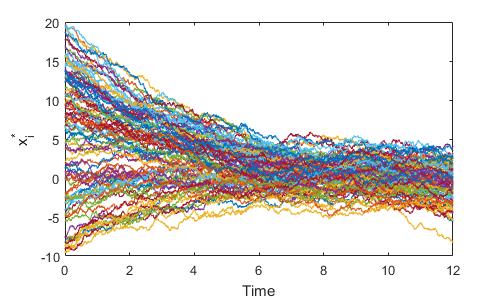}}
%\hspace{1in}
\subfigure[]{
\label{fig:subfig:b} %% label for second subfigure
\includegraphics[width=2.4in]{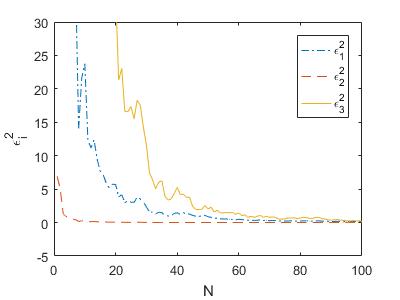}}
%\hspace{1in}
%\subfigure[]{
%\label{fig:subfig:c} %% label for second subfigure
%\includegraphics[width=1.4in]{pN-phat2.jpg}}
\caption{(a) is the trajectories of $x^*_i$, $i=1,\cdots,100$ and (b) is the curves of $\varepsilon_i^2$, $i=1,2,3$ when time interval is $[0,12]$.}
%\label{figb} %% label for entire figure
\end{figure}

\section{Conclusion}
This paper has analyzed the social optima in a class of LQ mean field control problem. We obtain the decentralized form of the optimal controls for the leader and $N$ followers. By Ricatti equation method, we discuss the solvability of the FBSDE. Finally, a Stackelberg equilibrium theorem is established. For future work, one can extend the results of this paper to the hierarchical control with many leaders case.

\appendix  %This command ends the counting of sections.
\section{Proof of Lemma \ref{lma_6.1}}
By \eqref{barxi_barpi} and \eqref{realized state}, we have
%\begin{equation*}
%\left\{
%\begin{aligned}
%&d(x^*)^{(N)}=[A(x^*)^{(N)}-BR^{-1}B^Tp^{(N)}+C(x^*)^{(N)} +Fx_0^*]dt+\frac{1}{N}\sum_{i=1}^{N}DdW_i,\\
%&d\bar{x}^{(N)}=[A\bar{x}^{(N)}-BR^{-1}B^Tp^{(N)}+C\hat{x}+F\bar{x}_0]dt +\frac{1}{N}\sum_{i=1}^{N}DdW_i, \\
%&(x^*)^{(N)}(0)=\frac{1}{N}\sum_{i=1}^{N}\xi_{i}, \quad \bar{x}^{(N)}(0)=\frac{1}{N}\sum_{i=1}^{N}\xi_{i}.
%\end{aligned}
%\right.
%  \end{equation*}
\begin{equation}\label{x*_barx_p_N}%\resizebox{\linewidth}{!}{$
\left\{
\begin{aligned}
&d(x^*)^{(N)}=[A(x^*)^{(N)}-BR^{-1}B^Tp^{(N)}+C(x^*)^{(N)} +Fx_0^*]dt+\frac{1}{N}\sum_{i=1}^{N}DdW_i,\ (x^*)^{(N)}(0)=\frac{1}{N}\sum_{i=1}^{N}\xi_{i},\\
&d\bar{x}^{(N)}=[A\bar{x}^{(N)}-BR^{-1}B^Tp^{(N)}+C\hat{x}+F\bar{x}_0]dt +\frac{1}{N}\sum_{i=1}^{N}DdW_i, \quad \bar{x}^{(N)}(0)=\frac{1}{N}\sum_{i=1}^{N}\xi_{i},\\
&dp^{(N)}=-[A^Tp^{(N)}+\chi_1]dt+\zeta_0dW_0 +\frac{1}{N}\sum_{i=1}^{N}\zeta_idW_i,
\quad p^{(N)}(T)=Gx^{(N)}(T)+\chi_2.\\
\end{aligned}
\right.
%$}
  \end{equation}
To prove Lemma \ref{lma_6.1}, we need the following two lemmas.
\begin{lmm}\label{lma_a.1}
Assume that \text{\rm(\textbf{A1})-(\textbf{A4})} hold. %and the CC systems \eqref{ccsys1_1} and \eqref{ccsys2_1} admit solutions.
Let $\bar{x}^{(N)}=\frac{1}{N}\sum_{i=1}^{N}\bar{x}_i$ and $p^{(N)}=\frac{1}{N}\sum_{i=1}^{N}p_i$. Then
%$\sup_{0\leq t\leq T}\mathbb{E}\|\bar{x}^{(N)}-\hat{x}\|^2=O(\frac{1}{N})$, $\sup_{0\leq t\leq T}\mathbb{E}\|p^{(N)}-\hat{p}\|^2=O(\frac{1}{N})$.
\begin{equation*}
\begin{aligned}
\sup_{0\leq t\leq T}\mathbb{E}\|\bar{x}^{(N)}-\hat{x}\|^2=O(\frac{1}{N}), \quad \sup_{0\leq t\leq T}\mathbb{E}\|p^{(N)}-\hat{p}\|^2=O(\frac{1}{N}).
\end{aligned}
  \end{equation*}
\end{lmm}
\begin{proof}
%{\it Proof }
%By \eqref{barxi_barpi}, we have
%\begin{equation*}
%\left\{
%\begin{aligned}
%&d\bar{x}^{(N)}=[A\bar{x}^{(N)}-BR^{-1}B^Tp^{(N)}+C\hat{x}+F\bar{x}_0]dt +\frac{1}{N}\sum_{i=1}^{N}DdW_i, \quad \bar{x}^{(N)}(0)=\frac{1}{N}\sum_{i=1}^{N}\xi_{i},\\
%&dp^{(N)}=-[A^Tp^{(N)}+\chi_1]dt+\zeta_0dW_0 +\frac{1}{N}\sum_{i=1}^{N}\zeta_idW_i,
%\quad p^{(N)}(T)=Gx^{(N)}(T)+\chi_2.\\
%\end{aligned}
%\right.
%  \end{equation*}
Combining \eqref{x*_barx_p_N} and \eqref{lim_x_p_fbsde_1}, we can obtain
\begin{equation*}
\left\{
\begin{aligned}
&d\mu_1=[A\mu_1-BR^{-1}B^T\mu_2]dt +\frac{1}{N}\sum_{i=1}^{N}DdW_i, \quad \mu_1(0)=\frac{1}{N}\sum_{i=1}^{N}\xi_{i}-\hat{\xi},\\
&d\mu_2=-[A^T\mu_2+Q\mu_1+\frac{1}{N}\sum_{i=1}^{N}\zeta_idW_i, \quad \mu_2(T)=G\mu_1,\\
\end{aligned}
\right.
  \end{equation*}
where $\mu_1=\bar{x}^{(N)}-\hat{x}$ and $\mu_2=p^{(N)}-\hat{p}$. Denote $\mu_2=P\mu_1+\varphi, \ t\in [0,T]$, where $P\in C^1([0,T];\mathbb{S}^{n})$ is the solution of the following Ricatti equation
and $\varphi\in C^1([0,T];\mathbb{R}^{n})$ satisfies
\begin{equation*}
\left\{
\begin{aligned}
&\dot{P}+PA-PBR^{-1}B^TP+A^TP+Q=0, \ t\in [0,T],\quad
P(T)=G,\\
%\end{aligned}
%\right.
%\end{equation*}
%and $\varphi\in C^1([0,T];\mathbb{R}^{n})$ satisfies
%\begin{equation*}
%\left\{
%\begin{aligned}
&d\varphi=-(A-BR^{-1}B^TP)^T\varphi dt+\frac{1}{N}\sum_{i=1}^{N}(PD-\zeta_i)dW_i, \ t\in [0,T],\quad
 \varphi(T)=0.
\end{aligned}
\right.
\end{equation*}
This is a standard Ricatti equation and the latter BSDE has a unique solution $\varphi=0, t\in[0,T]$. Thus  $\mu_2=P\mu_1$ and
\begin{equation*}
d\mu_1=[A-BR^{-1}B^TP]\mu_1dt +\frac{1}{N}\sum_{i=1}^{N}DdW_i.
\end{equation*}
By Cauchy-Schwarz inequality and Burkholder-Davis-Gundy's inequality, we have
\begin{equation*}
\begin{aligned}
&\sup_{0\leq t\leq T}\mathbb{E}\|\mu_1\|^2=\sup_{0\leq t\leq T}\mathbb{E}\bigg\|\int_{0}^{t}(A-BR^{-1}B^TP)\mu_1ds +\int_{0}^{t}\frac{1}{N}\sum_{i=1}^{N}DdW_i\bigg\|^2\\
\leq&2\sup_{0\leq t\leq T}\mathbb{E}\bigg\|\int_{0}^{t}(A-BR^{-1}B^TP)\mu_1ds\bigg\|^2 +2\sup_{0\leq t\leq T}\mathbb{E}\bigg\|\int_{0}^{t}\frac{1}{N}\sum_{i=1}^{N}DdW_i\bigg\|^2\\
\leq&2K\bigg\{\sup_{0\leq t\leq T}\mathbb{E}\int_{0}^{t}\|\mu_1\|^2 ds +\frac{1 }{N^2}\sum_{i=1}^{N}\mathbb{E}\int_{0}^{T}\|D\|^2 ds\bigg\}=2K\sup_{0\leq t\leq T}\mathbb{E}\int_{0}^{t}\|\mu_1\|^2 ds+O(\frac{1}{N}),
\end{aligned}
\end{equation*}
where constant $K$ is independent of $N$. Then, by Gronwall's inequality and $\mu_2=P\mu_1$, we obtain
\begin{equation*}
\sup_{0\leq t\leq T}\mathbb{E}\|\mu_1\|^2=O(\frac{1}{N}),\quad\sup_{0\leq t\leq T}\mathbb{E}\|\mu_2\|^2=O(\frac{1}{N}).
\end{equation*}
The lemma follows.
%\qed
\end{proof}
\begin{lmm}\label{lma_a.2}
Assume that \text{\rm(\textbf{A1})-(\textbf{A4})} hold. %and the CC systems \eqref{ccsys1_1} and \eqref{ccsys2_1} admit solutions.
Let $(x^*)^{(N)}=\frac{1}{N}\sum_{i=1}^{N}x^*_i$. Then
\begin{equation*}
\begin{aligned}
  &\sup_{0\leq t\leq T}\mathbb{E}\|x_0^*-\bar{x}_0\|^2=O(\frac{1}{N}), \quad \sup_{0\leq t\leq T}\mathbb{E}\|(x^*)^{(N)}-\bar{x}^{(N)}\|^2=O(\frac{1}{N}).
\end{aligned}
  \end{equation*}
\end{lmm}
\begin{proof}
%{\it Proof }
%By \eqref{barxi_barpi} and \eqref{realized state}, we have
%\begin{equation*}
%\left\{
%\begin{aligned}
%&d(x^*)^{(N)}=[A(x^*)^{(N)}-BR^{-1}B^Tp^{(N)}+C(x^*)^{(N)} +Fx_0^*]dt+\frac{1}{N}\sum_{i=1}^{N}DdW_i,\\
%&d\bar{x}^{(N)}=[A\bar{x}^{(N)}-BR^{-1}B^Tp^{(N)}+C\hat{x}+F\bar{x}_0]dt +\frac{1}{N}\sum_{i=1}^{N}DdW_i, \\
%&(x^*)^{(N)}(0)=\frac{1}{N}\sum_{i=1}^{N}\xi_{i}, \quad \bar{x}^{(N)}(0)=\frac{1}{N}\sum_{i=1}^{N}\xi_{i}.
%\end{aligned}
%\right.
%  \end{equation*}
Denote $\mu_3=x_0^*-\bar{x}_0$ and $\mu_4=(x^*)^{(N)}-\bar{x}^{(N)}$. By \eqref{x*_barx_p_N}, we can obtain
\begin{equation*}
%\left\{
\begin{aligned}
&d\left(
  \begin{array}{c}
    \mu_3  \\
    \mu_4  \\
  \end{array}
\right)=\bigg[\left(
  \begin{array}{cc}
    A_0 & C_0 \\
    F & A+C \\
  \end{array}
\right)\left(
  \begin{array}{c}
    \mu_3  \\
    \mu_4  \\
  \end{array}
\right)+\left(
  \begin{array}{c}
    C_0  \\
    C  \\
  \end{array}
\right)\mu_1\bigg]dt,\quad
\left(
  \begin{array}{c}
    \mu_3  \\
    \mu_4  \\
  \end{array}
\right)(0)=\left(
  \begin{array}{c}
    0  \\
    0  \\
  \end{array}
\right).
\end{aligned}
%\right.
\end{equation*}
For some constant $K$ which is independent of $N$ such that
\begin{equation*}
\begin{aligned}
&\sup_{0\leq t\leq T}\mathbb{E}\bigg\|\left(
  \begin{array}{c}
    \mu_3  \\
    \mu_4  \\
  \end{array}
\right)\bigg\|^2=\sup_{0\leq t\leq T}\mathbb{E}\bigg\|\int_{0}^{t}\bigg[\left(
  \begin{array}{cc}
    A_0 & C_0 \\
    F & A+C \\
  \end{array}
\right)\left(
  \begin{array}{c}
    \mu_3  \\
    \mu_4  \\
  \end{array}
\right)+\left(
  \begin{array}{c}
    C_0  \\
    C  \\
  \end{array}
\right)\mu_1\bigg]ds \bigg\|^2\\
%&\leq\sup_{0\leq t\leq T}\mathbb{E}\int_{0}^{t}\bigg\|\left(
%  \begin{array}{cc}
%    A_0 & C_0 \\
%    F & A+C \\
%  \end{array}
%\right)\left(
%  \begin{array}{c}
%    \mu_3  \\
%    \mu_4  \\
%  \end{array}
%\right)+\left(
%  \begin{array}{c}
%    C_0  \\
%    C  \\
%  \end{array}
%\right)\mu_1\bigg\|^2 ds \\
\leq&2K\bigg\{\sup_{0\leq t\leq T}\mathbb{E}\int_{0}^{t}\bigg\|\left(
  \begin{array}{c}
    \mu_3  \\
    \mu_4  \\
  \end{array}
\right)\bigg\|^2 ds +\sup_{0\leq t\leq T}\mathbb{E}\int_{0}^{t}\|\mu_1\|^2 ds\bigg\}=2K\sup_{0\leq t\leq T}\mathbb{E}\int_{0}^{t}\bigg\|\left(
  \begin{array}{c}
    \mu_3  \\
    \mu_4  \\
  \end{array}
\right)\bigg\|^2 ds+O(\frac{1}{N}).
\end{aligned}
\end{equation*}
By Gronwall's inequality, one can obtain
\begin{equation*}
\sup_{0\leq t\leq T}\mathbb{E}\bigg\|\left(
  \begin{array}{c}
    \mu_3  \\
    \mu_4  \\
  \end{array}
\right)\bigg\|^2=O(\frac{1}{N}).
\end{equation*}
Thus, the lemma follows.
%\qed\\~\\
\end{proof}
\begin{proof}[Proof of Lemma \ref{lma_6.1}] %\emph{of Lemma \ref{lma_6.1}}\quad
%{\it Proof of Lemma \ref{lma_6.1}. }
Since
\begin{equation*}
\begin{aligned}
\|(x^*)^{(N)}-\hat{x}\|^2=\|(x^*)^{(N)}-\bar{x}^{(N)}+\bar{x}^{(N)} -\hat{x}\|^2\leq 2\|(x^*)^{(N)}-\bar{x}^{(N)}\|^2+2\|\bar{x}^{(N)} -\hat{x}\|^2.
\end{aligned}
\end{equation*}
Combining Lemma \ref{lma_a.1} and Lemma \ref{lma_a.2}, it leads to
\begin{equation*}
\mathbb{E}\int_{0}^{T}\|(x^*)^{(N)}-\hat{x}\|^2dt+ \mathbb{E}\int_{0}^{T}\|p^{(N)}-\hat{p}\|^2dt+ \mathbb{E}\int_{0}^{T}\|x^*_0-\bar{x}_0\|^2dt\leq T\cdot O(\frac{1}{N})=O(\frac{1}{N}).
\end{equation*}
The lemma follows.
%\qed
\end{proof}
\section{Proof of Lemma \ref{lma_6.2}}
\begin{proof} %\textbf{}\quad
%{\it Proof }
By \eqref{ccsys2_1}, \eqref{x*_barx_p_N}, \eqref{realized state} and
%, we have:
%\begin{equation}\label{x0*_xN*}
%\left\{
%\begin{aligned}
%&dx_0^*=[A_0x_0^*-B_0(\alpha R_0)^{-1}B_0^Ty_0^*+C_0(x^*)^{(N)}]dt +D_0dW_0,\\
%&d(x^*)^{(N)}=[A(x^*)^{(N)}-BR^{-1}B^Tp^{(N)}+C(x^*)^{(N)} +Fx_0^*]dt+\frac{1}{N}\sum_{i=1}^{N}DdW_i,\\
%&x_0^*(0)=\xi_0,\quad (x^*)^{(N)}(0)=\frac{1}{N}\sum_{i=1}^{N}\xi_{i},\quad i=1,2,\cdots,N. \\
%\end{aligned}
%\right.
%\end{equation}
using a similar argument in Lemma \ref{lma_a.2}, one obtain that for some constant $K$ which is not dependent on $N$ such that
\begin{equation*}%\left\{
\begin{aligned}
%\mathbb{E}\sup_{0\leq t\leq T}\|x^*_i-\bar{x}_i\|^2&=\mathbb{E}\sup_{0\leq t\leq T}\bigg\|\int_{0}^{s}[A(x^*_i-\bar{x}_i)+C((x^*)^{(N)}-\bar{x})+F(x^*_0-\bar{x}_0)]ds \bigg\|^2\\
%&\leq K_5\mathbb{E}\sup_{0\leq t\leq T}\int_{0}^{s}\|x^*_i-\bar{x}_i\|^2+\|(x^*)^{(N)}-\bar{x}\|^2+\|x^*_0-\bar{x}_0\|^2 ds\\
%&=K_6\mathbb{E}\sup_{0\leq t\leq T}\int_{0}^{s}\|x^*_i-\bar{x}_i\|^2 ds+O(\frac{1}{N}).
&\sup_{0\leq t\leq T}\mathbb{E}\|x_0^*\|^2\leq K, \quad \sup_{0\leq t\leq T}\mathbb{E}\|(x^*)^{(N)}\|^2\leq K.
\end{aligned}%\right.
\end{equation*}
By Cauchy-Schwarz inequality and Burkholder-Davis-Gundy's inequality, we obtain
\begin{equation*}
\begin{aligned}
&\sup_{0\leq t\leq T}\mathbb{E}\|x^*_i\|^2=\sup_{0\leq t\leq T}\mathbb{E}\bigg\|\int_{0}^{t}(Ax^*_i-BR^{-1}B^Tp_i+C(x^*)^{(N)}+Fx^*_0)ds +DdW_i \bigg\|^2\\
\leq&\sup_{0\leq t\leq T}\mathbb{E}\int_{0}^{t}2\|Ax^*_i\|^2+2\|BR^{-1}B^Tp_i\|^2 +2\|C(x^*)^{(N)}\|^2+2\|Fx^*_0\|^2ds+\sup_{0\leq t\leq T}\int_{0}^{t}2\|DdW_i\|^2\\
\leq& 2\mathbb{E}\int_{0}^{T}\|Ax^*_i\|^2+\|BR^{-1}B^Tp_i\|^2+\|C(x^*)^{(N)}\|^2 +\|Fx^*_0\|^2ds+\int_{0}^{T}\|D\|^2ds\\
\leq&2A^2\sup_{0\leq t\leq T}\mathbb{E}\int_{0}^{T}\|x^*_i\|^2ds+K,
\end{aligned}
\end{equation*}
where constants $K$ is independent of $N$. By Gronwall's inequality, we have
\begin{equation*}
\sup_{1\leq i\leq N}\big[\sup_{0\leq t\leq T}\mathbb{E}\|x^*_i\|^2\big]\leq K,
\end{equation*}
where $K$ is not dependent on $N$. Then, according to Cauchy-Schwarz inequality, Burkholder-Davis-Gundy's inequality and the above discussion, we have
\begin{equation*}
\begin{aligned}
&\mathcal{J}_{soc}^{(N)}(u_0^*; u^*)=\alpha N\mathcal{J}_0(u_0^*;u^*)+\sum_{i=1}^{N}\mathcal{J}_i(u_0^*;u^*)\\
=&\alpha N\mathbb{E}\bigg\{\int_{0}^{T}\big[\|x_0^*- \Theta_0(x^*)^{(N)}-\eta_0\|_{Q_0}^{2}+\|-(\alpha R_0)^{-1}B_0^Ty_0^*\|_{R_0}^{2}\big]dt+ \|x_0^*(T)\\
&- \hat{\Theta}_0(x^*)^{(N)}(T)-\hat{\eta}_0\|_{G_0}^{2}\bigg\}+\sum_{i=1}^{N}\mathbb{E}\bigg\{\int_{0}^{T}\big[\|x_i^*- \Theta (x^*)^{(N)}(t)-\Theta_1x_0^*-\eta\|_{Q}^{2}\\
&+\|-R^{-1}B^Tp_i\|_{R}^{2}\big]dt+\|x_i^*(T)- \hat{\Theta} (x^*)^{(N)}(T)-\hat{\Theta}_1x_0^*(T)-\hat{\eta}\|_{G}^{2}\bigg\}\leq NK,\\
%&
\end{aligned}
\end{equation*}
where $K$ is independent of $N$. The lemma follows.
%\qed
\end{proof}
\section{Proof of Proposition \ref{prop_6.2} and \ref{prop_6.3}}
\begin{proof}[Proof of Proposition \ref{prop_6.2}] %\emph{of Proposition \ref{prop_6.2}}\quad
%{\it Proof of Proposition \ref{prop_6.2}. }
Since
\begin{equation*}
\left\{
\begin{aligned}
&d\delta x_0=(A_0\delta x_0+C_0\delta x^{(N)})dt,\quad \delta x_0(0)=0,\\
&d\delta x^{(N)}=[(A+C)\delta x^{(N)}+\frac{B}{N}\delta u_i]dt,\quad \delta x^{(N)}(0)=0,
\end{aligned}
\right.
\end{equation*}
and by Proposition \ref{prop_6.1},
%\begin{equation*}
%\alpha N\|u_0\|_{L^2}^2+\|u\|_{L^2}^2\leq \inf_{(u_0;u)}J_{soc}^{(N)}(u_0; u)\leq NK, (\alpha>0),
%\end{equation*}
%(this will be rigorously proved in Section 6).
we have $\|\delta u_i\|_{L^2}^2\leq K$, $K$ is independent of $N$. Using Cauchy-Schwarz inequality, it follows that
\begin{equation*}
    \begin{aligned}
\mathbb{E}\sup_{0\leq s\leq t}\|\delta x^{(N)}\|^2  =& \mathbb{E}\sup_{0\leq s\leq t}\Big\|\int_{0}^{s}[(A + C)\delta x^{( N)} + \frac{B}{N}\delta u_i  ]dr\Big\|^2
      \\\leq& K\mathbb{E}\int_{0}^{t}\|\delta x^{( N)}\|^2 dr  + \frac{1}{N^2}\|K\|^2dr
      \leq K\mathbb{E}\int_{0}^{t}\|\delta x^{( N)}\|^2 dr  + O\Big(\frac{1}{N^2}\Big),
    \end{aligned}
  \end{equation*}
where $K$ is independent of $N$. By Gronwall's inequality %$\mathbb{E}\sup_{0\leq t\leq T}\|\delta x^{(N)}\|^2 = O(\frac{1}{N^2})$.
  \begin{equation*}
    \begin{aligned}
      \mathbb{E}\sup_{0\leq t\leq T}\|\delta x^{(N)}\|^2 = O\Big(\frac{1}{N^2}\Big).
    \end{aligned}
  \end{equation*}
For $\delta x_0$, we have
\begin{equation*}
    \begin{aligned}
      \mathbb{E}\sup_{0\leq s\leq t}\|\delta x_0\|^2  = & \mathbb{E}\sup_{0\leq s\leq t}\Big\|\int_{0}^{s}[A_0\delta x_0+C_0\delta x^{( N)}]dr\Big\|^2
      \leq K\mathbb{E}\int_{0}^{t}\|\delta x_0\|^2 dr  + O\Big(\frac{1}{N^2}\Big),
    \end{aligned}
  \end{equation*}
where $K$ is independent of $N$. By Gronwall's inequality %$\mathbb{E}\sup_{0\leq t\leq T}\|\delta x_0\|^2 = O(\frac{1}{N^2})$.
  \begin{equation*}
    \begin{aligned}
      \mathbb{E}\sup_{0\leq t\leq T}\|\delta x_0\|^2 = O\Big(\frac{1}{N^2}\Big).
    \end{aligned}
  \end{equation*}
Moreover,
  \begin{equation*}
    \begin{aligned}
 \Theta^{T}Q(\bar{x}_i-\Theta \bar{x}^{(N)}-\Theta_1\bar{x}_0-\eta)\leq \frac{1}{N}\inf_{(u_0;u)}J_{soc}^{(N)}(u_0; u)\leq K.
   \end{aligned}
  \end{equation*}
%Similarly, $$\Theta_1^{T}Q(\bar{x}_i-\Theta\bar{x}^{(N)} -\Theta_1\bar{x}_0-\eta)\leq K,$$
%$$\hat{\Theta}^{T}G(\bar{x}_i(T) -\hat{\Theta}\bar{x}^{(N)}(T) -\hat{\Theta}_1\bar{x}_0(T)-\hat{\eta})\leq K,$$ $$\hat{\Theta}_1^{T}G(\bar{x}_i(T) -\hat{\Theta}\bar{x}^{(N)}(T) -\hat{\Theta}_1\bar{x}_0(T)-\hat{\eta})\leq K.$$
Similarly, $\Theta_1^{T}Q(\bar{x}_i-\Theta\bar{x}^{(N)} -\Theta_1\bar{x}_0-\eta)$,
$\hat{\Theta}^{T}G(\bar{x}_i(T) -\hat{\Theta}\bar{x}^{(N)}(T) -\hat{\Theta}_1\bar{x}_0(T)-\hat{\eta})$, $\hat{\Theta}_1^{T}G(\bar{x}_i(T) -\hat{\Theta}\bar{x}^{(N)}(T) -\hat{\Theta}_1\bar{x}_0(T)-\hat{\eta})$ are bounded.
Thus,
  \begin{equation*}
    \begin{aligned}
&\langle \Theta^{T}Q(\bar{x}_i-\Theta \bar{x}^{(N)}-\Theta_1\bar{x}_0-\eta),\delta x^{(N)}\rangle+ \langle\Theta_1^{T}Q(\bar{x}_i-\Theta\bar{x}^{(N)} -\Theta_1\bar{x}_0-\eta),\delta x_0\rangle+\langle \hat{\Theta}^{T}G(\bar{x}_i(T) -\hat{\Theta}\bar{x}^{(N)}(T)\\
& -\hat{\Theta}_1\bar{x}_0(T)-\hat{\eta}),\delta x^{(N)}(T)\rangle+\langle \hat{\Theta}_1^{T}G(\bar{x}_i(T) -\hat{\Theta}\bar{x}^{(N)}(T) -\hat{\Theta}_1\bar{x}_0(T)-\hat{\eta}),\delta x_0(T)\rangle=o(1).
   \end{aligned}
  \end{equation*}
The proposition follows.
%\qed\\~\\
\end{proof}
\begin{proof}[Proof of Proposition \ref{prop_6.3}] %\emph{of Proposition \ref{prop_6.3}}\quad
%{\it Proof of Proposition \ref{prop_6.3}. }
Since
\begin{equation*}
\left\{
\begin{aligned}
&d(\sum_{j\neq i}\delta x_j)=[A(\sum_{j\neq i}\delta x_j)+C\frac{N-1}{N}\delta x^{(N)} + F(N-1)\delta x_0]dt,\quad (\sum_{j\neq i}\delta x_j)(0)=0,\\
&d(N\delta x_0)=[A_0(N\delta x_0)+C_0(N\delta x^{(N)})]dt,\quad (N\delta x_0)(0)=0,\\
&d(N\delta x_j)=[A(N\delta x_j)+C(N\delta x^{(N)})+F(N\delta x_0)]dt,\quad (N\delta x_j)(0)=0.\\
\end{aligned}
\right.
\end{equation*}
According to equations in \eqref{dx_dagger}, one can obtain
\begin{equation*}
\left\{
\begin{aligned}
&d(N\delta x_j-\sum_{j\neq i}\delta x_j)=[A(N\delta x_j-\sum_{j\neq i}\delta x_j)+\frac{C}{N}\delta x^{( N)} + F\delta x_0]dt,\ (N\delta x_j-\sum_{j\neq i}\delta x_j)(0)=0,\\
&d(N\delta x_0-\delta x_0^{\dagger})=[A_0(N\delta x_0-\delta x_0^{\dagger})+C_0(N\delta x_j-\sum_{j\neq i}\delta x_j)+C_0(N\delta x_j-d\delta x^{\dagger})]dt,\ (N\delta x_0-\delta x_0^{\dagger})(0)=0,\\
&d(N\delta x_j-\delta x^{\dagger})=[(A+C)(N\delta x_j-\delta x^{\dagger})+C(N\delta x_j-\sum_{j\neq i}\delta x_j)+F(N\delta x_0-\delta x_0^{\dagger})]dt,\ (N\delta x_j-\delta x^{\dagger})(0)=0.\\
%& \quad \quad
\end{aligned}
\right.
\end{equation*}
Combing with the results in Proposition \ref{prop_6.2}, we have
\begin{equation*}
\begin{aligned}
\mathbb{E}\sup_{0\leq s\leq t}\|N\delta x_j-\sum_{j\neq i}\delta x_j\|^2  =&  \mathbb{E}\sup_{0\leq s\leq t}\Big\|\int_{0}^{s}[A(N\delta x_j-\sum_{j\neq i}\delta x_j) + \frac{C}{N}\delta x^{( N)} + F\delta x_0  ]dr\Big\|^2\\
% \\\leq& K\mathbb{E}\int_{0}^{t}\|\delta x^{( N)}\|^2 dr  + \frac{1}{N^2}\|c\|^2dr
\leq& K\mathbb{E}\int_{0}^{t}\|N\delta x_j-\sum_{j\neq i}\delta x_j\|^2 dr  + O\Big(\frac{1}{N^2}\Big),
\end{aligned}
\end{equation*}
where constant $K$ is independent of $N$. By Gronwall's inequality
\begin{equation*}
\begin{aligned}
\mathbb{E}\sup_{0\leq t\leq T}\|N\delta x_j-\sum_{j\neq i}\delta x_j\|^2 = O(\frac{1}{N^2}).
\end{aligned}
\end{equation*}
Since $N\delta x_0-\delta x_0^{\dagger}$ and $N\delta x_j-\delta x^{\dagger}$ are coupled, we have
\begin{equation*}
\begin{aligned}
&\mathbb{E}\sup_{0\leq s\leq t}\Bigg\|\left(
  \begin{array}{c}
    N\delta x_0-\delta x_0^{\dagger} \\
    N\delta x_j-d\delta x^{\dagger}  \\
  \end{array}
\right)\Bigg\|^2\\  = & \mathbb{E}\sup_{0\leq s\leq t}\Bigg\|\int_{0}^{s}\Bigg[\left(
  \begin{array}{cc}
    A_0 & C_0\\
    F & A+C \\
  \end{array}
\right)\left(
  \begin{array}{c}
    N\delta x_0-\delta x_0^{\dagger} \\
    N\delta x_j-d\delta x^{\dagger}  \\
  \end{array}
\right) + \left(
  \begin{array}{c}
    C_0 \\
    C \\
  \end{array}
\right)(N\delta x_j-\sum_{j\neq i}\delta x_j)\Bigg]dr\Bigg\|^2\\
% \\\leq& K\mathbb{E}\int_{0}^{t}\|\delta x^{( N)}\|^2 dr  + \frac{1}{N^2}\|c\|^2dr
\leq& K\mathbb{E}\int_{0}^{t}\Bigg\|\left(
  \begin{array}{c}
    N\delta x_0-\delta x_0^{\dagger} \\
    N\delta x_j-d\delta x^{\dagger}  \\
  \end{array}
\right)\Bigg\|^2 dr  + O\Big(\frac{1}{N^2}\Big),
\end{aligned}
\end{equation*}
where real-valued matrix $K$ is independent of $N$. By Gronwall's inequality $$\mathbb{E}\sup_{0\leq t\leq T}\Bigg\|\left(
  \begin{array}{c}
    N\delta x_0-\delta x_0^{\dagger} \\
    N\delta x_j-d\delta x^{\dagger}  \\
  \end{array}
\right)\Bigg\|^2 = O\Big(\frac{1}{N^2}\Big).$$
%Thus, $\|N\delta x_0-\delta x_0^{\dagger}\|$, $\|\sum_{j\neq i}\delta x_j-N\delta x_j\|$ and $\|N\delta x_j-\delta x^{\dagger}\|$ are $L^2$-convergence. %also the same as $\|x^{(N)} - \hat{x}\|$. The rigorous proof will be shown in section 6.
The proposition follows.
%\qed
\end{proof}
%
%\begin{acknowledgement}
%The authors wish to thank the anonymous referees for their valuable comments. The first author was supported by RGC Grants PolyU 153005/14P, 153275/16P. The second author was supported by the National Natural Science Foundation of China under Grant 61773241.
%\end{acknowledgement}

%\subsection{...}
%\subsection{...}
%...
%\subsection{...}
%...
%%-----------------------------
%%      your bibliography
%%-----------------------------

\begin{thebibliography}{99}%\vspace{-0.1cm}

\bibitem{bbs2010}
T. Ba\c{s}ar, A. Bensoussan and S. P. Sethi, Differential games with mixed leadership: the open-loop solution. \emph{Appl. Math. Comput.} \textbf{217} (2010) 972-979.%\vspace{-0.1cm}

\bibitem{bo1999}
T. Ba\c{s}ar and G. J. Olsder, \emph{Dynamic Noncooperative Game Theory}. SIAM, Philadelphia (1999). %\vspace{-0.1cm}

\bibitem{bp2014}
M. Bardi and F. S. Priuli, Linear-quadratic $N$-person and mean-field games with ergodic cost. \emph{SIAM J. Control Optim.} \textbf{52} (2014) 3022-3052.%\vspace{-0.1cm}

\bibitem{cd2004}
 C. T. Bauch and D. J. D. Earn, Vaccination and the theory of games. \emph{P. Natl. Acad. Sci. USA.} \textbf{101} (2004) 13391-13394.%\vspace{-0.1cm}

\bibitem{dr2012}
D. Bauso and R. Pesenti, Team theory and person-by-person optimization with binary decisions. \emph{SIAM J. Control Optim.} \textbf{50} (2012) 3011-3028.%\vspace{-0.1cm}

\bibitem{bcy2016a}
A. Bensoussan, M. H. M. Chau and S. C. P. Yam, Mean field games with a dominating player. \emph{Appl. Math. Optim.} \textbf{74} (2016) 91-128.%\vspace{-0.1cm}

\bibitem{bfy2013}
A. Bensoussan, J. Frehse and P. Yam, \emph{Mean Field Games and Mean Field Type Control Theory}. Springer, New York (2013). %\vspace{-0.1cm}

\bibitem{c2014}
P. E. Caines, Mean field games, in {\it Encyclopedia of Systems and Control}, Ed. T. Samad and J. Baillieul, Springer-Verlag, Berlin (2014).%\vspace{-0.1cm}

\bibitem{cd2013}
R. Carmona and F. Delarue, Probabilistic analysis of mean-field games. \emph{SIAM J. Control Optim.} \textbf{51} (2013) 2705-2734.%\vspace{-0.1cm}

\bibitem{cdl2016}
R. Carmona, F. Delarue and D. Lacker, Mean field games with common noise. \emph{Ann. Probab.} \textbf{44} (2016) 3740-3803.%\vspace{-0.1cm}

\bibitem{cw2016}
R. Carmona and P. Wang, Finite state mean field games with major and minor players. \emph{arXiv: 1610.05408.}%\vspace{-0.1cm}

\bibitem{cz2014}
R. Carmona and X. Zhu, A probabilistic approach to mean field games with major and minor players. \emph{Ann. Appl. Probab.} \textbf{26} (2016) 1535-1580%\vspace{-0.1cm}


\bibitem{cptd2012}
R. Couillet, S. M. Perlaza, H. Tembine and M. Debbah, Electrical vehicles in the smart grid: a mean field game analysis. \emph{IEEE. J. Sel. Area. Comm.} \textbf{30} (2012) 1086-1096.%\vspace{-0.1cm}

\bibitem{fc2016} D. Firoozi and P. E. Caines, Mean field game $\varepsilon$-Nash equilibria for partially observed optimal execution problems in finance. \emph{Proc. the IEEE 55th Conference on Decision and Control}. Las Vegas (2016) 268-275.%\vspace{-0.1cm}

\bibitem{gmm2012}
G. Gnecco, M. Sanguineti and M. Gaggero, Suboptimal solutions to team optimization problems with stochastic information structure. \emph{SIAM J. Optimiz.} \textbf{22} (2012) 212-243.%\vspace{-0.1cm}

\bibitem{tg1973} T. Groves, Incentives in teams. \emph{Econometrica}, \textbf{41} (1973) 617-631.%\vspace{-0.1cm}

\bibitem{hch1972} Y. C. Ho and K. C. Chu, Team decision theory and information structures in optimal control Part I. \emph{IEEE Trans. Automat. Contr.} \textbf{17} (1972) 15-22. %\vspace{-0.1cm}

\bibitem{jm2017}
J. Huang and M. Huang, Robust mean field linear-quadratic-Gaussian games with unknown $L^2$-disturbance. \emph{SIAM J. Control Optim.} \textbf{55} (2017) 2811-2840.%\vspace{-0.1cm}

\bibitem{hl2018}
J. Huang and N. Li, Linear quadratic mean-field game for stochastic delayed systems. \emph{IEEE Trans. Automat. Contr.} \textbf{63} (2018) 2722-2729.%\vspace{-0.1cm}

\bibitem{hww2016}
J. Huang, S. Wang and Z. Wu, Backward-forward linear-quadratic mean-field games with major and minor agents. \emph{Probability, Uncertainty and Quantitative Risk} \textbf{1} (2016) 1-27.%\vspace{-0.1cm}

\bibitem{hcm2003a} M. Huang, P. E. Caines and R. P. Malham\'{e}, Individual and mass behaviour in large population stochastic wireless power control problems: centralized and Nash equilibrium solutions. \emph{Proc. 42nd IEEE International Conference on Decision and Control}. Maui (2003) 98-103.%\vspace{-0.1cm}

%\bibitem{hcm2007a} M. Huang, P. E. Caines and R. P. Malham\'{e} (2007). Large-population cost-coupled LQG problems with non-uniform agents: individual-mass behavior and decentralized $\varepsilon$-Nash equilibria. \emph{IEEE Transactions on Automatic Control}, \textbf{52(9)}, 1560-1571.\vspace{-0.2cm}

%\bibitem{hmc2006} M. Huang, R. P. Malham\'{e} and P. E. Caines (2006). Large population stochastic dynamic games: closed-loop McKean-Vlasov systems and the Nash certainty equivalence principle. \emph{Communication in Information and Systems}, \textbf{6}, 221-251.\vspace{-0.2cm}

%\bibitem{h2010} M. Huang (2010). Large-population LQG games involving a major player: the Nash certainty equivalence principle. \emph{SIAM Journal on Control and Optimization}, \textbf{48(5)}, 3318-3353.\vspace{-0.2cm}

\bibitem{hcm2012} M. Huang, P. E. Caines and R. P. Malham\'{e}, Social optima in mean-field LQG control: centralized and decentralized strategies. \emph{IEEE Trans. Automat. Contr.} \textbf{57} (2012) 1736-1751.%\vspace{-0.1cm}

\bibitem{hn2016} M. Huang and S. L. Nguyen, Linear-quadratic mean field teams with a major agent. \emph{Proc. IEEE 55th Conference on Decision and Control}. Las Vegas (2016) 6958-6963.%\vspace{-0.1cm}

\bibitem{km2014}
A. C. Kizilkale and R. P. Malhame, Collective target tracking mean field control for Markovian jump-driven models of electric water heating loads. \emph{Proc. the 19th IFAC World Congress}. Cape Town, South Africa (2014) 1867-1972.%\vspace{-0.1cm}

\bibitem{ll2007} J. Lasry and P. Lions, Mean field games. \emph{Jpn. J. Math.} \textbf{2} (2007) 229-260.%\vspace{-0.1cm}

\bibitem{lz2008}
T. Li and J. Zhang, Asymptotically optimal decentralized control for large population stochastic multiagent systems. \emph{IEEE Trans. Automat. Contr.} \textbf{53} (2008) 1643-1660.%\vspace{-0.1cm}

\bibitem{ljz2019}
Y. N. Lin, X. S. Jiang, and W. H. Zhang, An open-loop Stackelberg strategy for the linear quadratic mean-field stochastic differential game. \emph{IEEE Trans. Automat. Contr.} \textbf{64} (2019) 97-110.%\vspace{-0.1cm}

\bibitem{my1999} J. Ma and J. Yong, \emph{Forward-Backward Stochastic Differential Equations and their Applications}.  Springer-Verlag, Berlin (1999).%\vspace{-0.1cm}

\bibitem{mzzgb2013}
S. Maharjan, Q. Zhu, Y. Zhang, S. Gjessing and T. Basar, Dependable demand response management in the smart grid: a Stackelberg game approach. \emph{IEEE Trans. Smart Grid.} \textbf{4} (2013) 120-132.%\vspace{-0.1cm}

\bibitem{jm1955}
J. Marschak, Elements for a theory of teams. \emph{Manage Sci.} \textbf{1} (1955) 127-137.%\vspace{-0.1cm}

\bibitem{mb2018}
J. Moon and T. Ba\c{s}ar, Linear quadratic mean field Stackelberg differential games. \emph{Automatica} \textbf{97} (2018) 200-213.%\vspace{-0.1cm}

\bibitem{nh2012}
S. L. Nguyen and M. Huang, Linear-quadratic-Gaussian mixed games with continuum-parametrized minor players. \emph{SIAM J. Control Optim.} \textbf{50} (2012) 2907-2937.%\vspace{-0.1cm}

\bibitem{nc2013}
M. Nourian and P. E. Caines, $\epsilon$-Nash mean field game theory for nonlinear stochastic dynamical systems with major and minor agents. \emph{SIAM J. Control Optim.} \textbf{51} (2013)  3302-3331.%\vspace{-0.1cm}

%\bibitem{ra1962}
%R. Radner (1962). Team decision problems. \emph{The Annals of Mathematical Statistics}, \textbf{33(3)}, 857-881.\vspace{-0.1cm}

\bibitem{swx2016}
J. Shi, G. Wang and J. Xiong, Leader-follower stochastic differential game with asymmetric information and applications. \emph{Automatica} \textbf{63} (2016) 60-73.

\bibitem{sc1973}
M. Simaan and J. Cruz, A Stackelberg solution for games with many players. \emph{IEEE Trans. Automat. Contr.} \textbf{18} (1973) 322-324.%\vspace{-0.1cm}

\bibitem{tzb2014}
H. Tembine, Q. Zhu and T. Ba\c{s}ar, Risk-sensitive mean-field games. \emph{IEEE Trans. Automat. Contr.} \textbf{59} (2014) 835-850.%\vspace{-0.1cm}

\bibitem{bj2012a}
B. Wang and J. Zhang, Mean field games for large population multiagent systems with Markov jump parameters.  \emph{SIAM J. Control Optim.} \textbf{50} (2012) 2308-2334.%\vspace{-0.1cm}

%\bibitem{bj2012b}
%B. Wang and J. Zhang (2012). Distributed control of multi-agent systems with random parameters and a major agent.  \emph{Automatica}, \textbf{48(9)}, 2093-2106.\vspace{-0.2cm}

\bibitem{bj2017a}
B. Wang and J. Zhang, Hierarchical mean field games for multiagent systems with tracking-type costs: distributed $\varepsilon$-Stackelberg equilibria.  \emph{IEEE Trans. Automat. Contr.} \textbf{59} (2014) 2241-2247.%\vspace{-0.1cm}

\bibitem{bj2017b}
B. Wang and J. Zhang, Social optima in mean field linear-quadratic-Gaussian models with Markov jump parameters.  \emph{SIAM J. Control Optim.} \textbf{55} (2017) 429-456.%\vspace{-0.1cm}

\bibitem{wh2017}
B. Wang and J. Huang, Social optima in robust mean field LQG control.  \emph{Proc. the 11th Asian Control
Conference}. Gold Coast, QLD (2017) 2089-2094.%\vspace{-0.1cm}



%\bibitem{ra1962}
%R. Radner (1962). Team decision problems. \emph{The Annals of Mathematical Statistics}, \textbf{33(3)}, 857-881.\vspace{-0.2cm}




\bibitem{wbvr2008}
G. Y. Weintraub, C. L. Benkard and B. V. Roy, Markov perfect industry dynamics with many firms. \emph{Econometrica} \textbf{76} (2008) 1375-1411.%\vspace{-0.1cm}

%\bibitem{lz1999}
%A. Lim and X. Y. Zhou (1999). Stochastic optimal LQR control with integral quadratic constraints and indefinite control weights. \emph{IEEE Transactions on Automatic Control}, \textbf{44(7)}, 1359-1369.\vspace{-0.2cm}

\bibitem{y1999} J. Yong,  Linear forward-backward stochastic differential equations. \emph{Appl. Math. Optim.} \textbf{39} (1999) 93-119.%\vspace{-0.1cm}

\bibitem{y02} J. Yong, A leader-follower stochastic linear quadratic differential game. \emph{SIAM J. Control Optim.} \textbf{41} (2002) 1015-1041.%\vspace{-0.1cm}

\bibitem{yz1999} J. Yong and X. Y. Zhou, \emph{Stochastic Controls: Hamiltonian Systems and HJB Equations.} Springer-Verlag. New York (1999). %\vspace{-0.1cm}


\end{thebibliography}
\end{document}